\documentclass[a4paper]{amsart}
\usepackage{amssymb} 
\usepackage[T1]{fontenc}
\usepackage[latin1]{inputenc}
\usepackage{amsfonts}
\usepackage{amsxtra}
\usepackage{ae}
\usepackage[all]{xy}
\usepackage{enumerate}

\include{diagram}

\newcommand*{\ket}{\rangle}
\newcommand*{\bra}{\langle}
\newcommand*{\Comp}{\mathfrak{Comp}}

\newcommand*{\D}{\mathcal{D}}
\newcommand*{\E}{\mathcal{E}}

\renewcommand*{\S}{\mathcal{S}}
\newcommand*{\cotimes}{\hat{\otimes}}

\newcommand*{\Fine}{\mathfrak{Fine}}
\newcommand*{\Bound}{\mathfrak{Bound}}

\DeclareMathOperator{\LSMod}{-\mathsf{Mod}}
\DeclareMathOperator{\SComodR}{\mathsf{Comod}-}

\DeclareMathOperator{\Hom}{Hom}

\DeclareMathOperator{\id}{id}

\newenvironment{bnum}
{\begin{list}{}
    {\setlength{\labelwidth}{15pt}
     \setlength{\leftmargin}{\labelwidth}
    }
}
{\end{list}}

\numberwithin{equation}{section}
\theoremstyle{change}
\newtheorem{theorem}{Theorem}[section]
\newtheorem{prop}[theorem]{Proposition}
\newtheorem{lemma}[theorem]{Lemma}

\newtheorem{definition}[theorem]{Definition}

\begin{document}

\title[Bornological quantum groups]{Bornological 
quantum groups}
\author{Christian Voigt}
\address{Institut for Mathematical Sciences\\
         University of Copenhagen\\
         Universitetsparken 5 \\
         2100 Copenhagen\\
         Denmark
}
\email{cvoigt@math.ku.dk}


\maketitle

\begin{abstract}
We introduce and study the concept of a bornological quantum group. This generalizes 
the theory of algebraic quantum groups in the sense of van Daele from the algebraic setting to 
the framework of bornological vector spaces. Working with bornological 
vector spaces, the scope of the latter theory can be extended considerably. In particular, 
the bornological theory covers smooth convolution algebras of arbitrary locally compact 
groups and their duals. Moreover Schwartz algebras of nilpotent Lie groups are bornological 
quantum groups in a natural way, and similarly one may consider algebras of functions on 
finitely generated discrete groups defined by various decay conditions. Another source of 
examples arises from deformation quantization in the sense of Rieffel. 
Apart from describing these examples we obtain some general results on 
bornological quantum groups. In particular, we construct the dual of a bornological quantum group and 
prove the Pontrjagin duality theorem. 
\end{abstract}

\section{Introduction}

The concept of a multiplier Hopf algebra introduced by van Daele \cite{vD1} extends the notion 
of a Hopf algebra to the setting of nonunital algebras. An important difference to the situation for ordinary Hopf 
algebras is that the comultiplication of a multiplier Hopf algebra $ H $ takes values in the multiplier algebra 
$ M(H \otimes H) $ and not in $ H \otimes H $ itself. Due to the occurence of multipliers, certain 
constructions with Hopf algebras have to be carried out more carefully in this context. 
Still, every multiplier Hopf algebra is equipped with a counit and 
an antipode satisfying analogues of the usual axioms. 
A basic example of a multiplier Hopf algebra is the algebra $ C_c(\Gamma) $ of compactly supported functions on a 
discrete group $ \Gamma $. This multiplier Hopf algebra is an ordinary Hopf algebra iff the 
group $ \Gamma $ is finite. \\
Algebraic quantum groups form a special class of multiplier Hopf algebras with particularly nice properties. 
Roughly speaking, an algebraic quantum group is a multiplier Hopf algebra with invertible antipode 
equipped with a Haar integral. Every 
algebraic quantum group admits a dual quantum group and the analogue of the Pontrjagin duality 
theorem holds \cite{vD2}. For instance, the multiplier Hopf algebra $ C_c(\Gamma) $ associated to a discrete group 
$ \Gamma $ is in fact an algebraic quantum group, its Pontrjagin dual being the complex group ring 
$ \mathbb{C}\Gamma $. More generally, all discrete and all compact quantum groups can be viewed as algebraic 
quantum groups. In addition, the class of algebraic 
quantum groups is closed under some natural operations 
including the construction of the Drinfeld double \cite{DvD}. 
Moreover, algebraic quantum groups give rise to examples of locally compact quantum 
groups \cite{KvD} illustrating nicely some general features of the latter. \\
However, due to the purely algebraic nature of the theory of multiplier Hopf algebras it is not possible to treat smooth convolution
algebras of Lie groups in this context, for instance. Roughly speaking, if viewed in terms of 
convolution algebras, the theory of algebraic quantum groups covers only the case of totally disconnected groups. 
Accordingly, the variety of quantum groups that can be described in this setting is obviously limited. 
It is thus desirable to have a more general setup then the one provided by algebraic quantum groups. \\
Motivated by these facts we introduce in this paper the concept of a bornological quantum group. 
The main idea is to replace the category of vector spaces underlying the definition of 
an algebraic quantum group by the category of bornological vector spaces. 
It is worth pointing out that bornological vector spaces provide the most natural setting to study various problems in noncommutative 
geometry and cyclic homology \cite{Meyerthesis}, \cite{Meyersmoothrep}, \cite{Meyercomb}, \cite{Voigtepch}, \cite{Voigtbs}. 
However, it seems that they have not been used in the context of quantum groups before. \\
The notion of a bornological vector space is dual to the concept of a locally convex vector space in a certain sense. 
Whereas the theory of locally convex vector spaces is based on the notion of an open subset, 
the key concept in the theory of bornological vector spaces is the notion of a bounded subset. 
It follows essentially from the definitions that both approaches are equivalent for normed spaces. 
In general, a locally convex vector space can be written canonically as a projective limit of normed spaces whereas a 
bornological vector space is the inductive limit of normed spaces in a canonical way. In particular both approaches yield natural extensions 
of the theory of normed spaces. It is an important fact that bornological and topological analysis 
are equivalent for Fr\'echet spaces \cite{Meyerborntop}. 
However, as soon as one moves beyond Fr\'echet spaces, the bornological approach
is superior to the topological one in many respects. In particular, the category of bornological vector spaces 
has much better algebraic properties than the category of locally convex vector spaces. 
It is thus in fact quite natural to work with bornological vector spaces in order to extend the theory of 
algebraic quantum groups. \\
Let us make some more specific comments on this paper. As far as the general theory of bornological quantum groups is concerned 
we follow the work by van Daele in the algebraic case. However, most of the constructions have to be rephrased in a more abstract way. 
Unfortunately, a large part of the arguments becomes less transparent at the same time. In addition, 
the calculations we have to perform are quite lengthy and annoying. 
On the other hand, with a basic knowledge of the theory of Hopf algebras it is easy to 
translate our considerations into a slightly imprecise but more appealing form using the Sweedler 
notation. \\
An important feature of the definition of a bornological quantum group is that it allows us to prove 
important general results while being quite simple at the same time. Although stronger assertions are 
possible at several points of the paper our setup seems to be sufficiently general for most purposes. \\
Let us now describe in detail how the paper is organized. For the convenience of the reader, we have 
collected some preliminary material on bornological vector spaces in section \ref{secborn}. Section 
\ref{secnondeg} contains the definition of
multiplier algebras of essential bornological algebras and a description of their basic properties. 
In section \ref{secqg} we define bornological quantum groups. 
Moreover we prove that a bornological quantum group can be characterized as a generalized 
Hopf algebra in the sense that there exists a counit and an antipode 
satisfying axioms familiar from the theory of Hopf algebras. 
In section \ref{secmodhaar} we study modular properties of the Haar 
integral of a bornological quantum group. This part of the theory is completely parallel 
to the algebraic case. 
Section \ref{secmodcomod} contains the definition of essential modules and comodules 
over a bornological quantum group and some related considerations. 
In section \ref{secdual} we construct the dual quantum group of a bornological quantum group 
and prove the Pontrjagin duality theorem. 
Using Pontrjagin duality we show in section \ref{secdualmod} that the category of essential modules over a 
bornological quantum group is isomorphic to the category of essential comodules over the dual quantum group 
and vice versa. As a consequence we obtain in particular a duality result for morphisms of bornological quantum groups. 
In the remaining part of the paper we consider some basic examples in order to illustrate the general theory. First we study 
the case of Lie groups in section \ref{seclie}. More precisely, we show that the smooth group algebra $ \D(G) $ 
of a Lie group $ G $ as well as the algebra of smooth functions $ C^\infty_c(G) $ 
with compact support are bornological quantum groups. 
Using the structure theory of locally compact groups, these considerations can actually be extended to arbitrary 
locally compact groups. In a similar way we study in section \ref{secschwartz} algebras of Schwartz functions on 
abelian and nilpotent Lie groups and their associated bornological quantum groups. 
Moreover, we consider algebras of functions on finitely generated discrete groups defined by various decay conditions 
with respect to the word length. These algebras can be viewed as certain completions of the complex group ring. 
Finally, in section \ref{secrieffel} we describe bornological quantum groups arising from 
deformation quantization in the sense of Rieffel \cite{Rieffelcqg}, \cite{Rieffelncqg}. \\
I would like to thank R. Meyer for some helpful comments. 

\section{Bornological vector spaces}\label{secborn}

In this section we review basic facts from the theory of bornological vector spaces. 
More information can be found in \cite{H-L1}, \cite{H-L2}, \cite{Meyerthesis}, \cite{Meyerborntop}. 
Throughout we work over the complex numbers. \\
A bornological vector space is a vector space $ V $ together with a collection $ \mathfrak{S}(V) $ of subsets 
of $ V $ satisfying certain conditions. These conditions can be viewed as an abstract reformulation of the properties of 
bounded subsets in a locally convex vector space. Following \cite{Meyerthesis}, we call a subset 
$ S $ of a bornological vector space $  V $ small iff it is contained in the bornology $ \mathfrak{S}(V) $. Throughout the 
paper we assume that all bornologies are convex. \\
As already indicated, the guiding example of a bornology is the collection of bounded subsets of a locally convex vector space. 
We write $ \Bound(V) $ for the bornological vector space associated to a locally convex vector space $ V $ in this way. 
One obtains another bornological vector space $ \Comp(V) $ by considering all precompact subsets 
of $ V $ instead. Recall that in a complete space $ V $ a set $ S \subset V $ is precompact 
iff its closure is compact. In certain situations the precompact bornology has nicer properties than the 
bounded bornology. Finally, one may view an arbitrary vector space $ V $ as a bornological vector space by considering 
the fine bornology $ \Fine(V) $. The fine bornology consists precisely of the 
bounded subsets of finite dimensional subspaces of $ V $. \\
Returning to the general theory, recall that a subset $ S $ of a complex vector space is called a disk if it is circled and convex. 
To a disk $ S \subset V $ one associates 
the semi-normed space $ \bra S \ket $ which is defined as the linear span of $ S $ endowed with the 
semi-norm $ \| \cdot \|_S $ given by the Minkowski functional. The disk $ S $ is called completant if $ \bra S \ket $ is 
a Banach space. A bornological vector space is called complete if each 
small set is contained in a completant small disk $ T $. Throughout the paper we work only with complete 
bornological vector spaces. For simplicity, we will not mention this explicitly in the sequel. \\
A linear map $ f: V \rightarrow W $ between bornological vector spaces is called bounded if it maps 
small subsets to small subsets. The space of bounded linear maps from $ V $ to $ W $ is 
denoted by $ \Hom(V,W) $. There is a natural bornology on $ \Hom(V,W) $ which consist of all sets $ L $ of linear maps such that 
$ L(S) \subset W $ is small for all small sets $ S \subset V $. In contrast, in the setting of locally convex vector 
spaces there are many different topologies on spaces of continuous linear maps. \\
We point out that the Hahn-Banach theorem does not hold for bornological vector spaces. 
In general there need not exist any nonzero bounded linear functional on a bornological vector space. 
A bornological vector space $ V $ is called 
regular if the bounded linear functionals on $ V $ separate points. The regularity of the underlying 
bornological vector space of a bornological quantum group will be guaranteed by the faithfulness of the Haar 
functional. Note also that all examples of bornological vector spaces arising from locally convex vector spaces 
are regular. \\
In the category of bornological vector spaces direct sums, direct 
products, projective limits and inductive limits exist. These constructions are characterized
by universal properties. Every complete bornological vector space can be written in a canonical way as a direct 
limit $ V \cong \varinjlim \bra S \ket $ of Banach spaces where the limit is taken over all 
small completant disks $ S $ in $ V $. In this way analysis in bornological vector spaces reduces
to analysis in Banach spaces. For instance, a sequence 
in a bornological vector space converges iff there exists a small disk $ S \subset V $ such that 
the sequence is contained in $ \bra S \ket $ and converges in this Banach space in the usual sense. \\
There exists a natural tensor product in the category of bornological vector spaces. More precisely, the 
bornological tensor product $ V \cotimes W $ is characterized by the universal property that 
bounded bilinear maps $ V \times W \rightarrow X $ correspond to bounded linear maps $ V \cotimes W \rightarrow X $. 
The bornological tensor product is associative and commutative and there is a natural adjunction isomorphism
$$
\Hom(V \cotimes W, X) \cong \Hom(V, \Hom(W,X)) 
$$
for all bornological vector spaces $ V, W, X $. 
This relation is one of the main reasons that the category of bornological vector spaces is much  
better adapted for algebraic constructions than the category of locally convex spaces. 
Note that the completed projective tensor product in the category of 
locally convex spaces does not have a right adjoint functor because it does not commute 
with direct sums. \\
Throughout the paper we will use the leg numbering convention for maps defined on tensor products. 
For instance, if $ f: V \cotimes W \rightarrow V \cotimes W $ 
is a bounded linear map and $ U $ is some bornological vector space we write $ f^{23} $ for the map 
$ \id \cotimes f $ defined on $ U \cotimes V \cotimes W $. Moreover, we 
sometimes write $ \id_{(n)} $ to indicate that we consider the identity map on an 
$ n $-fold tensor product. \\
A bornological algebra is a complete bornological 
vector space $ A $ with an associative multiplication given as a bounded linear map 
$ \mu: A \cotimes A \rightarrow A $. A homomorphism between 
bornological algebras is a bounded linear map $ f: A \rightarrow B $ 
which is compatible with multiplication. Remark that bornological algebras are not assumed to have a unit. 
The tensor product $ A \cotimes B $ of two bornological algebras is a bornological algebra in a 
natural way. \\
A (left) $ A $-module over a bornological algebra $ A $ is a bornological vector space $ M $ together with a bounded linear map 
$ \lambda: A \hat{\otimes} M \rightarrow M $ satisfying the 
axiom $ \lambda(\id \hat{\otimes} \lambda) = \lambda(\mu \hat{\otimes} \id) $ for 
an action. A homomorphism $ f: M \rightarrow N $ of $ A $-modules 
is a bounded linear map commuting with the action of $ A $. Such homomorphisms will also called 
$ A $-module maps. \\
Let us return to the basic examples of bornological vector spaces mentioned above. 
It follows immediately from the definitions that 
all linear maps $ f: V \rightarrow W $ from a fine space $ V $ into any 
bornological vector space $ W $ are bounded. In particular there is 
a fully faithful functor $ \mathfrak{Fine} $ from the category of complex vector 
spaces into the category of bornological vector spaces. 
This embedding is compatible with tensor products. If $ V_1 $ and 
$ V_2 $ are fine spaces the completed bornological tensor product 
$ V_1 \hat{\otimes} V_2 $ is the algebraic tensor product 
$ V_1 \otimes V_2 $ equipped with the fine bornology. 
In particular, every algebra $ A $ over the complex numbers can be viewed as a bornological 
algebra with the fine bornology. \\
In the case of Fr\'echet spaces a linear map $ f: V \rightarrow W $ is bounded for the bounded or 
the precompact bornologies iff it is continuous. 
Hence the functors $ \mathfrak{Bound} $ and $ \mathfrak{Comp} $ 
from the category of Fr\'echet spaces into the category of bornological vector spaces are fully faithful. 
The following theorem describes the completed bornological tensor product 
of Fr\'echet spaces with the precompact bornology and is proved in ~\cite{Meyerthesis}. 
\begin{theorem} \label{Frechettensor}
Let $ V $ and $ W $ be Fr\'echet spaces and let 
$ V \hat{\otimes}_\pi W $ be their completed projective tensor product. 
Then there is a natural isomorphism
\begin{align*}
&\mathfrak{Comp}(V) \hat{\otimes} \mathfrak{Comp}(W) \cong 
\Comp(V \hat{\otimes}_\pi W)
\end{align*}
of complete bornological vector spaces. 
\end{theorem}
In our setup for the theory of bornological quantum groups we will use the approximation property 
in order to avoid certain 
analytical problems with completed tensor products. The approximation property in the setting 
of bornological vector spaces has been studied by Meyer \cite{Meyerborntop}. 
Let us explain some definitions and results as far as they are needed for our purposes. \\
A subset $ S $ of a bornological vector space $ V $ is called compact if 
it is a compact subset of the Banach space $ \bra T \ket $ for some small disk $ T \subset V $. 
By definition, a bounded linear map $ f: V \rightarrow W $ can be approximated uniformly on compact subsets by 
finite rank operators if for all compact disks $ S \subset V $ there exists a sequence $ (f_n)_{n \in \mathbb{N}} $ of 
finite rank operators $ f_n: V \rightarrow W $ such that $ f_n $ converges to $ f $ in $ \Hom(\bra S \ket, W) $.  
An operator $ f: V \rightarrow W $ is of finite rank if it is contained in the image of 
the natural map from the uncompleted tensor product $ W \otimes V' $ into $ \Hom(V,W) $ 
where $ V' = \Hom(V, \mathbb{C}) $ is the dual space of $ V $. 
\begin{definition}\label{defapprox}
Let $ V $ be a complete bornological vector space. Then $ V $ has the approximation property if the identity 
map of $ V $ can be approximated uniformly on compact subsets by finite rank operators. 
\end{definition}
The following result from \cite{Meyerborntop} explains the relation between the bornological approximation 
property and Grothendieck's approximation property for locally convex vector spaces \cite{Grothendieck}. 
\begin{theorem}
Let $ V $ be a Fr\'echet space. Then $ \Comp(V) $ has the approximation property iff $ V $ has 
the approximation property as a locally convex vector space. 
\end{theorem}
We will need the following two properties of bornological vector spaces satisfying the approximation property. 
\begin{lemma}\label{appinj}
Let $ H $ be a bornological vector space satisfying the approximation property and let $ \iota: V \rightarrow W $ 
be an injective bounded linear map. Then the induced bounded linear map $ \id \cotimes \iota: H \cotimes V \rightarrow H \cotimes W $ 
is injective as well.
\end{lemma}
\proof Let $ x \in H \cotimes V $ be a nonzero element. We have to show that $ (\id \cotimes \iota)(x) $ is 
nonzero as well. The element $ x $ is contained in 
the completant disked hull $ R $ of $ S \otimes T $ for some small disks $ S \subset H $ and $ T \subset V $. 
Choose a sequence $ f_n: H \rightarrow H $ of finite rank operators approximating the identity uniformly 
on $ S $. Then $ f_n \cotimes \id $ approximates the identity uniformly on $ R $. 
Consider the commutative diagram 
$$
\xymatrix{
H \cotimes V\; \ar@{->}[r]^{\; \id \cotimes \iota} \ar@{->}[d]^{f_n \cotimes \id} & H \cotimes W \ar@{->}[d]^{f_n \cotimes \id} \\
H_n \otimes V \ar@{->}[r]^{\id \otimes \iota} & H_n \otimes W 
}
$$
where $ H_n $ denotes the image of $ f_n $. Observe that the spaces $ H_n \otimes V $ and $ H_n \otimes W $ are 
complete since $ H_n $ is finite dimensional. The bottom horizontal arrow 
is injective by assumption. Hence it suffices to show that $ (f_n \cotimes \id)(x) \in H_n \otimes V $ is nonzero for 
some $ n $. However, $ (f_n \cotimes \id)(x) $ converges to 
$ x $ in $ H \cotimes V $ since $ f_n \cotimes \id: H \cotimes V \rightarrow H \cotimes V $ 
converges uniformly to the identity on $ R $. This yields the claim. \qed 
\begin{lemma}\label{apphom}
Let $ H $ be a bornological vector space satisfying the approximation property and let $ V $ be 
an arbitrary bornological vector space. Then the canonical linear map $ \iota: H \cotimes V' \rightarrow \Hom(V, H) $ 
is injective.
\end{lemma}
\proof Let $ x \in H \cotimes V' $ be a nonzero element. We shall show that $ \iota(x) $ is nonzero 
as well. Again, $ x $ is contained in 
the completant disked hull $ R $ of $ S \otimes T $ for some small disks $ S \subset H $ and $ T \subset V' $. 
Choose a sequence $ f_n: H \rightarrow H $ of finite rank operators approximating the identity uniformly 
on $ S $ and consider the commutative diagram 
$$
\xymatrix{
H \cotimes V'\; \ar@{->}[r]^{\!\!\!\!\!\!\iota} \ar@{->}[d]^{f_n \cotimes \id} &\Hom(V,H) \ar@{->}[d]^{(f_n)_*} \\
H_n \cotimes V' \ar@{->}[r]^{\!\!\!\!\!\!\!\!\!\cong} & \Hom(V, H_n)
}
$$
where $ H_n $ is the image of $ f_n $ as above. In the same way as in the proof of lemma \ref{appinj} one 
obtains the assertion. \qed 

\section{Multiplier algebras}\label{secnondeg} 

In this section we prove basic results on multiplier algebras of bornological algebras that will be needed in the  
sequel. \\
The theory of multiplier Hopf algebras is an extension of the theory of Hopf algebras to the case where the 
underlying algebras do not have an identity element. Similarly, in our setting we will have to work with non-unital bornological 
algebras. However, in order to obtain a reasonable theory, it is necessary to impose some conditions on the multiplication. 
We will work with bornological algebras that are essential in the following sense. 
\begin{definition} 
A bornological algebra $ H $ is called essential if the multiplication map induces 
an isomorphism $ H \cotimes_H H \cong H $.
\end{definition}
In order to avoid trivialities we shall always assume that essential bornological algebras are different from zero. 
Clearly, every unital bornological algebra is essential. If $ H $ has an approximate identity \cite{Meyersmoothrep} 
then $ H $ is essential iff the multiplication $ H \cotimes H \rightarrow H $ is a bornological quotient 
map. A bornological algebra $ H $ has an approximate identity if for every small subset $ S $ of $ H $ 
there is a sequence $ (u_n)_{n \in \mathbb{N}} $ in $ H $ such that $ u_n x $ and $ x u_n $ converge uniformly 
to $ x $ for every $ x \in S $. 
We will not require the existence of approximate identites in the general definition of a bornological quantum group. 
However, in many cases it is easy to check that approximate identities do indeed exist. 
\begin{definition} 
Let $ H $ be a bornological algebra. An $ H $-module $ V $ 
is called essential if the canonical map $ H \cotimes_H V \rightarrow V $ is an isomorphism. 
\end{definition}
An analogous definition can be given for right modules. In particular, an essential algebra $ H $ is an essential
left and right module over itself. \\
We shall now discuss multipliers. A left multiplier for a bornological algebra $ H $ is a bounded linear map 
$ L: H \rightarrow H $ such that $ L(fg) = L(f)g $ for all $ f,g \in H $. Similarly, a right multiplier 
is a bounded linear map $ R: H \rightarrow H $ such that $ R(fg) = fR(g) $ for all $ f,g \in H $. 
We let $ M_l(H) $ and $ M_r(H) $ be the spaces of left and right multipliers, respectively. 
These spaces are equipped with the subspace bornology of $ \Hom(H,H) $ and become bornological algebras 
with multiplication given by composition of maps. The multiplier algebra $ M(H) $ of a bornological algebra $ H $ is the space of all 
pairs $ (L,R) $ where $ L $ is a left muliplier and $ R $ is a right multiplier 
for $ H $ such that $ fL(g) = R(f) g $ for all $ f,g \in H $. The bornology and algebra structure of $ M(H) $ are 
inherited from $ M_l(H) \oplus M_r(H) $. 
There is a natural homomorphism $ \iota: H \rightarrow M(H) $. By construction, $ H $ is a left and right $ M(H) $-module 
in a natural way. \\
Let $ H $ and $ K $ be bornological algebras and let $ f: H \rightarrow M(K) $ be a 
homomorphism. Then $ K $ is a left and right $ H $-module in an obvious way. 
We say that the homomorphism $ f: H \rightarrow M(K) $ is essential if it turns 
$ K $ into an essential left and right $ H $-module. That is, for the corresponding module structures we have 
$ H \cotimes_H K \cong K \cong K \cotimes_H H $ in this case. Note that the identity map $ \id: H \rightarrow H $ 
defines an essential homomorphism $ H \rightarrow M(H) $ iff the bornological algebra $ H $ is essential. 
\begin{lemma}
Let $ H $ be a bornological algebra and 
let $ f: H \rightarrow M(K) $ be an essential homomorphism into the multiplier algebra of an essential bornological algebra $ K $. Then there 
exists a unique unital homomorphism $ F: M(H) \rightarrow M(K) $ such that $ F \iota = f $ where 
$ \iota: H \rightarrow M(H) $ is the canonical map.  
\end{lemma}
\proof We obtain a bounded linear map $ F_l: M_l(H) \rightarrow M_l(K) $ by 
\begin{equation*} 
M_l(H) \cotimes K \cong
 \xymatrix{
     M_l(H) \cotimes H \cotimes_H K \; \ar@{->}[r]^{\;\;\;\;\;\;\;\; \mu \cotimes \id }
      & H \cotimes_H K }
\cong K 
\end{equation*}
and accordingly a map $ F_r: M_r(H) \rightarrow M_r(K) $ by 
\begin{equation*} 
K \cotimes M_r(H) \cong
 \xymatrix{
     K \cotimes_H H \cotimes M_r(H) \; \ar@{->}[r]^{\;\;\;\;\;\;\;\; \id \cotimes \mu }
      & K \cotimes_H H }
\cong K.
\end{equation*}
It is straightforward to check that $ F((L,R)) = (F_l(L),F_r(R)) $ defines a unital homomorphism 
$ F: M(H) \rightarrow M(K) $ such that $ F \iota = f $. Uniqueness of $ F $ follows from the 
fact that $ f(H) \cdot K \subset K $ and $ K \cdot f(H) \subset K $ are dense subspaces. \qed 
\begin{lemma}\label{lemnondeg3}
Let $ H_1, H_2 $ be essential bornological algebras and let 
$ f_1: H_1 \rightarrow M(K_1) $ and $ f_2: H_2 \rightarrow M(K_2) $ be essential 
homomorphisms into the multiplier algebras of bornological algebras $ K_1 $ and $ K_2 $. Then the induced homomorphism 
$ f_1 \cotimes f_2 : H_1 \cotimes H_2 \rightarrow M(K_1 \cotimes K_2) $ is essential. 
\end{lemma}
\proof By assumption the maps $ f_1 $ and $ f_2 $ induce isomorphisms $ H_1 \cotimes_{H_1} K_1 \cong K_1 $ and 
$ H_2 \cotimes_{H_2} K_2 \cong K_2 $. Let us show that the natural bounded 
linear map 
$$
\beta: (H_1 \cotimes H_2) \cotimes_{(H_1 \cotimes H_2)} (K_1 \cotimes K_2) \rightarrow H_1 \cotimes_{H_1} K_1 \cotimes H_2 \cotimes_{H_2} K_2
$$
is an isomorphism. We observe that
\begin{align*}
h_1 l_1 x \otimes h_2 l_2 \otimes &k_1  \otimes k_2 - h_1 l_1 \otimes h_2 l_2 \otimes x k_1 \otimes k_2 \\
&= h_1 l_1 x \otimes h_2 l_2 \otimes k_1  \otimes k_2 - h_1 \otimes h_2 \otimes l_1 x k_1 \otimes l_2 k_2 \\
& \qquad + h_1 \otimes h_2 \otimes l_1 x k_1 \otimes l_2 k_2 - h_1 l_1 \otimes h_2 l_2 \otimes x k_1 \otimes k_2
\end{align*}
is zero in $ (H_1 \cotimes H_2) \cotimes_{(H_1 \cotimes H_2)} (K_1 \cotimes K_2) $.  
Using a similar formula for tensor relations over $ H_2 $ we see that the map 
$$
H_1 \otimes H_1 \otimes K_1 \otimes H_2 \otimes H_2 \otimes K_2 \rightarrow (H_1 \cotimes H_2) \cotimes_{(H_1 \cotimes H_2)} K_1 \cotimes K_2
$$
induced by multiplication in $ H_1 $ and $ H_2 $ and a flip of the tensor factors descends to a bounded linear map 
$$
H_1 \cotimes_{H_1} K_1 \cotimes H_2 \cotimes_{H_2} K_2 
\rightarrow (H_1 \cotimes H_2) \cotimes_{(H_1 \cotimes H_2)} (K_1 \cotimes K_2)
$$
which is inverse to the map $ \beta $. Hence the map $ f_1 \cotimes f_2 $ induces an isomorphism 
$ (H_1 \cotimes H_2) \cotimes_{H_1 \cotimes H_2} (K_1 \cotimes K_2) \cong K_1 \cotimes K_2 $. In a similar way one obtains the 
isomorphism $ (K_1 \cotimes K_2) \cotimes_{H_1 \cotimes H_2}(H_1 \cotimes H_2) \cong K_1 \cotimes K_2 $. \qed \\
Following the terminology of van Daele \cite{vD1}, we say that a bornological algebra $ H $ is nondegenerate if 
$ fg = 0 $ for all $ g \in H $ implies $ f = 0 $ and $ fg = 0 $ for all $ f $ implies $ g = 0 $. 
These conditions can be reformulated by saying that the natural maps 
$$
H \rightarrow M_l(H), \quad H \rightarrow M_r(H)
$$
are injective. In particular, for a nondegenerate bornological algebra the canonical map $ H \rightarrow M(H) $ is injective. \\
Nondegeneracy of a bornological algebra is a consequence of the existence of a faithful linear functional in
the following sense. 
\begin{definition} 
Let $ H $ be a bornological algebra. A bounded linear functional $ \omega: H \rightarrow \mathbb{C} $ is called faithful if
$ \omega(fg) = 0 $ for all $ g $ implies $ f = 0 $ and $ \omega(fg) = 0 $ for all $ f $ implies $ g = 0 $. 
\end{definition}
Remark that a bornological algebra $ H $ equipped with a faithful bounded linear functional is regular in the sense that 
bounded linear functionals separate the points of $ H $. 
\begin{lemma}\label{appmult}
Let $ H_1 $ and $ H_2 $ be bornological algebras satisfying the approximation property equipped with faithful bounded linear functionals 
$ \phi_1 $ and $ \phi_2 $, respectively. Then $ \phi_1 \cotimes \phi_2 $ is a faithful linear functional on $ H_1 \cotimes H_2 $. 
\end{lemma}
\proof Let us show that the canonical map $ H_1 \cotimes H_2 \rightarrow M_l(H_1 \cotimes H_2) $ is injective. 
We assume that $ x $ is in the kernel of this map. Since $ \phi_1 $ is faithful the map $ \mathcal{F}: H_1 \rightarrow H_1' $ given by 
$ \mathcal{F}(f)(g) = \phi_1(fg) $ is injective. Consider the chain of maps 
\begin{equation*} 
\xymatrix{
H_1 \cotimes H_2 \; \ar@{->}[r]^{\mathcal{F} \cotimes \id} & H_1' \cotimes H_2 \ar@{->}[r] & \Hom(H_1, H_2)
}  
\end{equation*}
where the second arrow is the obvious one. The first of these maps is injective according to lemma \ref{appinj}, the second map is 
injective according to lemma \ref{apphom}. 
Since $ x(g \otimes h) = 0 $ for all $ g \in H_1 $ and $ h \in H_2 $ and $ \phi_2 $ is faithful we and see that the image of $ x $ in 
$ \Hom(H_1, H_2) $ is zero. Hence $ x = 0 $ as well. The assertion concerning right multipliers 
is proved in a similar way. \qed 

\section{Bornological quantum groups}\label{secqg} 

In this section we introduce the notion of a bornological quantum group. 
Moreover we prove that every bornological quantum group is equipped with a counit and an invertible antipode. \\
In the sequel we assume that $ H $ is an 
essential bornological algebra satisfying the approximation property. Moreover we suppose that $ H $ is equipped with a faithful bounded 
linear functional. Remark that we may thus view $ H $ as a subset of the multiplier algebra $ M(H) $. We will do this 
frequently without further mentioning. Taking into account 
lemma \ref{appmult}, an analogous statement applies to tensor powers of $ H $. \\
First we have to discuss the concept of a comultiplication on $ H $. Let $ \Delta: H \rightarrow M(H \cotimes H) $ be a homomorphism. 
The left Galois maps $ \gamma_l, \gamma_r: H \cotimes H \rightarrow M(H \cotimes H) $ for $ \Delta $ are defined by 
$$
\gamma_l(f\otimes g) = \Delta(f)(g \otimes 1),\qquad \gamma_r(f \otimes g) = \Delta(f) (1 \otimes g). 
$$
Similarly, the right Galois maps $ \rho_l, \rho_r: H \cotimes H \rightarrow M(H \cotimes H) $ for $ \Delta $ are defined by 
$$
\rho_l(f \otimes g) = (f \otimes 1) \Delta(g), \qquad \rho_r(f \otimes g) = (1 \otimes f) \Delta(g). 
$$
These maps, or rather their appropriate analogues, play an important role in the algebraic as well as the analytic 
theory of quantum groups \cite{vD1}, \cite{vD2}, \cite{BaSka2}.
Our terminology is motivated from the fact that they also occur in the theory of Hopf-Galois 
extensions \cite{Montgomery}. \\
Assume in addition that the homomorphism $ \Delta: H \rightarrow M(H \cotimes H) $ is essential. Then $ \Delta $ is called coassociative 
if
$$
(\Delta \cotimes \id) \Delta = (\id \cotimes \Delta)\Delta 
$$
where both sides are viewed as maps from $ H $ to $ M(H \cotimes H \cotimes H) $. Remark that these maps are well-defined according to 
lemma \ref{lemnondeg3}. 
\begin{definition}
An essential homomorphism $ \Delta: H \rightarrow M(H \cotimes H) $ is 
called a comultiplication if it is coassociative. 
\end{definition}
An essential algebra homomorphism $ f: H \rightarrow M(K) $ between bornological algebras equipped with comultiplications is called a 
coalgebra homomorphism if $ \Delta f = (f \cotimes f)\Delta $. \\
We need some more terminology. The opposite algebra $ H^{op} $ of $ H $ is the space $ H $ equipped with 
the opposite multiplication. That is, the multiplication $ \mu^{op} $ in $ H^{op} $ is defined by $ \mu^{op} = \mu \tau $ where 
$ \mu: H \cotimes H \rightarrow H $ is the multiplication in $ H $ and $ \tau: H \cotimes H \rightarrow H \cotimes H $ 
is the flip map given by $ \tau(f \otimes g) = g \otimes f $. An algebra antihomomorphism between $ H $ and $ K $ 
is an algebra homomorphism $ \phi: H \rightarrow K^{op} $. Equivalently, an algebra antihomomorphism can be viewed as 
an algebra homomorphism $ H^{op} \rightarrow K $. 
If $ \Delta: H \rightarrow M(H \cotimes H) $ is a comultiplication then $ \Delta $ also defines a comultiplication 
$ H^{op} \rightarrow M(H^{op} \cotimes H^{op}) $. We write $ \gamma_l^{op}, \gamma_r^{op}, \rho_l^{op}, \rho_r^{op} $ for the 
corresponding Galois maps. \\
Apart from changing the order of multiplication we may also reverse the order of a comultiplication. 
If $ \Delta: H \rightarrow M(H \cotimes H) $ is a comultiplication then the opposite comultiplication $ \Delta^{cop} $ is the essential homomorphism
from $ H $ to $ M(H \cotimes H) $ defined by $ \Delta^{cop} = \tau \Delta $. We write $ \gamma_l^{cop}, \gamma_r^{cop}, \rho_l^{cop}, \rho_r^{cop} $ for the 
Galois maps associated to this comultiplication. Moreover we write $ H^{cop} $ for $ H $ equipped with the opposite comultiplication. 
Using opposite comultiplications we obtain the notion of a coalgebra antihomomorphism. \\
We may also combine these procedures, that is, reverse both multiplication and comultiplication. 
The bornological algebra with comultiplication arising in this way is denoted by $ H^{op cop} = (H^{op})^{cop} $ and we write 
$ \gamma_l^{op cop}, \gamma_r^{op cop}, \rho_l^{op cop}, \rho_r^{op cop} $ for the corresponding Galois maps. \\
It is straightforward to check that the Galois maps of $ H, H^{op}, H^{cop} $ and $ H^{opcop} $ are related as follows. 
\begin{lemma} \label{Galopcop}
Let $ \Delta: H \rightarrow M(H \cotimes H) $ be a comultiplication. Then we have
\begin{align*}
\gamma_r &= \tau \gamma_l^{cop}, \quad \rho_l = \gamma_l^{op} \tau, \quad \rho_r = \tau \gamma_l^{opcop} \tau \\
\gamma_l &= \tau \gamma_r^{cop}, \quad \rho_r = \gamma_r^{op} \tau, \quad \rho_l = \tau \gamma_r^{opcop} \tau \\
\rho_r &= \tau \rho_l^{cop}, \quad \gamma_l = \rho_l^{op} \tau, \quad \gamma_r = \tau \rho_l^{opcop} \tau \\
\rho_l &= \tau \rho_r^{cop}, \quad \gamma_r = \rho_r^{op} \tau, \quad \gamma_l = \tau \rho_r^{opcop} \tau
\end{align*}
for the Galois maps of $ H, H^{op}, H^{cop} $ and $ H^{opcop} $. 
\end{lemma}
Observe that the previous relations can also be rewritten in the form
\begin{align*}
\gamma_l^{cop} &= \tau \gamma_r, \quad \gamma_l^{op} = \rho_l \tau, \quad \gamma_l^{opcop} = \tau \rho_r \tau \\
\gamma_r^{cop} &= \tau \gamma_l, \quad \gamma_r^{op} = \rho_r \tau, \quad \gamma_r^{opcop} = \tau \rho_l \tau \\
\rho_l^{cop} &= \tau \rho_r, \quad \rho_l^{op} = \gamma_l \tau, \quad \rho_l^{opcop} = \tau \gamma_r \tau \\
\rho_r^{cop} &= \tau \rho_l, \quad \rho_r^{op} = \gamma_r \tau, \quad \rho_r^{opcop} = \tau \gamma_l \tau.  
\end{align*}
As a consequence, the Galois maps for $ H $ may be expressed in terms of the maps $ \gamma_l, \gamma_l^{op}, \gamma_l^{cop} $ and 
$ \gamma_l^{opcop} $ and vice versa. Of course, there are similar statements for $ \gamma_r, \rho_l $ and $ \rho_r $. 
This basic observation will be used frequently below. \\
Let $ \Delta: H \rightarrow M(H \cotimes H) $ be a comultiplication such that all Galois maps associated 
to $ \Delta $ define bounded linear maps from $ H \cotimes H $ into itself. If $ \omega $ is a bounded linear 
functional on $ H $ we define for every $ f \in H $ a multiplier $ (\id \cotimes \omega)\Delta(f) \in M(H) $ by 
\begin{eqnarray*}
(\id \cotimes \omega)\Delta(f) \cdot g= (\id \cotimes \omega)\gamma_l(f \otimes g) \\
g \cdot (\id \cotimes \omega)\Delta(f) = (\id \cotimes \omega)\rho_l(g \otimes f). 
\end{eqnarray*}
To check that this is indeed a two-sided multiplier observe that  
$$
(f \otimes 1)\gamma_l(g \otimes h) = \rho_l(f \otimes g)(h \otimes 1)
$$
for all $ f,g,h \in H $. 
In a similar way we define $ (\omega \cotimes \id)\Delta(f) \in M(H) $ by 
\begin{eqnarray*}
(\omega \cotimes \id)\Delta(f) \cdot g= (\id \otimes \omega)\gamma_r(f \otimes g) \\
g \cdot (\omega \cotimes \id)\Delta(f) = (\id \otimes \omega)\rho_r(g \otimes f).
\end{eqnarray*}
\begin{definition}
Let $ \Delta: H \rightarrow M(H \cotimes H) $ be a comultiplication such that all Galois maps associated to $ \Delta $ 
define bounded linear maps from $ H \cotimes H $ into itself. \\
A bounded linear functional $ \phi: H \rightarrow \mathbb{C} $ is called left invariant if 
\begin{equation*}
(\id \cotimes \phi)\Delta(f) = \phi(f) 1
\end{equation*}
for all $ f \in H $. Similarly, a bounded linear functional $ \psi: H \rightarrow \mathbb{C} $ is called right invariant if 
\begin{equation*}
(\psi \cotimes \id)\Delta(f) = \psi(f) 1
\end{equation*}
for all $ f \in H $. 
\end{definition}
Let us now give the definition of a bornological quantum group.
\begin{definition}\label{bqgdef}
A bornological quantum group is an essential bornological algebra $ H $ satisfying the approximation property 
together with a comultiplication $ \Delta: H \rightarrow M(H \cotimes H) $ such that all 
Galois maps associated to $ \Delta $ are isomorphisms and a faithful left invariant functional $ \phi: H \rightarrow \mathbb{C} $ . \\
A morphism between bornological quantum groups $ H $ and $ K $ is an essential algebra homomorphism 
$ \alpha: H \rightarrow M(K) $ such that $ (\alpha \cotimes \alpha)\Delta = \Delta \alpha $. 
\end{definition}
To be precise, the Galois maps in a bornological quantum group are supposed to yield bornological isomorphisms of $ H \cotimes H $ 
into itself. The left invariant functional $ \phi $ is also referred to as the left Haar functional. \\ 
Our definition of a bornological quantum group is equivalent to the definition of an algebraic quantum group in the 
sense of van Daele \cite{vD2} provided the underlying bornological vector space carries the fine bornology. The only difference in this case
is that we have included faithfulness of the Haar functional in the definition. 
\begin{lemma} \label{hopflemma1} 
Let $ H $ be a bornological quantum group. Then 
$$
(\rho_l \cotimes \id)(\id \cotimes \gamma_r) = (\id \cotimes \gamma_r)(\rho_l \cotimes \id)
$$
where both sides are viewed as maps from $ H \cotimes H \cotimes H $ into itself. 
\end{lemma}
\proof Using associativity and coassociativity we compute 
\begin{align*}
(\rho_l \cotimes& \id)(\id \cotimes \gamma_r) = 
(\mu \cotimes \id \cotimes \id)(\id \cotimes \Delta \cotimes \id) (\id \cotimes \id \cotimes \mu)(\id \cotimes \Delta \cotimes \id) \\
&= (\mu \cotimes \id \cotimes \id)(\id \cotimes \id \cotimes \id \cotimes \mu) (\id \cotimes \Delta \cotimes \id \cotimes \id)
(\id \cotimes \Delta \cotimes \id) \\
&= (\mu \cotimes \id \cotimes \id)(\id \cotimes \id \cotimes \id \cotimes \mu) (\id \cotimes \id \cotimes \Delta \cotimes \id)
(\id \cotimes \Delta \cotimes \id) \\
&= (\id \cotimes \id \cotimes \mu)(\id \cotimes \Delta \cotimes \id) (\mu \cotimes \id \cotimes \id)
(\id \cotimes \Delta \cotimes \id) \\
&= (\id \cotimes \gamma_r)(\rho_l \cotimes \id)
\end{align*}
which yields the claim. \qed \\
The following theorem provides an alternative description of bornological quantum groups.
\begin{theorem} \label{bqchar}
Let $ H $ be an essential bornological algebra satisfying the approximation property
and let $ \Delta: H \rightarrow M(H \cotimes H) $ be a comultiplication such that all associated Galois maps 
define bounded linear maps from $ H \cotimes H $ to itself. Moreover assume that 
$ \phi: H \rightarrow \mathbb{C} $ is a faithful left invariant functional. 
Then $ H $ is a bornological quantum group iff there exist an essential algebra homomorphism  $ \epsilon: H \rightarrow \mathbb{C} $ and 
a linear isomorphism $ S: H \rightarrow H $ which is both an algebra antihomomorphism and a coalgebra antihomomorphism such that 
\begin{equation*}
(\epsilon \cotimes \id)\Delta = \id = (\id \cotimes \epsilon)\Delta
\end{equation*}
and
\begin{equation*}
\mu(S \cotimes \id) \gamma_r = \epsilon \cotimes \id, \qquad \mu(\id \cotimes S) \rho_l = \id \cotimes \epsilon.
\end{equation*}
In this case the maps $ \epsilon $ and $ S $ are uniquely determined. 
\end{theorem}
\proof The proof follows the discussion in \cite{vD1}. Along the way we obtain some formulas which are also useful in 
other situations. \\
Let us first assume that there exist maps $ \epsilon $ and $ S $ satisfying the above conditions. 
Following the traditional terminology, these maps will be called the counit and the antipode of $ H $. 
We claim that the inverse $ \gamma^{-1}_r $ of $ \gamma_r $ is given by 
$$
\gamma^{-1}_r = (S^{-1} \cotimes \id)\gamma_r^{cop}(S \cotimes \id). 
$$
Using that $ S $ is a coalgebra antihomomorphism we obtain the equality 
\begin{align*}
(S^{-1} \cotimes \id)&\gamma_r^{cop}(S \cotimes \id) = 
(S^{-1} \cotimes \id)(\id \cotimes \mu)(\Delta^{cop} \cotimes \id)(S \cotimes \id) \\
&= (\id \cotimes \mu)(S^{-1} \cotimes \id \cotimes \id)(\tau \Delta S \cotimes \id) \\
&= (\id \cotimes \mu)(\id \cotimes S \cotimes \id)(\Delta \cotimes \id)
\end{align*}
where both sides are viewed as maps from $ H \cotimes H $ to $ M(H \cotimes H) $. In particular, the image of the last map is contained in 
$ H \cotimes H $. We compute 
\begin{align*}
\mu_{(2)}(\gamma_r^{-1} \gamma_r &\cotimes \id_{(2)}) = 
\mu_{(2)}(\id \cotimes \mu \cotimes \id_{(2)})(\id \cotimes S \cotimes \id_{(3)})(\Delta \cotimes \id_{(3)})(\gamma_r \cotimes \id_{(2)}) \\
&= (\mu \cotimes \mu)(\id \cotimes \tau \cotimes \id)(\id \cotimes \mu \cotimes \id_{(2)})(\id \cotimes S \cotimes \id_{(3)})(\Delta \cotimes \id_{(3)})
(\gamma_r \cotimes \id_{(2)}) \\
&= (\id \cotimes \mu)(\id \cotimes \mu \cotimes \id)(\id \cotimes S \cotimes \id_{(2)})
(\gamma_l \cotimes \id_{(2)})(\id \cotimes \tau \cotimes \id)(\gamma_r \cotimes \id_{(2)}) \\
&= (\id \cotimes \mu)(\id \cotimes \mu \cotimes \id)(\id \cotimes S \cotimes \id_{(2)})(\id \cotimes \gamma_r \cotimes \id)
(\id \cotimes \tau \cotimes \id) \gamma_l^{13} \\
&= (\id \cotimes \mu)(\id \cotimes \epsilon \cotimes \id_{(2)})(\id \cotimes \tau \cotimes \id)\gamma_l^{13} \\
&= (\id \cotimes \mu)(\id_{(2)} \cotimes \epsilon \cotimes \id)\gamma_l^{13} \\
&= \mu_{(2)}(\id \cotimes \epsilon \cotimes \id_{(3)})(\Delta \cotimes \id_{(3)}) = \mu_{(2)}
\end{align*}
where we write $ \mu_{(2)} $ for the multiplication in the tensor product $ H \cotimes H $. Similarly we have 
\begin{align*}
\mu_{(2)}&(\id_{(2)} \cotimes \gamma_r \gamma_r^{-1}) = 
\mu_{(2)}(\id_{(3)} \cotimes \mu)(\id_{(2)} \cotimes \Delta \cotimes \id)(\id_{(2)} \cotimes \gamma_r^{-1}) \\
&= (\mu \cotimes \mu)(\id \cotimes \tau \cotimes \id)(\id_{(3)} \cotimes \mu)(\id_{(2)}  \cotimes \Delta \cotimes \id)
(\id_{(2)} \cotimes \gamma_r^{-1}) \\
&= (\mu \cotimes \mu)(\id \cotimes \tau \cotimes \id)(\id_{(3)} \cotimes \mu)(\id_{(2)} \cotimes \Delta \cotimes \id) (\id_{(3)} \cotimes \mu)\\
&\qquad (\id_{(3)} \cotimes S \cotimes \id)(\id_{(2)} \cotimes \Delta \cotimes \id) \\
&= (\mu \cotimes \mu)(\id_{(3)} \cotimes \mu)(\id \cotimes \tau \cotimes \id_{(2)})(\id_{(4)} \cotimes \mu)(\id_{(4)} \cotimes S \cotimes \id) \\
&\qquad (\id_{(2)} \cotimes \Delta \cotimes \id_{(2)})(\id_{(2)} \cotimes \Delta \cotimes \id) \allowdisplaybreaks \\ 
&= (\mu \cotimes \mu)(\id_{(2)} \cotimes \mu \cotimes \mu)(\id_{(4)} \cotimes S \cotimes \id)(\id \cotimes \tau \cotimes \id_{(3)})\\
&\qquad (\id_{(2)} \cotimes \Delta \cotimes \id_{(2)})(\id_{(2)} \cotimes \Delta \cotimes \id) \allowdisplaybreaks \\
&= (\mu \cotimes \mu)(\id_{(3)} \cotimes \mu)(\id_{(3)} \cotimes S \cotimes \id)(\id_{(2)} \cotimes \mu \cotimes \id_{(2)})\\
&\qquad(\id_{(3)} \cotimes \Delta \cotimes \id)(\id \cotimes \tau \cotimes \id_{(2)})(\id_{(2)} \cotimes \Delta \cotimes \id) \\
&= (\mu \cotimes \mu)(\id_{(2)} \cotimes \mu \cotimes \id)(\id_{(3)} \cotimes S \cotimes \id)(\id_{(2)} \cotimes \rho_l \cotimes \id)\\
&\qquad(\id \cotimes \tau \cotimes \id_{(2)})(\id_{(2)} \cotimes \Delta \cotimes \id) \\
&= (\mu \cotimes \mu)(\id_{(3)} \cotimes \epsilon \cotimes \id)(\id \cotimes \tau \cotimes \id_{(2)})(\id_{(2)} \cotimes \Delta \cotimes \id) \\
&= \mu_{(2)}(\id_{(3)} \cotimes \epsilon \cotimes \id)(\id_{(2)} \cotimes \Delta \cotimes \id) = \mu_{(2)}
\end{align*}
which shows that $ \gamma_r $ is an isomorphism. Remark that we did not use the linear functional $ \phi $ in this discussion. \\
In order to treat the other Galois maps one could perform similar calculations. We proceed in a different way and show first that 
the given counit and invertible antipode for $ H $ provide as with counits and antipodes for $ H^{op}, H^{cop} $ and $ H^{opcop} $ as well. 
More precisely, observe that the counit $ \epsilon $ satisfies 
$$
(\epsilon \cotimes \id)\Delta^{cop} = \id = (\id \cotimes \epsilon)\Delta^{cop}
$$
which means that $ \epsilon $ is a counit for $ H^{cop} $ and $ H^{opcop} $. 
Using lemma \ref{Galopcop} we obtain 
$$
\mu^{op}(S \cotimes \id) \gamma_r^{opcop} = \mu \tau (S \cotimes \id) \tau \rho_l \tau = \mu (\id \cotimes S) \rho_l \tau = \epsilon \cotimes \id
$$ 
and 
$$
\mu^{op}(\id \cotimes S) \rho_l^{opcop} = \mu \tau (\id \cotimes S) \tau \gamma_r \tau = \mu (\id \cotimes S) \gamma_r \tau = \id \cotimes \epsilon
$$
which shows that $ S $ is an antipode for $ H^{opcop} $. We compute 
\begin{align*}
\mu(\mu^{op}\cotimes \id)(&S^{-1} \cotimes \id \cotimes \id) (\gamma_r^{op} \cotimes \id) \\
&= \mu(\mu^{op} \cotimes \id)(S^{-1} \cotimes S^{-1} \cotimes \id)(\id \cotimes S \cotimes \id)(\gamma_r^{op} \cotimes \id) \\
&= \mu (S^{-1} \cotimes \id)(\mu \cotimes \id)(\id \cotimes S \cotimes \id)(\gamma_r^{op} \cotimes \id) \\
&= S^{-1}\mu \tau (\mu \cotimes \id)(\id \cotimes S \cotimes S)(\gamma_r^{op} \cotimes \id) \\
&= S^{-1}\mu (\id \cotimes \mu)(\tau \cotimes \id)(\id \cotimes \tau)(\id \cotimes S \cotimes S)(\gamma_r^{op} \cotimes \id) \allowdisplaybreaks \\
&= S^{-1}\mu (\mu \cotimes \id)(S \cotimes \id \cotimes S)(\tau \cotimes \id)(\id \cotimes \tau)(\gamma_r^{op} \cotimes \id) \allowdisplaybreaks\\
&= S^{-1}\mu(\mu \cotimes \id)(S \cotimes \id \cotimes S)(\id \cotimes \gamma_r^{op})(\tau \cotimes \id)(\id \cotimes \tau) \allowdisplaybreaks \\
&= S^{-1}\mu(\mu \cotimes \id)(\id \cotimes S \cotimes \id)(\rho_l \cotimes \id)(S \cotimes \id \cotimes S)(\tau \cotimes \id)(\id \cotimes \tau) \\
&= S^{-1}\mu(\id \cotimes \epsilon \cotimes \id)(S \cotimes \id \cotimes S)(\tau \cotimes \id)(\id \cotimes \tau) \\
&= S^{-1}\mu(\epsilon \cotimes \id \cotimes \id)(\id \cotimes S \cotimes S)(\id \cotimes \tau) \\
&= \mu(S^{-1} \cotimes S^{-1})(S \cotimes S)(\epsilon \cotimes \id \cotimes \id) \\
&= \mu(\epsilon \cotimes \id \cotimes \id)
\end{align*}
which yields 
$$
\mu^{op}(S^{-1} \cotimes \id) \gamma_r^{op} = \epsilon \cotimes \id.
$$ 
Similarly we have
\begin{align*}
\mu(\id \cotimes \mu^{op})(&\id \cotimes \id \cotimes S^{-1}) (\id \cotimes \rho_l^{op}) \\
&= \mu(\id \cotimes \mu^{op})(\id \cotimes S^{-1} \cotimes S^{-1})(\id \cotimes S \cotimes \id)(\id \cotimes \rho_l^{op}) \\
&= \mu (\id \cotimes S^{-1})(\id \cotimes \mu)(\id \cotimes S \cotimes \id)(\id \cotimes \rho_l^{op}) \\
&= S^{-1}\mu \tau (\id \cotimes \mu)(S \cotimes S \cotimes \id)(\id\cotimes \rho_l^{op}) \\
&= S^{-1}\mu (\mu \cotimes \id)(\id \cotimes \tau)(\tau \cotimes \id)(S \cotimes S \cotimes \id)(\id \cotimes \rho_l^{op}) \allowdisplaybreaks\\
&= S^{-1}\mu (\id \cotimes \mu)(S \cotimes \id \cotimes S)(\id \cotimes \tau)(\tau \cotimes \id)(\id \cotimes \rho_l^{op}) \allowdisplaybreaks \\
&= S^{-1}\mu(\id \cotimes \mu)(S \cotimes \id \cotimes S)(\rho_l^{op} \cotimes \id)(\id \cotimes \tau)(\tau \cotimes \id) \allowdisplaybreaks \\
&= S^{-1}\mu(\id \cotimes \mu)(\id \cotimes S \cotimes \id)(\id \cotimes \gamma_r)(S \cotimes \id \cotimes S)(\id \cotimes \tau)(\tau \cotimes \id) \\
&= S^{-1}\mu(\id \cotimes \id \cotimes \epsilon)(S \cotimes S \cotimes \id)(\tau \cotimes \id) \\
&= \mu(S^{-1} \cotimes S^{-1})(S \cotimes S)(\id \cotimes \id \cotimes \epsilon) \\
&= \mu(\id \cotimes \id \cotimes \epsilon)
\end{align*}
which implies 
$$
\mu^{op}(\id \cotimes S^{-1}) \rho_l^{op} = \id \cotimes \epsilon.
$$ 
Hence $ S^{-1} $ is an antipode for $ H^{op} $. As above it follows that $ S^{-1} $ is also an antipode for $ H^{cop} = (H^{op})^{opcop} $. 
We may now apply our previous argument for the Galois map $ \gamma_r $ to $ H^{op}, H^{cop} $ and $ H^{opcop} $ and use lemma \ref{Galopcop}
to see that $ \gamma_l, \rho_l $ and $ \rho_r $ are isomorphisms as well. This shows that $ H $ is a bornological quantum group. \\
Conversely, let us assume that $ H $ is a bornological quantum group and construct the maps $ \epsilon $ and $ S $. We 
begin with the counit $ \epsilon $. 
Choose an element $ h \in H $ such that $ \phi(h) = 1 $ and set 
$$
\epsilon(f) = \phi(\mu \rho_l^{-1}(h \otimes f)).
$$
This yields obviously a bounded linear map $ \epsilon: H \rightarrow \mathbb{C} $. 
Using 
\begin{equation}\label{gammalinear} 
\gamma_r(\id \cotimes \mu) = (\id \cotimes \mu) \gamma_r
\end{equation} 
we easily see that the formula 
$$
E(f) \cdot g = \mu \gamma_r^{-1}(f\otimes g)
$$
defines a left multiplier $ E(f) $ of $ H $. Actually, we obtain a bounded linear map 
$ E: H \rightarrow M_l(H) $ in this way. Using lemma \ref{hopflemma1} we obtain
\begin{align*}
(\id \cotimes \mu)(\id \cotimes &E \cotimes \id)(\rho_l \cotimes \id)(\id \cotimes \gamma_r) = 
(\id \cotimes \mu)(\id \cotimes \gamma_r^{-1})(\rho_l \cotimes \id)(\id \cotimes \gamma_r) \\
&= (\id \cotimes \mu)(\id \cotimes \gamma_r^{-1})(\id \cotimes \gamma_r)(\rho_l \cotimes \id) \\
&= (\id \cotimes \mu)(\rho_l \cotimes \id) \\
&= (\id \cotimes \mu)(\mu \cotimes \id \cotimes \id) (\id \cotimes \Delta \cotimes \id) \\
&= (\mu \cotimes \id)(\id \cotimes \id \cotimes \mu) (\id \cotimes \Delta \cotimes \id) \\
&= (\mu \cotimes \id)(\id \cotimes \gamma_r). 
\end{align*}
Since $ \gamma_r $ and $ \rho_l $ are isomorphisms this implies 
\begin{equation}\label{s2e1}
(\id \cotimes \mu)(\id \cotimes E \cotimes \id) = \mu \rho_l^{-1} \cotimes \id.
\end{equation}
Evaluating equation (\ref{s2e1}) on a tensor $ h \otimes f \otimes g $ where $ h $ is chosen as above and applying
$ \phi \cotimes \id $ we get 
$$
E(f) \cdot g = (\phi \cotimes \id)(h \otimes E(f)\cdot g) = \phi(\mu\rho_l^{-1}(h \otimes f)) g = \epsilon(f)g 
$$
and hence 
\begin{equation}\label{Ee}
E(f) = \epsilon(f)1
\end{equation}
in $ M_l(H) $ for every $ f \in H $.
This shows in particular that we could have used any nonzero bounded linear functional in order to define $ \epsilon $. 
To obtain equation (\ref{Ee}) we did not use the fact that $ \phi $ is left invariant and faithful. \\
According to equation (\ref{Ee}) and the definition of $ E $ we have
\begin{equation}\label{epsilongammar}
(\epsilon \cotimes \id)\gamma_r = \mu(E \cotimes \id)\gamma_r = \mu. 
\end{equation}
Equation (\ref{s2e1}) yields 
$$ 
g \epsilon(f) \otimes h = g \otimes \epsilon(f) h = \mu \rho_l^{-1}(g \otimes f) \otimes h 
$$ 
for all $ f,g \in H $ which implies 
\begin{equation}\label{s2e0}
g \epsilon(f) = \mu \rho_l^{-1}(g \otimes f). 
\end{equation}
This is equivalent to
\begin{equation}\label{s2e2}
(\id \cotimes \epsilon)\rho_l = \mu  
\end{equation}
since $ \rho_l $ is an isomorphism. \\
Let us now show that $ \epsilon $ is an algebra homomorphism. 
We have 
\begin{equation}\label{s2e3}
\rho_l(\id \cotimes \mu) = (\mu \cotimes \id)(\id \cotimes \mu_{(2)})(\id \cotimes \Delta \cotimes \Delta) = \mu_{(2)}(\rho_l \cotimes \Delta) 
\end{equation}
because $ \Delta $ is an algebra homomorphism. According to this relation and equation (\ref{s2e2}) we get
\begin{align*}
(\id \cotimes \epsilon)&\mu_{(2)}(\id \cotimes \id \cotimes \Delta)(\rho_l \cotimes \id) = 
(\id \cotimes \epsilon)\rho_l(\id \cotimes \mu) = \mu(\id \cotimes \mu) \\
&= \mu(\mu \cotimes \id) = \mu(\id \cotimes \epsilon \cotimes \id)
(\rho_l \cotimes \id)
\end{align*}
and since $ \rho_l $ is an isomorphism this implies
$$
(\id \cotimes \epsilon)\mu_{(2)}(\id \cotimes \id \cotimes \Delta) = \mu(\id \cotimes \epsilon \cotimes \id).
$$
Now observe $ \mu_{(2)}(\id \cotimes \id \cotimes \Delta) = (\id \cotimes \mu) \rho_l^{13}$
and hence  
$$
(\id \cotimes \epsilon)(\id \cotimes \mu) \rho_l^{13} = \mu(\id \cotimes \epsilon \cotimes \id) = 
(\id \cotimes \epsilon)(\id \cotimes \epsilon \cotimes \id)\rho_l^{13}
$$
where we use equation (\ref{s2e2}). We deduce 
$$
(\id \cotimes \epsilon)(\id \cotimes \mu) = (\id \cotimes \epsilon)(\id \cotimes \epsilon \cotimes \id)
$$
which implies 
$$
\epsilon(fg) = \epsilon(f) \epsilon(g) 
$$
for all $ f,g \in H $. Thus $ \epsilon $ is an algebra homomorphism. \\
Using this fact and equation (\ref{s2e2}) we calculate 
$$
\mu(\id_{(2)}\cotimes \epsilon)(\id \cotimes \gamma_r) = (\id \cotimes \epsilon)(\id \cotimes \mu)(\rho_l \cotimes \id) = 
(\id \cotimes \mu)(\id \cotimes \epsilon \cotimes \epsilon)(\rho_l \cotimes \id) = \mu \cotimes \epsilon 
$$
which implies 
\begin{equation}\label{epsilongammainv1}
(\id \cotimes \epsilon)\gamma_r = \id \cotimes \epsilon. 
\end{equation}
Analogously one has 
\begin{equation}\label{epsilongammainv2}
(\epsilon \cotimes \id)\rho_l = \epsilon \cotimes \id
\end{equation}
as a consequence of equation (\ref{epsilongammar}). \\
It is easy to see that the map $ \epsilon $ is nonzero. To check that $ \epsilon $ is nondegenerate 
we define a bounded linear map $ \sigma: \mathbb{C} \rightarrow H \cotimes_H \mathbb{C} $ by $ \sigma(1) = k \otimes 1 $ 
where $ k \in H $ is an element satisfying $ \epsilon(k) = 1 $. 
Using equation (\ref{epsilongammar}) and equation (\ref{epsilongammainv1}) we obtain 
$$
\sigma (\epsilon \otimes \id)(f \otimes 1) = \epsilon(f)k \otimes 1 = \mu \gamma_r^{-1}(f \otimes k) \otimes 1 = (\id \cotimes \epsilon)\gamma_r^{-1}(f \otimes k) 
= f \otimes \epsilon(k) 1 = f \otimes 1 
$$
which implies $ H \cotimes_H \mathbb{C} \cong \mathbb{C} $. In a similar way one checks $ \mathbb{C} \cotimes_H H \cong \mathbb{C} $ 
using equation (\ref{epsilongammainv2}). This shows that $ \epsilon $ is nondegenerate. \\ 
According to equation (\ref{s2e2}) we thus have 
\begin{equation} \label{epsdeltar}
(\id \cotimes \epsilon)\Delta = \id
\end{equation}
and using equation (\ref{epsilongammar}) we get 
\begin{equation} \label{epsdeltal}
(\epsilon \cotimes \id) \Delta = \id. 
\end{equation}
Conversely, the last equation implies $ (\epsilon \cotimes \id)\gamma_r = \mu $ which in turn determines $ \epsilon $ uniquely 
since $ \gamma_r $ is an isomorphism. \\
Now we shall construct the antipode. It is easy to check that the formulas 
$$
S_l(f) \cdot g = (\epsilon \cotimes \id) \gamma_r^{-1}(f \otimes g), \qquad 
g \cdot S_r(f) = (\id\cotimes \epsilon) \rho_l^{-1}(g \otimes f)
$$
define a left multiplier $ S_l(f) $ and a right multiplier $ S_r(f) $ of $ H $ for every $ f \in H $. In this way we obtain bounded linear maps 
$ S_l: H \rightarrow M_l(H) $ and $ S_r: H \rightarrow M_r(H) $. \\
Let us show that $ S_l $ is an algebra antihomomorphism. 
Using lemma \ref{hopflemma1} and equation (\ref{s2e2}) we get 
\begin{align*}
(\id \cotimes &\mu)(\id \cotimes S_l \cotimes \id)(\rho_l \cotimes \id) = 
(\id \cotimes \mu)(\id \cotimes S_l \cotimes \id)(\id \cotimes \gamma_r)(\rho_l \cotimes \id)(\id \cotimes \gamma_r^{-1}) \\
&= (\id \cotimes \epsilon \cotimes \id)(\id \cotimes \gamma_r^{-1})(\id \cotimes \gamma_r)(\rho_l \cotimes \id)(\id \cotimes \gamma_r^{-1}) \\
&= (\id \cotimes \epsilon \cotimes \id)(\rho_l \cotimes \id)(\id \cotimes \gamma_r^{-1}) \\
&= (\mu \cotimes \id)(\id \cotimes \gamma_r^{-1}). 
\end{align*}
Applying the multiplication map $ \mu $ to this equation yields
\begin{equation}\label{s2e4}
\mu(\id \cotimes \mu)(\id \cotimes S_l \cotimes \id)(\rho_l \cotimes \id) = \mu(\id \cotimes \mu \gamma_r^{-1}) = 
\mu(\id \cotimes \epsilon \cotimes \id)
\end{equation}
where we use equation \ref{epsilongammar}. 
According to equation (\ref{s2e3}), equation (\ref{s2e4}) and the fact that $ \epsilon $ is an algebra homomorphism 
and another application of equation (\ref{s2e4}) we obtain
\begin{align*}
\mu(\id \cotimes &\mu)(\id \cotimes S_l \cotimes \id)(\mu_{(2)} \cotimes \id)(\rho_l \cotimes \Delta \cotimes \id) \\
&= \mu(\id \cotimes \mu)(\id \cotimes S_l \cotimes \id)(\rho_l \cotimes \id)(\id \cotimes \mu \cotimes \id)\\
&= \mu(\id \cotimes \epsilon \cotimes \id)(\id \cotimes \mu \cotimes \id) \\
&= \mu(\id \cotimes \epsilon \cotimes \epsilon \cotimes \id) \\
&= \mu(\id \cotimes \mu)(\id \cotimes S_l \cotimes \id)(\rho_l \cotimes \id)(\id \cotimes \id \cotimes \epsilon \cotimes \id).
\end{align*}
Since $ \rho_l $ is an isomorphism this yields due to equation (\ref{s2e4})
\begin{align*}
\mu(&\id \cotimes \mu)(\id \cotimes S_l \cotimes \id)(\id \cotimes \mu \cotimes \id)\rho_l^{13} \\
&= \mu(\id \cotimes \mu)(\id \cotimes S_l \cotimes \id)(\mu_{(2)} \cotimes \id)(\id \cotimes \id \cotimes \Delta \cotimes \id) \\
&= \mu(\id \cotimes \mu)(\id \cotimes S_l \cotimes \id)(\id \cotimes \id \cotimes \epsilon \cotimes \id) \\
&= \mu (\id \cotimes \epsilon \cotimes \id)(\id_{(2)} \cotimes \mu)(\id \cotimes \tau \cotimes \id) 
(\id \cotimes S_l \cotimes \id_{(2)}) \\
&=  \mu (\id \cotimes \mu)(\id \cotimes S_l \cotimes \id)(\id_{(2)} \cotimes \mu)
(\rho_l \cotimes \id_{(2)})(\id \cotimes \tau \cotimes \id) 
(\id \cotimes S_l \cotimes \id_{(2)}) \\
&=  \mu (\id \cotimes \mu)(\id \cotimes S_l \cotimes \id)(\id_{(2)} \cotimes \mu)(\id_{(2)} \cotimes S_l \cotimes \id)(\id  \cotimes \tau \cotimes \id)
\rho_l^{13} 
\end{align*}
and hence 
\begin{equation*}
\mu(\id \cotimes \mu)(\id \cotimes S_l \cotimes \id)(\id \cotimes \mu \cotimes \id)
=  \mu (\id \cotimes \mu)(\id \cotimes \mu \cotimes \id)(\id \cotimes \tau \cotimes \id)
(\id \cotimes S_l \cotimes S_l \cotimes \id).
\end{equation*}
Since the algebra $ H $ is nondegenerate we obtain
\begin{equation}\label{s2e4a}
S_l(fg) = S_l(g)S_l(f)
\end{equation}
for all $ f,g \in H $ as claimed. \\
For the map $ S_r $ we do an analogous calculation. We have 
\begin{align*}
(\mu \cotimes \id)&(\id \cotimes S_r \cotimes \id)(\id \cotimes \gamma_r) = (\mu \cotimes \id)(\id \cotimes S_r \cotimes \id)(\rho_l \cotimes \id)
(\id \cotimes \gamma_r)(\rho_l^{-1} \cotimes \id) \\
&= (\id \cotimes \epsilon \cotimes \id)(\id \cotimes \gamma_r)(\rho_l^{-1} \cotimes \id) 
= (\id \cotimes \mu)(\rho_l^{-1} \cotimes \id)
\end{align*}
and applying $ \mu $ yields
\begin{equation}\label{s2e6}
\mu(\mu \cotimes \id)(\id \cotimes S_r \cotimes\id)(\id \cotimes \gamma_r) = \mu(\mu \rho_l^{-1} \cotimes \id) = 
\mu(\id \cotimes \epsilon \cotimes \id). 
\end{equation}
As above one may proceed to show that $ S_r $ is an algebra antihomomorphism. We shall instead first show 
that $ (S_l(f), S_r(f)) $ is a two-sided multiplier of $ H $ for every $ f \in H $. 
By the definition of $ S_r $ we have $ \mu(\id \cotimes S_r) = (\id \cotimes \epsilon) \rho_l^{-1} $ and
hence equation (\ref{s2e4}) implies 
\begin{equation*}
\mu(\id \cotimes \mu)(\id \cotimes S_l \cotimes \id) = \mu(\id \cotimes \epsilon \cotimes \id)(\rho_l^{-1} \cotimes \id) 
= \mu(\mu \cotimes \id)(\id \cotimes S_r \cotimes \id)
\end{equation*}
which is precisely the required identity.  
We can now use equation (\ref{s2e4a}) to obtain that $ S_r $ is an algebra antihomomorphism. 
If $ S: H \rightarrow M(H) $ denotes the linear map 
determined by $ S_l $ and $ S_r $ we have thus showed so far that 
$ S: H \rightarrow M(H) $ is a bounded algebra antihomomorphism. \\
Let us define for $ f \in H $ an element $ \bar{S}_l(f) \in M_l(H) $ and an element $ \bar{S}_r(f) \in M_r(H) $ by 
$$
\bar{S}_l(f) \cdot g = (\epsilon \cotimes \id) \gamma_l^{-1} \tau(f \otimes g), \qquad 
g \cdot \bar{S}_r(f) = (\id \cotimes \epsilon) \rho_r^{-1} \tau(g \otimes f).
$$
According to lemma \ref{Galopcop} we have $ \gamma_l^{-1} \tau = (\gamma_r^{cop})^{-1} $ and 
$ \rho_r^{-1} \tau = (\rho_l^{cop})^{-1} $. The discussion above applied to $ H^{cop} $ shows that 
$ \bar{S}_l $ and $ \bar{S}_r $ determine a bounded algebra antihomomorphism $ \bar{S}: H \rightarrow M(H) $. \\
Our next goal is to prove that $ S $ and $ \bar{S} $ actually define bounded linear maps from $ H $ into itself which are inverse to 
each other. In order to do this observe
$$
(\id \cotimes \mu)(\tau \cotimes \id) = (\id \cotimes \mu)(\tau \cotimes \id)(\gamma_r\cotimes \id)(\gamma_r^{-1} \cotimes \id) 
= \mu_{(2)}(\Delta^{cop} \cotimes \id_{(2)})(\gamma_r^{-1} \cotimes \id)
$$ 
which implies 
\begin{align*}
(\mu &\cotimes \id)(\id \cotimes \bar{S} \cotimes \id)(\id \cotimes \id \cotimes \mu)(\id \cotimes \tau \cotimes \id) \\
&= (\mu \cotimes \id)(\id \cotimes \bar{S} \cotimes \id)(\id \cotimes \mu_{(2)})(\id \cotimes \Delta^{cop} \cotimes \id_{(2)})
(\id \cotimes \gamma_r^{-1} \cotimes \id) \\
&= (\mu \cotimes \id)(\id \cotimes \bar{S} \cotimes \id) (\id \cotimes \mu \cotimes \id)(\id \cotimes \tau \cotimes \id) \\
&\qquad (\id_{(2)} \cotimes \gamma_r^{cop})(\id \cotimes \tau \cotimes \id)(\id \cotimes \gamma_r^{-1} \cotimes \id) \\
&= (\mu \cotimes \id)(\id \cotimes \mu \cotimes \id)(\id \cotimes \bar{S} \cotimes \bar{S} \cotimes \id)
(\id_{(2)} \cotimes \gamma_r^{cop} \cotimes \id)(\id \cotimes \tau \gamma_r^{-1} \cotimes \id) 
\end{align*}
since $ \bar{S} $ is an algebra antihomomorphism. Applying $ \mu $ to this equation yields 
\begin{align*}
\mu(\id &\cotimes \mu)(\id \cotimes \mu \cotimes \id)(\id \cotimes \bar{S} \cotimes \id_{(2)})(\id \cotimes \tau \cotimes \id) \\
&= \mu(\mu \cotimes \id)(\id \cotimes \bar{S} \cotimes \id)(\id \cotimes \id \cotimes \mu)(\id \cotimes \tau \cotimes \id) \\
&= \mu (\mu \cotimes \id)(\id \cotimes \mu \cotimes \id)(\id \cotimes \bar{S} \cotimes \bar{S} \cotimes \id)
(\id_{(2)} \cotimes \gamma_r^{cop})(\id \cotimes \tau \gamma_r^{-1} \cotimes \id) \\
&= \mu (\id \cotimes \mu)(\id_{(2)} \cotimes \mu)(\id_{(2)} \cotimes \bar{S} \cotimes \id)
(\id_{(2)} \cotimes \gamma_r^{cop})(\id \cotimes \bar{S} \cotimes \id_{(2)})(\id \cotimes \tau \gamma_r^{-1} \cotimes \id) \\
&= \mu (\id \cotimes \mu)(\id_{(2)} \cotimes \epsilon \cotimes \id)(\id \cotimes \bar{S} \cotimes \id_{(2)})
(\id \cotimes \tau \cotimes \id)(\id \cotimes \gamma_r^{-1} \cotimes \id) \\
&= \mu (\id \cotimes \mu)(\id \cotimes \bar{S} \cotimes \id)(\id \cotimes \epsilon \cotimes \id_{(2)})(\id \cotimes \gamma_r^{-1} \cotimes \id) \\
&= \mu (\id \cotimes \mu)(\id \cotimes \bar{S} \cotimes \id)(\id \cotimes \mu \cotimes \id)(\id \cotimes S \cotimes \id_{(2)}) 
\end{align*}
where we use the definitions of $ \bar{S} $ and $ S $. As a consequence we obtain the relation 
\begin{equation}\label{s2e7}
\mu(\bar{S} \cotimes \id)\tau = \bar{S} \mu(S \cotimes \id)
\end{equation}
where both sides are viewed as maps from $ H \cotimes H $ into $ M(H) $. 
Choose $ k \in H $ such that $ \epsilon(k) = 1 $. Equation (\ref{s2e7}) together with the definition of $ S $ 
yields the relation
$$
\bar{S}(f) = \bar{S}(f) \epsilon(k) = \bar{S} \mu(S \cotimes \id)\gamma_r(k \otimes f) = \mu(\bar{S} \cotimes \id)\tau \gamma_r(k \otimes f)
$$
for all $ f \in H $. This shows that $ \bar{S} $ defines a bounded linear map from $ H $ to $ H $. 
Replacing $ H $ by $ H^{cop} $ we see that $ S $ may be viewed as a bounded linear map from $ H $ to $ H $ as well. 
Since $ \bar{S} $ is an algebra antihomomorphism 
equation (\ref{s2e7}) then yields 
$$
\mu(\bar{S} \cotimes \id) = \mu(\bar{S} \cotimes \bar{S}S)
$$
and hence 
$$
\mu(\mu \cotimes \id)(\id \cotimes \bar{S} \cotimes \id) = \mu(\mu \cotimes \id)(\id \cotimes \bar{S} \cotimes \bar{S}S). 
$$
Since $ \mu(\id \cotimes \bar{S}) = (\id \cotimes \epsilon)\rho_r^{-1} \tau $ by the definition of $ \bar{S} $ 
we get $ \bar{S}S = \id $. Analogously one obtains $ S \bar{S} = \id $. 
Equation (\ref{s2e6}) and equation (\ref{s2e4}) yield 
\begin{equation*}
\mu(S \cotimes \id) \gamma_r = \epsilon \cotimes \id, \qquad \mu(\id \cotimes S) \rho_l = \id \cotimes \epsilon 
\end{equation*}
as desired. Moreover these equations determine the map $ S $ uniquely. \\
Let us show that $ S $ is a coalgebra antihomomorphism. Since $ \Delta $ is an algebra homomorphism we have 
\begin{align*}
\gamma_r(\mu &\cotimes \id) = (\id \cotimes \mu)(\Delta \cotimes \id)(\mu \cotimes \id) \\
&= (\id \cotimes \mu)(\mu_{(2)} \cotimes \id)(\Delta \cotimes \Delta \cotimes \id) \\
&= (\id \cotimes \mu)(\gamma_l \cotimes \id)(\id \cotimes \gamma_r)
\end{align*}
and using equation (\ref{gammalinear}) we get
\begin{equation}\label{coop111}
(\mu \cotimes \id)(\id \cotimes \gamma_r^{-1}) = (\id \cotimes \mu)(\gamma_r^{-1} \cotimes \id)(\gamma_l \cotimes \id).
\end{equation}
According to lemma \ref{hopflemma1}, equation (\ref{s2e2}) and the definition of $ S $ we have 
\begin{align*}
(\mu \cotimes \id)(\id &\cotimes \gamma_r^{-1}) = (\mu \cotimes \id)(\rho_l^{-1} \cotimes \id)(\id \cotimes \gamma_r^{-1})(\rho_l \cotimes \id) \\
&= (\id \cotimes \epsilon \cotimes \id)(\id \cotimes \gamma_r^{-1})(\rho_l \cotimes \id) \\
&= (\id \cotimes \mu)(\id \cotimes S \cotimes \id)(\rho_l \cotimes \id).
\end{align*}
Together with equation (\ref{coop111}) we obtain
\begin{equation}\label{coop11}
\gamma_r^{-1} \gamma_l = (\id \cotimes S)\rho_l
\end{equation}
and equation (\ref{coop11}) applied to $ H^{cop} $ yields 
\begin{equation}\label{coop12}
\gamma_l^{-1} \gamma_r = (\gamma_r^{cop})^{-1} \tau \tau \gamma_l^{cop} = (\id \cotimes S^{-1})\rho_l^{cop} = (\id \cotimes S^{-1})\tau \rho_r.  
\end{equation}
Equations (\ref{coop11}) and (\ref{coop12}) imply 
\begin{equation}\label{coop1}
\rho_r(S \cotimes S) = (S \cotimes \id)\tau \rho_l^{-1} (S \cotimes \id).
\end{equation}
According to the definition of $ S $ and equation (\ref{epsilongammar}) we have 
\begin{align*}
(\mu \cotimes \id)(\id \cotimes S&\cotimes \id)(\id \cotimes \gamma_r)(\rho_l \cotimes \id) = 
(\mu \cotimes \id)(\id \cotimes S \cotimes \id)(\rho_l \cotimes \id)(\id \cotimes \gamma_r) \\
&= (\id \cotimes \epsilon \cotimes \id)(\id \cotimes \gamma_r) = (\id \cotimes \mu) 
\end{align*}
and hence 
\begin{align*}
(\id \cotimes \mu)(\rho_l^{-1}&\cotimes \id)(\tau \cotimes \id)(\id \cotimes S \cotimes \id) \\
&= (\mu \cotimes \id)(\id \cotimes S \cotimes \id)(\id \cotimes \gamma_r)(\tau \cotimes \id)(\id \cotimes S \cotimes \id) \\
&= (\mu \cotimes \id)(S \cotimes S \cotimes \id)(\id \cotimes \gamma_r)(\tau \cotimes \id) \\
&= (S \cotimes \id)(\mu \cotimes \id)(\tau \cotimes \id)(\id \cotimes \gamma_r)(\tau \cotimes \id) \\
&= (\id \cotimes \mu)(S \cotimes \id_{(2)})(\gamma_l \cotimes \id)
\end{align*}
since $ S $ is an algebra antihomomorphism. Consequently we obtain
\begin{equation}\label{coop2}
(S \cotimes \id) \gamma_l = \rho_l^{-1}(S \cotimes \id) \tau.
\end{equation}
Equations (\ref{coop1}) and (\ref{coop2}) yield
\begin{align*}
\rho_r(S \cotimes S) &\tau = (S \cotimes \id)\tau \rho_l^{-1} (S \cotimes \id)\tau \\
&= (S \cotimes S)(\id \cotimes S^{-1})\tau \rho_l^{-1}(S \cotimes \id) \tau = (S \cotimes S) \tau \gamma_l
\end{align*}
and using that $ S $ is an algebra antihomomorphism we get
\begin{align*}
(\id \cotimes \mu^{op})(\Delta &\cotimes \id)(S \cotimes \id)(\id \cotimes S) = (\id \cotimes \mu)(\id \cotimes \tau)(\Delta \cotimes \id)(S \cotimes S) \\
&= (\id \cotimes \mu)(\tau \cotimes \id)(\id \cotimes \Delta)\tau(S \cotimes S) \\
&= \rho_r(S \cotimes S) \tau \\
&= (S \cotimes S)\tau \gamma_l \allowdisplaybreaks \\
&= (S \cotimes S)(\id \cotimes \mu)(\tau \cotimes \id)(\Delta \cotimes \id) \allowdisplaybreaks \\
&= (S \cotimes \id)(\id \cotimes \mu)(\id \cotimes S \cotimes S)(\id \cotimes \tau)(\tau \cotimes \id)(\Delta \cotimes \id) \\
&= (\id \cotimes \mu)(\id \cotimes \tau)(S \cotimes S \cotimes S)(\tau \cotimes \id)(\Delta \cotimes \id) \\
&= (\id \cotimes \mu^{op})(S \cotimes S \cotimes \id)(\tau \cotimes \id)(\Delta \cotimes \id)(\id \cotimes S). 
\end{align*}
This implies 
$$
\Delta S = (S \cotimes S) \tau \Delta
$$
and shows that $ S $ is a coalgebra antihomomorphism. Of course $ S^{-1} $ is a coalgebra antihomomorphism as well. 
We have thus shown that there exist unique maps $ S $ and $ \epsilon $ with 
the desired properties. This finishes the proof. \qed \\
Recall that a morphism of bornological quantum groups is an essential algebra homomorphism $ \alpha: H \rightarrow M(K) $ 
which is also a coalgebra homomorphism. 
\begin{prop}\label{morphism}
Every morphism $ \alpha: H \rightarrow M(K) $ of bornological quantum groups is automatically compatible with the counits and the antipodes. 
\end{prop}
\proof Note that the Galois maps associated to the comultiplication of $ H $ extend to bounded linear maps from $ M(H) \cotimes M(H) $ into 
$ M(H \cotimes H) $. Moreover observe that the relation 
$ (\epsilon \cotimes \id) \gamma_r = \mu $ obtained in equation (\ref{epsilongammar}) still holds when we consider both sides as 
maps form $ M(H) \cotimes M(H) $ to $ M(H) $. \\
Since $ \alpha: H \rightarrow M(K) $ is an algebra homomorphism and a coalgebra homomorphism we have 
\begin{equation}\label{Galoismorphism}
\gamma_r (\alpha \cotimes \alpha) = (\alpha \cotimes \alpha) \gamma_r
\end{equation}
where both sides are viewed as maps from $ H \cotimes H $ into $ M(K \cotimes K) $. Hence we obtain
\begin{equation*}
(\epsilon \cotimes \id)(\alpha \cotimes \alpha) \gamma_r = (\epsilon \cotimes \id)\gamma_r(\alpha \cotimes \alpha) = \mu(\alpha \cotimes \alpha) = \alpha \mu 
= (\epsilon \cotimes \alpha)\gamma_r
\end{equation*}
where all maps are considered to be defined on $ H \cotimes H $ with values in $ M(K) $. We conclude 
$$
(\epsilon \alpha) \cotimes \alpha = \epsilon \cotimes \alpha
$$
because $ \gamma_r $ is an isomorphism. Since $ \alpha $ is nondegenerate this shows $ \epsilon \alpha = \epsilon $ which means that $ \alpha $ 
is compatible with the counits. \\
The arguments given in the proof of theorem \ref{bqchar} show that the inverses of the Galois maps of $ H $ can be 
described explicitly using the antipode $ S $ and its inverse. It follows that these maps 
are defined on $ M(H) \cotimes M(H) $ in a natural way. \\
With this in mind and using 
\begin{align*}
(\mu \cotimes \id)(\id \cotimes &\rho_r)(\id \cotimes \mu \cotimes \id)(\id \cotimes \tau \cotimes \id)(\rho_l \cotimes \id_{(2)}) \\
&= \mu_{(2)}(\id_{(2)} \cotimes \mu_{(2)})(\id_{(2)} \cotimes \Delta \cotimes \Delta)(\id \cotimes \tau \cotimes \id)
\end{align*}
we compute on $ H \cotimes H \cotimes M(H) \cotimes H $ 
\begin{align*}
(\mu \cotimes \id)&(\id \cotimes \mu \cotimes \id)(\id_{(2)} \cotimes \gamma_r^{-1}) = 
(\mu \cotimes \id)(\mu \cotimes \id_{(2)})(\id_{(2)} \cotimes \gamma_r^{-1}) \\
&= (\mu \cotimes \id)(\id \cotimes \gamma_r^{-1})(\id \cotimes \epsilon \cotimes \id_{(2)})(\rho_l \cotimes \id_{(2)}) \\
&= (\mu \cotimes \id)(\id \cotimes \gamma_r^{-1})(\id \cotimes \tau)(\id \cotimes \mu \cotimes \id)(\id \cotimes S \cotimes \id_{(2)}) \\
&\qquad (\id \cotimes \gamma_r \cotimes \id)(\id_{(2)} \cotimes \tau)(\rho_l \cotimes \id_{(2)}) \allowdisplaybreaks\\
&= (\mu \cotimes \id)(\id \cotimes \gamma_r^{-1})(\id \cotimes \tau)(\id \cotimes \mu \cotimes \id)(\id \cotimes S \cotimes \id_{(2)}) \\
&\qquad (\rho_l \cotimes \id_{(2)})(\id \cotimes \gamma_r \cotimes \id)(\id_{(2)} \cotimes \tau) \allowdisplaybreaks\\
&= (\mu \cotimes \id)(\id_{(2)} \cotimes S)(\id \cotimes \rho_r)(\id \cotimes S^{-1} \cotimes \id)(\id \cotimes \mu \cotimes \id)\\
&\qquad (\id \cotimes S \cotimes \id_{(2)}) (\rho_l \cotimes \id_{(2)})(\id \cotimes \gamma_r \cotimes \id)
(\id_{(2)} \cotimes \tau) \allowdisplaybreaks \\
&= (\id \cotimes S)(\mu \cotimes \id)(\id \cotimes \rho_r)(\id \cotimes \mu \cotimes \id)(\id \cotimes \tau \cotimes \id)\\
&\qquad (\rho_l \cotimes \id_{(2)})(\id_{(2)} \cotimes S^{-1} \cotimes \id)(\id \cotimes \gamma_r \cotimes \id)
(\id_{(2)} \cotimes \tau) \allowdisplaybreaks \\
&= (\id \cotimes S)\mu_{(2)}(\id_{(2)} \cotimes \mu_{(2)})(\id_{(2)} \cotimes \Delta \cotimes \Delta)(\id \cotimes \tau \cotimes \id)  \\
&\qquad (\id_{(2)} \cotimes S^{-1} \cotimes \id)(\id \cotimes \gamma_r \cotimes \id)(\id_{(2)} \cotimes \tau) \allowdisplaybreaks \\
&= (\mu \cotimes \id)(\id_{(2)} \cotimes \mu)(\id_{(2)} \cotimes S \cotimes S)(\id_{(2)} \cotimes \tau)(\id \cotimes \tau \cotimes \id)
(\id_{(2)} \cotimes \mu_{(2)}) \\
&\qquad (\id_{(2)} \cotimes \Delta \cotimes \Delta)(\id \cotimes \tau \cotimes \id) (\id_{(2)} \cotimes S^{-1} \cotimes \id)
(\id \cotimes \gamma_r \cotimes \id)(\id_{(2)} \cotimes \tau) \allowdisplaybreaks  \\
&= (\mu \cotimes \id)(\id_{(2)} \cotimes \mu)(\id_{(2)} \cotimes S \cotimes \id)(\id \cotimes \mu_{(2)} \cotimes \id)
(\id \cotimes \Delta \cotimes \Delta \cotimes \id) \\
&\qquad (\id_{(2)} \cotimes \tau)(\id \cotimes \gamma_r \cotimes \id)(\id_{(2)} \cotimes \tau) \\
&= (\mu \cotimes \id)(\id \cotimes \gamma_r^{-1})(\id \cotimes \mu \cotimes \id)\gamma_r^{24}.
\end{align*}
This yields
$$
(\mu \cotimes \id)(\id \cotimes \gamma_r^{-1}) = \gamma_r^{-1}(\mu \cotimes \id)\gamma_r^{13}
$$
on $ H \cotimes M(H) \cotimes H $ 
which in turn implies
\begin{equation}\label{alphaS}
(\mu \cotimes \mu)(\id \cotimes \gamma_r^{-1} \cotimes \id) = \gamma_r^{-1}(\mu \cotimes \id)\gamma_r^{13}(\id_{(2)} \cotimes \mu)
\end{equation}
if both sides are viewed as maps from $ H \cotimes M(H) \cotimes M(H) \cotimes H $ into $ H \cotimes H $. 
Moreover we have
\begin{align*}
\gamma_r(\mu &\cotimes \mu)(\id \cotimes \alpha \cotimes \alpha \cotimes \id)= (\id \cotimes \mu)(\id_{(2)} \cotimes \mu)
(\Delta \cotimes \id_{(2)})(\mu \cotimes \id_{(2)})(\id \cotimes \alpha \cotimes \alpha \cotimes \id) \\
&= (\id \cotimes \mu)(\id_{(2)} \cotimes \mu)(\mu_{(2)} \cotimes \id_{(2)})(\Delta \cotimes \Delta \cotimes \id_{(2)})
(\id \cotimes \alpha \cotimes \alpha \cotimes \id) \\
&= \mu_{(2)}(\Delta \cotimes \id_{(2)})(\id_{(2)} \cotimes \mu)(\id_{(2)} \cotimes \mu \cotimes \id)
(\id \cotimes \alpha \cotimes \alpha \cotimes \alpha \cotimes \id)(\id \cotimes \Delta \cotimes \id_{(2)}) \\
&= \mu_{(2)}(\Delta \cotimes \id_{(2)})(\id_{(2)} \cotimes \mu)
(\id \cotimes \alpha \cotimes \alpha \cotimes \id)(\id \cotimes \gamma_r \cotimes \id) \\
&= (\mu \cotimes \id)\gamma_r^{13}(\id_{(2)} \cotimes \mu)(\id \cotimes \alpha \cotimes \alpha \cotimes \id)(\id \cotimes \gamma_r \cotimes \id) 
\end{align*}
as maps from $ K \cotimes H \cotimes H \cotimes K $ into $ K \cotimes K $ which yields according to equation (\ref{alphaS}) 
\begin{align*}
(\mu &\cotimes \mu)(\id \cotimes \alpha \cotimes \alpha \cotimes \id)= 
\gamma_r^{-1} \gamma_r(\mu \cotimes \mu)(\id \cotimes \alpha \cotimes \alpha \cotimes \id) \\
&= \gamma_r^{-1}(\mu \cotimes \id)\gamma_r^{13}(\id_{(2)} \cotimes \mu)(\id \cotimes \alpha \cotimes \alpha \cotimes \id)(\id \cotimes \gamma_r \cotimes \id) \\
&= (\mu \cotimes \mu)(\id \cotimes \gamma_r^{-1} \cotimes \id)(\id \cotimes \alpha \cotimes \alpha \cotimes \id)
(\id \cotimes \gamma_r \cotimes \id) 
\end{align*}
and we deduce 
\begin{equation}\label{Galoisinversemorphism}
\gamma_r^{-1} (\alpha \cotimes \alpha) = (\alpha \cotimes \alpha) \gamma_r^{-1}.
\end{equation}
Note that this assertion does not immediately follow from equation (\ref{Galoismorphism}) since the 
map $ \gamma_r^{-1} $ is not defined on the multiplier algebra $ M(K \cotimes K) $. \\
Now remark that the relation
\begin{equation*}
\mu(S \cotimes \id) = (\epsilon \cotimes \id)\gamma_r^{-1}
\end{equation*}
still holds if both sides are viewed as maps from $ M(H) \cotimes M(H) $ into $ M(H) $. Using this observation we obtain
the relations 
\begin{equation*}
\mu(S \cotimes \id)(\alpha \cotimes \alpha) = (\epsilon \cotimes \id)\gamma_r^{-1}(\alpha \cotimes \alpha)
\end{equation*}
and 
\begin{equation*}
\alpha\mu(S \cotimes \id) = (\epsilon \cotimes \alpha)\gamma_r^{-1}. 
\end{equation*}
These equations, together with the fact that $ \alpha $ is compatible with the counits and equation (\ref{Galoisinversemorphism}), yield
$$
\mu(\alpha S \cotimes \alpha) = \alpha \mu(S \cotimes \id) = (\epsilon \cotimes \id)(\alpha \cotimes \alpha)\gamma_r^{-1} = 
(\epsilon \cotimes \id)\gamma_r^{-1}(\alpha \cotimes \alpha) = \mu(S \alpha \cotimes \alpha).
$$
Since $ \alpha $ is nondegenerate it follows that the maps $ \alpha S $ and $ S \alpha $ coincide. 
Remark that we also have $ \alpha S^{-1} = S^{-1} \alpha $. 
In other words, $ \alpha $ is compatible with the antipodes. This finishes the proof. \qed \\
Proposition \ref{morphism} implies in particular another relation between the antipode and 
the counit in a bornological quantum group $ H $. More precisely, we have $ \epsilon S = \epsilon $ and $ \epsilon S^{-1} = \epsilon $ since
$ S: H \rightarrow H^{opcop} $ is an isomorphism of bornological quantum groups. 

\section{Modular properties of the integral}\label{secmodhaar}

In this section we discuss modular properties of the Haar functional on a bornological quantum group. 
The results are parallel to the ones for algebraic quantum groups. \\
Let $ H $ be a bornological quantum group. If $ \phi $ is a left invariant functional one can show 
easily that $ S(\phi) $ is right invariant where $ S(\phi)(f) = \phi(S(f)) $. In particular, there 
always exists a faithful right invariant functional on $ H $. 
\begin{prop}\label{haarnu}
Let $ H $ be a bornological quantum group and let $ \phi $ and $ \psi $ be faithful left and right invariant functionals on $ H $, respectively. 
There exists a bornological isomorphism $ \nu $ of $ H $ such that 
$$
\psi(hf) = \phi(h \nu(f))
$$
for all $ f,h \in H $. 
\end{prop}
\proof Choose $ c \in H $ with $ \phi(c) = 1 $. Using invariance we obtain 
$$
\psi(hf) \phi(c) = (\psi \cotimes \phi)(\Delta(h) \Delta(f)(1 \otimes c)) = (\psi \cotimes \phi)((1 \otimes h) \tau \gamma_l^{-1}(\Delta(f)(1 \otimes c))). 
$$
If we set 
$$ 
\nu(f) = (\psi \cotimes \id) \tau \gamma_l^{-1}(\Delta(f)(1 \otimes c)) = (\psi \cotimes \id) \tau \gamma_l^{-1} \gamma_r(f \otimes c)
$$ 
we thus get $ \psi(hf) = \phi(h\nu(f)) $ for all $ f, h \in H $. Remark that we also have $ \psi(f) = \phi(\nu(f)) $ for all $ f \in H $. 
In a similar way one shows that 
there is a bounded linear endomorphism $ \lambda $ of $ H $ such that $ \psi(h\lambda(f)) = \phi(hf) $. 
Since $ \phi $ and $ \psi $ are both faithful we deduce that the maps $ \nu $ and $ \lambda $ are inverse to each other. \qed 
\begin{prop} \label{prophaarunique} 
Let $ H $ be a bornological quantum group. Then the left Haar functional $ \phi $ on $ H $ is unique up to a scalar.
\end{prop}
\proof Let $ \phi_1 $ and $ \phi_2 $ be faithful left invariant functionals and let $ \psi $ be a faithful right invariant 
functional. We may choose $ f_1 $ and $ g $ such that $ \psi(gf_1) = 1 $. By proposition (\ref{haarnu}) there exist 
$ k $ and $ f_2 $ such that $ \phi_1(hf_1) = \psi(hk) = \phi_2(hf_2) $ for all $ h \in H $. Consider the element 
$ \rho_l^{-1}((1 \otimes h)\Delta(g)) $. Multiplying this element with $ f_1 $ and $ f_2 $ on the right and applying 
$ \psi \cotimes \phi_1 $ and $ \psi \cotimes \phi_2 $, respectivly, yields 
$$
\psi(gf_1)\phi_1(h) = (\psi \cotimes \phi_1)(\rho_l^{-1}((1 \otimes h)\Delta(g))(1 \otimes f_1))
$$
and 
$$
\psi(gf_2)\phi_2(h) = (\psi \cotimes \phi_2)(\rho_l^{-1}((1 \otimes h)\Delta(g))(1 \otimes f_2)) 
$$
using invariance. According to the choice of $ f_1 $ and $ f_2 $ this implies $ \phi_1(h) = \psi(gf_2) \phi_2(h) $ for all 
$ h \in H $ and yields the claim. \qed 
\begin{prop}\label{sigmadef}
Let $ H $ be a bornological quantum group. There exists a unique bounded algebra automorphism $ \sigma $ of $ H $ 
such that 
$$
\phi(f g) = \phi(g \sigma(f))
$$
for all $ f, g \in H $. Moreover $ \phi $ is invariant under $ \sigma $. 
\end{prop}
\proof Let $ \psi $ be a faithful right invariant functional on $ H $. Using the relation $ (\psi \cotimes \id)\gamma_r = \psi \cotimes \id $ and 
equation (\ref{coop11}) we get
\begin{equation}\label{sigma2}
S(\psi \cotimes \id) \rho_l = S(\psi \cotimes \id)(\id \cotimes S^{-1})\gamma_r^{-1} \gamma_l = (\psi \cotimes \id) \gamma_l. 
\end{equation}
Similarly, using $ (\id \cotimes \phi)\gamma_l = (\id \cotimes \phi) \tau $ and equation (\ref{coop12}) 
we obtain 
\begin{equation}\label{sigma1}
S(\id \cotimes \phi) \gamma_r = S(\id \cotimes \phi)\gamma_l \tau(S^{-1} \cotimes \id) \rho_r = (\id \cotimes \phi) \rho_r. 
\end{equation}
Due to equation (\ref{sigma2}) and equation (\ref{sigma1}) we have
\begin{align*}
(\psi \cotimes \phi)((h &\otimes x)(\id \cotimes S)\Delta(y)) = \phi(x S(\psi \cotimes \id) \rho_l (h \otimes y)) \\
&= \phi(x(\psi \cotimes \id)\gamma_l(h \otimes y)) \\
&= (\psi \cotimes \phi)((1 \otimes x) \Delta(h)(y \otimes 1)) \\
&= \psi((\id \cotimes \phi)\rho_r(x \otimes h) y) \\
&= \psi(S(\id \cotimes \phi) \gamma_r(x \otimes h) y) \\
&= (\psi \cotimes \phi)(S \cotimes \id) (\Delta(x)(y \otimes h)).
\end{align*}
If we set 
$$
f = (\psi \cotimes \id)((S \cotimes \id)\Delta(x)(y \otimes 1)) = (\psi \cotimes \id)(S \cotimes \id)\gamma_l(x \otimes y)
$$
and 
$$
g = (\id \cotimes \phi)((1 \otimes x)(\id \cotimes S)\Delta(y)) = (\id \cotimes \phi)(S^{-1} \cotimes \id)\tau \rho_l(\id \cotimes S)(x \otimes y)
$$
this means $ \psi(hg) = \phi(fh) $ for all $ h \in H $. Choosing $ c $ such that $ \psi(c) = 1 $ we see that the formula
$$
\eta(f) = (\id \cotimes \phi)(S^{-1} \cotimes \id)\tau \rho_l (\id \cotimes S) \gamma_l^{-1}(S^{-1} \cotimes \id)(c \otimes f)
$$
defines a bounded linear endomorphism $ \eta $ of $ H $ such that 
$$ 
\phi(f h) = \psi(h\eta(f)) 
$$ 
for all $ h \in H $. Using the map $ \nu $ obtained in proposition \ref{haarnu} we deduce that 
$ \sigma = \nu \eta $ is an endomorphism of $ H $ such that $ \phi(fh) = \phi(h \sigma(f)) $ for all $ f,h \in H $. 
A similar argument shows that the map $ \sigma $ is a bornological isomorphism. \\
Uniqueness of $ \sigma $ follows from the faithfulness of $ \phi $. Let us show that $ \sigma $ is multiplicative. 
We compute 
$$
\phi(h\sigma(f) \sigma(g)) = \phi(gh \sigma(f)) = \phi(fg h) = \phi(h \sigma(fg)) 
$$
for all $ f,g, h \in H $ which yields the claim since $ \phi $ is faithful. 
The last assertion follows from the relation
$$
\phi(f g) = \phi(g \sigma(f)) = \phi(\sigma(f) \sigma(g)) = \phi(\sigma(fg)) 
$$
for all $ f, g \in H $ and the fact that $ H \cdot H $ is dense in $ H $. \qed \\
Let $ \psi $ be a right Haar measure on $ H $. Then $ S(\psi) $ is a left Haar measure and we obtain
$$
\psi(S(f) S(g)) = S(\psi)(gf) = S(\psi)(\sigma^{-1}(f) g) = \psi(S(g) (S \sigma^{-1} S^{-1}) S(f))) 
$$
for all $ f, g \in H $ and thus 
$$
\psi(fg) = \psi(g \rho(f))
$$
where $ \rho = S \sigma^{-1} S^{-1} $. If we also consider $ S^{-1}(\psi) $ and use proposition \ref{prophaarunique} 
we get $ \rho = S^{-1} \sigma^{-1} S $ which yields the relation
$$
S^2 \sigma = \sigma S^2
$$
for the automorphism $ \sigma $. \\
Next we shall introduce the modular element of a bornological quantum group. 
\begin{prop} \label{deltaconstruction} 
Let $ H $ be a bornological quantum group. There exists a unique multiplier $ \delta \in M(H) $ 
such that 
$$
(\phi \cotimes \id)\Delta(f) = \phi(f) \delta
$$
for all $ f \in H $.
\end{prop}
\proof For every $ f \in H $ we have a multiplier $ \delta_f $ defined by $ \delta_f = (\phi \otimes \id)\Delta(f) $. 
We have to show that $ \delta_f = \phi(f) \delta $ for some $ \delta \in M(H)$. Remark that $ \delta $ is uniquely 
determined by this equation. For any element $ g \in H $ the functional $ \phi_g $ given by 
$$ 
\phi_g(f) = (\phi \cotimes \phi)(1 \otimes g)\Delta(f)) = \phi(g\delta_f) 
$$ 
is easily seen to be left invariant. According to proposition \ref{prophaarunique} there exists a scalar $ \lambda_g $ depending on $ g $ such that 
$ \phi_g(f) = \lambda_g \phi(f) $ for all $ f \in H $. This implies 
$$ 
\phi(g \delta_f) \phi(h) = \lambda_g \phi(f) \phi(h) = \phi(g \delta_h) \phi(f) 
$$ 
for all 
$ f,g,h \in H $. Hence $ \phi(g (\delta_f \phi(h) - \delta_h \phi(f)) = 0 $ for all $ g \in H $. Inserting $ g = xy $ 
and using that $ \phi $ is faithful we get $ y(\delta_f \phi(h) - \delta_h \phi(f)) = 0 $ for all $ y \in H $. 
Multiplying this equation with $ z $ from the right we obtain 
$ (\delta_f \phi(h) - \delta_h \phi(f))z = 0 $ for all $ z \in H $. Hence the multipliers $ \delta_f \phi(h) $ and 
$ \delta_h \phi(f) $ are equal. Choose $ h \in H $ satisfying $ \phi(h) = 1 $ and set $ \delta = \delta_h $. 
Then we obtain $ \delta_f = \phi(f) \delta $ for all $ f \in H $ as desired. This yields the claim. \qed \\
The multiplier $ \delta $ is called the modular element of $ H $. 
\begin{prop} \label{deltaprop}
The modular element $ \delta $ is invertible and satisfies the relations
\begin{equation*}
\Delta(\delta) = \delta \otimes \delta, \quad \epsilon(\delta) = 1, \quad S(\delta) = \delta^{-1}.
\end{equation*}
\end{prop}
\proof Apply $ \Delta $ to the defining formula in proposition \ref{deltaconstruction} to obtain 
$$
\phi(h) \Delta(\delta) = (\phi \cotimes \id \cotimes \id)(\id \cotimes \Delta)\Delta(h) = \delta \otimes (\phi \cotimes \id)\Delta(h) 
= \phi(h)\, \delta \otimes \delta. 
$$
Choosing $ h $ such that $ \phi(h) = 1 $ yields the first equation. Similarly, the relation $ \epsilon(\delta) = 1 $ follows by 
applying $ \epsilon $ to the formula $ (\phi \cotimes \id)\Delta(f) = \phi(f) \delta $. To prove the last relation observe that we have 
$$
\epsilon(h) \epsilon(\delta) f = \mu(S \cotimes \id)(\Delta(h\delta)(1\otimes f)) = \mu(S \cotimes \id)(\Delta(h)(\delta \otimes \delta f)) = 
S(\delta) \epsilon(h) \delta f 
$$
which implies $ f = S(\delta) \delta f $ for all $ f \in H $. By faithfulness of the multiplication we obtain 
$ f = f S(\delta) \delta $ for all $ f \in H $ as well and hence $ S(\delta) \delta = 1 $. Similarly one obtains 
$ \delta S(\delta) = 1 $ which shows that $ \delta $ is invertible with inverse $ S(\delta) $. \qed \\
Observe that we also have $ S^{-1}(\delta) = \delta^{-1} $. If $ \psi $ is a faithful right invariant functional then 
the formula
\begin{equation*}
(\id \cotimes \psi)\Delta(f) = \psi(f) \delta^{-1}
\end{equation*}
describes the corresponding modular relation. This follows from proposition \ref{deltaconstruction} and proposition \ref{deltaprop} 
using the left invariant functional $ \phi = S(\psi) $. 

\section{Modules and comodules}\label{secmodcomod}

In this section we discuss the concepts of an essential module and an essential comodule over a bornological quantum group. \\
We begin with the notion of an essential module. Actually, the definition of essential modules over bornological algebras was 
already given in section \ref{secnondeg}. 
\begin{definition}
Let $ H $ be a bornological quantum group. An essential $ H $-module is an $ H $-module $ V $ 
such that the module action induces a bornological isomorphism $ H \cotimes_H V \cong V $. 
A bounded linear map $ f: V \rightarrow W $ between essential $ H $-modules is called $ H $-linear 
if the diagram
$$
\xymatrix{
H \cotimes V \ar@{->}[r]^{\mu_V} \ar@{->}[d]^{\id \cotimes f} & 
       V \ar@{->}[d]^f \\
       H \cotimes W \ar@{->}[r]^{\mu_W} & 
       W  
     }
$$
is commutative. 
\end{definition}
If $ H $ and the $ H $-module $ V $ carry the fine bornology, then $ V $ is essential iff $ HV = V $. This follows 
easily from the fact that $ H $ has an approximate identity in this case. Modules satisfying 
the condition $ HV = V $ are called unital in \cite{DvDZ}. Hence unital modules over an algebraic quantum
group are essential. \\
We denote the category of essential $ H $-modules by $ H \LSMod $. 
By definition, the morphisms in $ H \LSMod $ are the bounded $ H $-linear maps. 
We will also speak of $ H $-equivariant maps or $ H $-module maps instead of $ H $-linear maps. \\
There are some basic constructions with essential modules. The direct sum of a family of essential $ H $-modules is again an
essential $ H $-module. Using the quantum group structure of $ H $ one obtains a tensor 
product in $ H \LSMod $. More precisely, let $ V $ and $ W $ be essential $ H $-modules. As in the proof of 
lemma \ref{lemnondeg3} one obtains that there is a natural isomorphism 
$$
V \cotimes W \cong (H \cotimes H) \cotimes_{H \cotimes H} V \cotimes W
$$
of $ H \cotimes H $-modules. Moreover $ H \cotimes H $ can be viewed as an essential 
$ H $-module using the comultiplication.
The $ H $-module structure on $ V \cotimes W $ is defined by the map 
$$
H \cotimes_H (V \cotimes W) \cong H \cotimes_H (H \cotimes H) \cotimes_{H \cotimes H} V \cotimes W \cong (H \cotimes H) \cotimes_{H \cotimes H} V \cotimes W \cong V \cotimes W
$$
which shows at the same time that $ V \cotimes W $ is again an essential $ H $-module. 
Remark that the trivial one-dimensional $ H $-module $ \mathbb{C} $ given by the counit $ \epsilon $ 
behaves like a unit with respect to the tensor product. \\
Dually to the concept of an essential module one has the notion of an essential comodule. Let $ H $ be a bornological quantum 
group, let $ V $ be a bornological vector space and 
let $ \Hom_H(H, V \cotimes H) $ be the space of bounded right $ H $-linear maps from $ H $ to $ V \cotimes H $. 
A coaction of $ H $ on $ V $ is a bounded linear map $ \eta: V \rightarrow \Hom_H(H, V \cotimes H) $ 
which is colinear in the following sense. By adjoint associativity, the map $ \eta $ can equivalently
be described as a bounded $ H $-linear map $ V \cotimes H \rightarrow V \cotimes H $. Then $ \eta $ 
is said to be $ H $-colinear if the latter map is an isomorphism and satisfies the relation 
$$
(\id \otimes \gamma_r) \eta_{12} (\id \otimes \gamma_r^{-1}) = \eta_{12} \eta_{13}
$$
where both sides are viewed as maps from $ V \cotimes H \cotimes H $ to itself. 
\begin{definition}
Let $ H $ be a bornological quantum group. An essential $ H $-comodule is a bornological vector space 
$ V $ together with a coaction $ \eta: V \rightarrow \Hom_H(H, V \cotimes H) $. A bounded linear map 
$ f: V \rightarrow W $ between essential comodules is called $ H $-colinear if the diagram
$$
\xymatrix{
V \cotimes H \ar@{->}[r]^{\eta_V} \ar@{->}[d]^{f \otimes \id} & 
       V \cotimes H \ar@{->}[d]^{f \otimes \id} \\
       W \cotimes H \ar@{->}[r]^{\eta_W} & 
       W \cotimes H 
     }
$$
is commutative. 
\end{definition}
We write $ \SComodR H $ for the category of essential comodules over $ H $ with $ H $-colinear maps as morphisms. 
More precisely, we have defined right comodules. There are analogous definitions for left comodules. Let us point 
out that corepresentations and comodules in the framework of multiplier Hopf algebras have been discussed in 
detail in \cite{vDZ1}, \cite{vDZ2}. \\
The most elementary example of a coaction is the trivial coaction $ \tau $ of $ H $ on $ V $. 
The map $ \tau: V \rightarrow \Hom_H(H, V \cotimes H) $ is given by $ \tau(v)(f) = v \otimes f $. 
Equivalently, the linear map $ V \cotimes H \rightarrow V \cotimes H $ corresponding to $ \tau $ is the 
identity. \\
As in the case of essential modules, there exists a tensor product in the category of essential comodules. 
Assume that $ \eta_V: V \cotimes H \rightarrow V \cotimes H $ and $ \eta_W: W \cotimes H \rightarrow W \cotimes H $ 
are essential comodules. Then the tensor product coaction $ \eta_{V \cotimes W} $ is defined as the composition 
$$
 \xymatrix{
     V \cotimes W \cotimes H \; \ar@{->}[r]^{\eta_W^{23}} & V \cotimes W \cotimes H \ar@{->}[r]^{\eta_V^{13}} & V \cotimes W \cotimes H 
     }
$$
It is clear that $ \eta_{V \cotimes W} $ is a right $ H $-linear isomorphism and a straightforward calculation shows that 
it is indeed a coaction. The trivial coaction on $ \mathbb{C} $ behaves like a unit with respect to the tensor 
product of comodules. \\ 
An important example of a coaction is the regular coaction of $ H $ on itself given by the comultiplication 
$ \Delta: H \rightarrow M(H \cotimes H) $. More precisely, the regular coaction is the map from $ H $ to 
$ \Hom_H(H, H \cotimes H) $ corresponding to the Galois map $ \gamma_r $. The relation 
$$
(\id \otimes \gamma_r) \gamma_r^{12} (\id \otimes \gamma_r^{-1}) = \gamma_r^{12} \gamma_r^{13}
$$
is easily verified. Let us remark that rewriting this equation in the form 
$$
\gamma_r^{23} \gamma_r^{12} = \gamma_r^{12} \gamma_r^{13} \gamma_r^{23}
$$
yields the pentagon equation of the Kac-Takesaki operator \cite{BaSka2}. \\
Consider the special case that the bornological quantum group $ H $ is unital. Then there is a natural 
isomorphism $ \Hom_H(H, V \cotimes H) \cong V \cotimes H $ and a coaction is the same thing as 
a bounded linear map $ \eta: V \rightarrow  V \cotimes H $ such that 
$ (\rho \cotimes \id) \rho = (\id \cotimes \Delta)\rho $. That is, for unital bornological quantum groups
the notion of a coaction is very similar to the concept of a coaction as it is used in the theory 
of Hopf algebras. \\
For later use we give the following definitions. An essential $ H $-module $ P $ over a bornological 
quantum group $ H $ is called projective if 
for every $ H $-linear map $ \pi: V \rightarrow W $ with bounded linear splitting $ \sigma: W \rightarrow V $ 
and every $ H $-linear map $ \xi: P \rightarrow W $ there exists an $ H $-linear map $ \zeta: P \rightarrow V $ 
such that $ \pi \zeta = \xi $. In this case we say that $ P $ satisfies the lifting property for linearly split surjections of 
$ H $-modules. In a completely analogous way one defines the notion of a projective essential $ H $-comodule. \\
We conclude this section by studying the functoriality of essential modules and comodules under morphisms of quantum groups. 
Let $ \alpha: H \rightarrow M(K) $ be a morphism of bornological quantum groups. If $ \lambda: K \cotimes V \rightarrow V $ is an essential $ K $-module 
structure on $ V $ then $ \alpha^*(\lambda) $ is the $ H $-module structure defined by  
$$
 \xymatrix{
     H \cotimes V \ar@{->}[r]^{\!\!\!\!\!\!\!\!\! \alpha \cotimes \id} & M(K) \cotimes V \ar@{->}[r]^{\!\!\!\!\!\!\!\!\!\!\! \cong} & 
     M(K) \cotimes K \cotimes_K V \ar@{->}[r]^{\;\;\mu \cotimes \id} & K \cotimes_K V \cong V 
     }
$$
and it is easy to check that $ V $ becomes an essential $ H $-module in this way. This construction is evidently compatible with 
module maps and thus yields a functor $ \alpha^*: K \LSMod \rightarrow H \LSMod $. A similar functor is obtained for right modules. \\
Conversely, let $ \eta: V \cotimes H \rightarrow V \cotimes H $ be an essential $ H $-comodule. We define a bounded linear map 
$ \alpha_*(\eta): V \cotimes K \rightarrow V \cotimes K $ by the commutative diagram 
$$
\xymatrix{
V \cotimes K \ar@{->}[r]^{\alpha_*(\eta)} \ar@{->}[d]^{\cong} & 
       V \cotimes K \ar@{->}[d]^{\cong} \\
       V \cotimes H \cotimes_H K \ar@{->}[r]^{\eta \cotimes \id} & 
       V \cotimes H \cotimes_H K
     }
$$
where we use that $ \eta $ is right $ H $-linear. It is evident that $ \alpha_*(\eta) $ is a right $ K $-linear isomorphism and 
one checks that the relation 
$$
(\id \otimes \gamma_r) \alpha_*(\eta)_{12} (\id \otimes \gamma_r^{-1}) = \alpha_*(\eta)_{12} \alpha_*(\eta)_{13}
$$ 
is satisfied. Hence $ \alpha_*(\eta) $ defines a coaction of $ K $ on $ V $. 
This construction is compatible with 
comodule maps and yields a functor $ \alpha_*: \SComodR H \rightarrow \SComodR K$. 
Again, there is a similar functor for left comodules. 

\section{The dual quantum group and Pontrjagin duality}\label{secdual}

In this section we construct the dual quantum group $ \hat{H} $ of a bornological quantum group $ H $. 
Moreover we prove the analogue 
of Pontrjagin duality in the context of bornological quantum groups. Unless further specified we assume that $ \phi $ is a left Haar functional 
on $ H $ and we let $ \psi $ be any right Haar functional. \\ 
Using the invariant functional $ \phi $ we define bounded linear maps $ \mathcal{F}_l $ and $ \mathcal{F}_r $ from $ H $ 
into the dual space $ H' = \Hom(H, \mathbb{C}) $ by 
\begin{equation*}
\mathcal{F}_l(f)(h) = \phi(hf), \qquad \mathcal{F}_r(f)(h) = \phi(fh).
\end{equation*}
Similarly, we obtain bounded linear maps $ \mathcal{G}_l $ and $ \mathcal{G}_r $ from $ H $ into $ H' $ by 
\begin{equation*}
\mathcal{G}_l(f)(h) = \psi(hf), \qquad \mathcal{G}_r(f)(h) = \psi(fh)
\end{equation*} 
and all these maps are injective by faithfulness. Using notation and results from section \ref{secmodhaar} we 
obtain the following statement. 
\begin{prop}\label{Fourierprop}
Let $ H $ be a bornological quantum group. Then 
$$
\mathcal{F}_r(f) = \mathcal{F}_l(\sigma(f)), \qquad \mathcal{G}_l(f) = \mathcal{F}_l(\nu(f)), \qquad \mathcal{G}_r(f) = \mathcal{G}_l(\rho(f))
$$
for all $ f \in H $. 
\end{prop} 
Due to proposition \ref{Fourierprop} the images of the maps $ \mathcal{F}_l, \mathcal{F}_r, \mathcal{G}_l, \mathcal{G}_r $ 
in $ H' $ coincide. Let us write $ \hat{H} $ for this space. Moreover, since the maps $ \sigma, \nu $ and $ \rho $ are isomorphisms 
we may use any of them to define a unique bornology 
on $ \hat{H} $ by transferring the bornology from $ H $. We will always view $ \hat{H} $ as a bornological vector space 
with this bornology and hence the maps $ \mathcal{F}_l, \mathcal{F}_r, \mathcal{G}_l, \mathcal{G}_r $ yield bornological isomorphisms
from $ H $ to $ \hat{H} $. In particular, the space $ \hat{H} $ satisfies again the approximation 
property. \\
We say that a bounded bilinear map $ b: U \times V \rightarrow W $ is nondegenerate if 
$ b(u,v) = 0 $ for all $ u \in U $ implies $ v = 0 $ and $ b(u,v) = 0 $ for all $ v \in V $ implies $ u = 0 $. 
Since $ H $ is a regular bornological vector space the canonical pairing between $ H $ and $ H' $ given by $ \bra f, \omega \ket = \omega(f) $ is nondegenerate. 
By construction of the space $ \hat{H} $  there is an obvious injective bounded linear map $ \hat{H} \rightarrow H' $ and we have a 
nondegenerate pairing between $ H $ and $ \hat{H} $ as well. The latter may be extended naturally to a pairing between $ M(H) $ and $ \hat{H} $ 
which is again nondegenerate. There are similar constructions for tensor powers of $ H $ and $ \hat{H} $. \\
In order to obtain a quantum group structure on $ \hat{H} $ our first aim is to define a multiplication. Consider the transpose map 
$ \Delta^*: M(H \cotimes H)' \rightarrow H' $ of the comultiplication given by 
$$
\Delta^*(\omega)(f) = \omega(\Delta(f)). 
$$
According to the previous remarks, $ \hat{H} \cotimes \hat{H} $ can be viewed as a linear subspace of $ M(H \cotimes H) $ and $ \Delta^* $ 
restricts to a map $ \hat{H} \cotimes \hat{H} \rightarrow H' $. We shall show that the latter actually yields a 
bounded linear map $ \hat{H} \cotimes \hat{H} \rightarrow \hat{H} $. 
In order to do this we define a bounded linear map $ m: H \cotimes H \rightarrow H $ by 
$$
m(f \otimes g) = (\id \cotimes \phi)\gamma_l^{-1}(f \otimes g). 
$$
Transferring this map according to the isomorphism $ \mathcal{F}_l: H \rightarrow \hat{H} $ 
we obtain a bounded linear map $ \hat{\mu}: \hat{H} \cotimes \hat{H} \rightarrow \hat{H} $ which we call 
the convolution product. 
One computes 
\begin{align*}
(\id \cotimes \mu)&(\mu \cotimes \id_{(2)})(S \cotimes \id_{(3)})(\id \cotimes \Delta \cotimes \id)(\mu_{(2)} \cotimes \id)(\Delta \cotimes \id_{(3)}) \\
&= (\mu \cotimes \id)(S \cotimes \id_{(2)})(\id \cotimes \mu)(\id \cotimes \Delta \cotimes \id)(\mu_{(2)} \cotimes \id)(\Delta \cotimes \id_{(3)}) \\
&= (\mu \cotimes \id)(S \cotimes \id_{(2)})(\mu \cotimes \id_{(2)})(\tau \cotimes \id_{(2)})(\id \cotimes \gamma_r \cotimes \id) \\
&\qquad (\id \cotimes \tau \cotimes \id)(\id_{(2)} \cotimes \gamma_r)(\id \cotimes \tau \cotimes \id)(\id_{(2)} \cotimes \gamma_r)(\tau \cotimes \id_{(2)}) \allowdisplaybreaks \\
&= (\mu \cotimes \id)(\mu \cotimes \id_{(2)})(S \cotimes S \cotimes \id_{(2)})(\id \cotimes \gamma_r \cotimes \id) \\
&\qquad (\id \cotimes \tau \cotimes \id)(\id_{(2)} \cotimes \gamma_r)(\id \cotimes \tau \cotimes \id)(\id_{(2)} \cotimes \gamma_r)(\tau \cotimes \id_{(2)}) \allowdisplaybreaks \\
&= (\mu \cotimes \id)(S \cotimes \id_{(2)})(\id \cotimes \mu \cotimes \id)(\id \cotimes S \cotimes \id_{(2)})(\id \cotimes \gamma_r \cotimes \id) \\
&\qquad (\id \cotimes \tau \cotimes \id)(\id_{(2)} \cotimes \gamma_r)(\id \cotimes \tau \cotimes \id)(\id_{(2)} \cotimes \gamma_r)(\tau \cotimes \id_{(2)}) \allowdisplaybreaks \\
&= (\mu \cotimes \id)(S \cotimes \id_{(2)})(\id \cotimes \epsilon \cotimes \id_{(2)})(\id \cotimes \tau \cotimes \id)\\
&\qquad (\id_{(2)} \cotimes \gamma_r)(\id \cotimes \tau \cotimes \id)(\id_{(2)} \cotimes \gamma_r)(\tau \cotimes \id_{(2)}) \allowdisplaybreaks \\
&= (\mu \cotimes \id)(S \cotimes \id_{(2)})(\id_{(2)} \cotimes \epsilon \cotimes \id)\\
&\qquad (\id_{(2)} \cotimes \gamma_r)(\id \cotimes \tau \cotimes \id) (\id_{(2)} \cotimes \gamma_r)(\tau \cotimes \id_{(2)}) \\
&= (\mu \cotimes \id)(\id_{(2)} \cotimes \mu)(S \cotimes \id_{(3)})(\id \cotimes \tau \cotimes \id)(\id_{(2)} \cotimes \gamma_r)(\tau \cotimes \id_{(2)}) \\
&= (\id \cotimes \mu)(\mu \cotimes \id_{(2)})(S \cotimes \id_{(3)})(\id \cotimes \rho_r \cotimes \id)(\tau \cotimes \id_{(2)}). 
\end{align*}
We conclude 
\begin{equation}
(\mu \cotimes \id)(S \cotimes \id_{(2)})(\id \cotimes \Delta)\mu_{(2)}(\Delta \cotimes \id_{(2)}) = 
(\mu \cotimes \id)(S \cotimes \id_{(2)})(\id \cotimes \rho_r)(\tau \cotimes \id). 
\end{equation}
Moreover we have
\begin{align*}
\tau(\mu \cotimes \id)(&\id \cotimes \gamma_l^{-1}) \\
&= \tau(\mu \cotimes \id)(\id_{(2)} \cotimes \mu)(\id_{(2)} \cotimes S^{-1}\cotimes \id)
(\id \cotimes \tau \cotimes \id)(\id \cotimes \Delta \cotimes \id)(\id \cotimes \tau) \\
&= \tau (\id \cotimes \mu)(\id \cotimes S^{-1} \cotimes \id)(\tau \cotimes \id)(\rho_r \cotimes \id)(\id \cotimes \tau) \\
&= (\mu \cotimes \id)(\tau \cotimes \id)(\id \cotimes S^{-1} \cotimes \id)(\id \cotimes \rho_r)(\tau \cotimes \id). 
\end{align*}
Hence we compute using invariance 
\begin{align*}
\hat{\mu}&(\mathcal{F}_l(f) \otimes \mathcal{F}_l(g))(h) = \mathcal{F}_l(m(f \otimes g))(h) = 
(\phi \cotimes \phi)(\mu \cotimes \id)(\id \cotimes \gamma_l^{-1})(h \otimes f \otimes g) \\
&= (\phi \cotimes \phi)(\mu \cotimes \id)(\tau \cotimes \id)(\id \cotimes S^{-1} \cotimes \id)(\id \cotimes \rho_r)(\tau \cotimes \id)(h \otimes f \otimes g) \\
&= (\phi \cotimes \phi)(S^{-1} \cotimes \id)(\mu \cotimes \id)(S \cotimes \id_{(2)})(\id \cotimes \rho_r)(\tau \cotimes \id)(h \otimes f \otimes g) \\
&= (\phi \cotimes \phi) (S^{-1} \cotimes \id)(\mu \cotimes \id)(S \cotimes \id_{(2)})
(\id \cotimes \Delta) \mu_{(2)} (\Delta \cotimes \id_{(2)})(h \otimes f \otimes g) \\
&= (\phi \cotimes \phi) \mu_{(2)} (\Delta \cotimes \id_{(2)})(h \otimes f \otimes g) \\
&= \Delta^*(\mathcal{F}_l(f) \otimes \mathcal{F}_l(g))(h)
\end{align*}
and deduce that $ \hat{\mu} $ can be identified with $ \Delta^* $. Using this statement one calculates 
\begin{align*}
\hat{\mu}(\hat{\mu} \cotimes \id) (\mathcal{F}_l(f) & \otimes \mathcal{F}_l(g) \otimes \mathcal{F}_l(h))(x) = 
(\hat{\mu}(\mathcal{F}_l(f) \otimes \mathcal{F}_l(g)) \cotimes \phi)(\Delta(x)(1 \otimes h)) \\
&= (\mathcal{F}_l(f) \cotimes \mathcal{F}_l(g) \cotimes \phi)((\Delta \cotimes \id)\Delta(x)(1 \otimes 1 \otimes h)) \\
&= (\phi \cotimes \phi \cotimes \phi)((\Delta \cotimes \id)\Delta(x)(f \otimes g \otimes h)) \\
&= (\phi \cotimes \phi \cotimes \phi)((\id \cotimes \Delta)\Delta(x)(f \otimes g \otimes h)) \\
&= (\phi \cotimes \mathcal{F}_l(g) \cotimes \mathcal{F}_l(h))((\id \cotimes \Delta)\Delta(x)(f \otimes 1 \otimes 1)) \\
&= (\phi \cotimes \hat{\mu}(\mathcal{F}_l(g) \otimes \mathcal{F}_l(h)))(\Delta(x)(f \otimes 1)) \\
&= \hat{\mu}(\id \cotimes \hat{\mu}) (\mathcal{F}_l(f) \otimes \mathcal{F}_l(g) \otimes \mathcal{F}_l(h))(x)
\end{align*}
which means that the convolution product $ \hat{\mu} $ is 
associative. Hence $ \hat{H} $ is a bornological algebra with convolution as 
multiplication. According to the above considerations we have 
\begin{equation}\label{hatmufl}
\hat{\mu}(\mathcal{F}_l(f) \otimes \mathcal{F}_l(g)) = \mathcal{F}_l((\id \cotimes \phi)\gamma_l^{-1}(f \otimes g))
\end{equation}
and an analogous calculation yields the formula
\begin{equation}\label{hatmugr}
\hat{\mu}(\mathcal{G}_r(f) \otimes \mathcal{G}_r(g)) = \mathcal{G}_r((\psi \cotimes \id)\rho_r^{-1}(f \otimes g))
\end{equation}
for the multiplication in $ \hat{H} $. Actually, equation (\ref{hatmugr}) may be obtained directly from the previous discussion 
applied to $ H^{opcop} $. \\
For later use we shall extend the multiplication of $ \hat{H} $ in the following way. According to equation (\ref{coop12}) and 
the fact that $ \phi $ is left invariant we have 
\begin{equation}\label{phiinvarproj}
(\id \cotimes \phi)\gamma_r = (\id \cotimes \phi)\gamma_l \tau (S^{-1} \cotimes \id)\rho_r = (\id \cotimes \phi)(S^{-1} \cotimes \id)\rho_r. 
\end{equation}
Using this observation we define a bounded linear map 
$ \hat{\mu}_l: H' \cotimes \hat{H} \rightarrow H' $  by 
\begin{equation}\label{hatmul}
\hat{\mu}_l(\omega \otimes \mathcal{F}_l(f))(x) = (\omega \cotimes \phi)\gamma_r(x \otimes f) = (\omega \cotimes \phi)(S^{-1} \cotimes \id)\rho_r(x \otimes f). 
\end{equation}
Inserting $ \omega = \epsilon $ we see that $ \hat{\mu}_l(\omega \otimes \mathcal{F}_l(f)) = 0 $ 
for all $ \omega $ implies $ f = 0 $. Conversely, assume $ \hat{\mu}_l(\omega \otimes \mathcal{F}_l(f)) = 0 $ 
for all $ f \in H $. Then we have $ (\omega \cotimes \phi)\gamma_r(h \otimes f) = 0 $ for all $ h, f \in H $ 
and since $ \gamma_r $ is an isomorphism this yields $ \omega = 0 $. Hence $ \hat{\mu}_l $ defines a nondegenerate 
pairing. \\
Similarly, according to equation (\ref{coop11}) we have 
\begin{equation}\label{psiinvarproj}
(\psi \cotimes \id)\rho_l = (\psi \cotimes \id)(\id \cotimes S^{-1}) \gamma_r^{-1} \gamma_l = (\psi \cotimes \id)(\id \cotimes S^{-1})\gamma_l
\end{equation}
and we define $ \hat{\mu}_r: \hat{H} \cotimes H' \rightarrow H' $ by 
\begin{equation}\label{hatmur}
\hat{\mu}_r(\mathcal{G}_r(f) \otimes \omega)(x) = (\psi \cotimes \omega) \rho_l(f \otimes x) = (\psi \cotimes \omega) (\id \cotimes S^{-1})\gamma_l(f \otimes x). 
\end{equation}
As above one sees that the pairing given by $ \hat{\mu}_r $ is nondegenerate. If restricted to $ \hat{H} \cotimes \hat{H} $ the maps 
$ \hat{\mu}_l $ and $ \hat{\mu}_r $ are equal to the multiplication map $ \hat{\mu} $. Moreover it is straightforward to 
check that the maps $ \hat{\mu}_l $ and $ \hat{\mu}_r $ are associative whenever this assertion makes sense.
In the sequel we will simply write $ \hat{\mu} $ for the maps $ \hat{\mu}_l $ and $ \hat{\mu}_r $, respectively. \\
Using the definition of the modular automorphism $ \sigma $ we obtain
$$
(\id \cotimes \phi)\rho_r(f \otimes x) = (\id \cotimes \phi)\gamma_r(x \otimes \sigma(f))
$$
and together with equation (\ref{phiinvarproj}) this yields the formula
\begin{equation}\label{hatmulalt}
\hat{\mu}(\omega \otimes \mathcal{F}_r(f))(x) = (\omega \cotimes \phi) \rho_r(f \otimes x) = (\omega \cotimes \phi)(S \cotimes \id)\gamma_r(f \otimes x) 
\end{equation}
for the multiplication $ \hat{\mu} $. In a similar way we obtain
\begin{equation}\label{hatmuralt}
\hat{\mu}(\mathcal{G}_l(f) \otimes \omega)(x) = (\psi \cotimes \omega)\gamma_l(x \otimes f) = (\psi \cotimes \omega) (\id \cotimes S)\rho_l(x \otimes f)
\end{equation}
using equation (\ref{psiinvarproj}). \\
Our next aim is to show that $ \hat{H} $ is a projective module over itself. In order to do this we study the regular coaction of $ H $ on itself given 
by $ \gamma_r $. 
\begin{prop}\label{Hproj}
The regular coaction of a bornological quantum group $ H $ on itself is a 
projective $ H $-comodule. 
\end{prop}
\proof Choose an element $ h \in H $ such that $ \phi(h) = 1 $ and define $ \nu: H \rightarrow H \cotimes H $ 
by 
$$
\nu(f) = \gamma_l(h \otimes f). 
$$
Then we have 
$$
(\id \cotimes \phi) \nu(f) = (\id \cotimes \phi)\gamma_l(h \otimes f) = \phi(h) f = f
$$
for all $ f \in H $ since $ \phi $ is left invariant. 
Hence the map $ \nu $ satisfies the equation 
\begin{equation} \label{nuinv}
(\id \cotimes \phi) \nu  = \id.
\end{equation} 
Let us moreover define $ \lambda: H \rightarrow H \cotimes H $ by $ \lambda = \gamma_r^{-1} \nu $. \\
Now assume that $ \pi: V \rightarrow W $ is a surjective map of $ H $-comodules with bounded linear splitting $ \sigma $ and let 
$ \xi: H \rightarrow W $ be an $ H $-colinear map. We define $ \zeta: H \rightarrow V $ 
as the composition 
\begin{equation*} 
 \xymatrix{
     H \; \ar@{->}[r]^{\!\!\!\lambda} & H \cotimes H \ar@{->}[r]^{\sigma \xi \cotimes \id} & V \cotimes H \ar@{->}[r]^{\; \eta_V} 
     & V \cotimes H \ar@{->}[r]^{\;\;\;\id \cotimes \phi} & V 
     }
\end{equation*}
where $ \eta_V $ is the coaction of $ V $. \\
Let us check that $ \zeta $ is $ H $-colinear. Using equation (\ref{coop11}) we obtain 
\begin{align*}
(\gamma_r^{-1} \cotimes \id)(&\gamma_l \cotimes \id)(\id \cotimes \gamma_r) = (\id \cotimes S \cotimes \id)(\rho_l \cotimes \id)(\id \cotimes \gamma_r) \\
&= (\id \cotimes S \cotimes \id)(\id \cotimes \gamma_r)(\rho_l \cotimes \id) \\
&= (\id \cotimes S \cotimes \id)(\id \cotimes \gamma_r)(\id \cotimes S^{-1} \cotimes \id)(\gamma_r^{-1} \cotimes \id)(\gamma_l \cotimes \id)
\end{align*}
and deduce
\begin{equation}\label{rhogammproj}
(\lambda \cotimes \id)\gamma_r = (\id \cotimes S \cotimes \id) (\id \cotimes \gamma_r) (\id \cotimes S^{-1} \cotimes \id) (\lambda \cotimes \id). 
\end{equation}
Since $ S $ is an algebra and coalgebra antihomomorphism we have 
$$
\tau(S \cotimes S)\gamma_r = \rho_l (S \cotimes S) \tau
$$
which yields 
\begin{equation}\label{gammainvalt}
(S \cotimes \id)\gamma_r(S^{-1} \cotimes \id) = \tau (S^{-1} \cotimes \id)\rho_l(S \cotimes \id)\tau.
\end{equation}
Using equation (\ref{coop12}) and equation (\ref{coop2}) we obtain 
\begin{equation}\label{rhogammhelp}
\rho_r = (S \cotimes \id)\tau \gamma_l^{-1} \gamma_r = \rho_l (S \cotimes \id) \gamma_r. 
\end{equation}
Since $ \eta_V $ is right $ H $-linear we calculate 
\begin{align*}
\eta_V^{13}(\id \cotimes \rho_r)&(\id \cotimes \mu \cotimes \id) = \eta_V^{13}(\id_{(2)} \cotimes \mu)(\id \cotimes \tau \cotimes \id)(\id_{(2)} \cotimes \rho_r) \\ &= (\id_{(2)} \cotimes \mu)\eta_V^{13}(\id \cotimes \tau \cotimes \id)(\id_{(2)} \cotimes \rho_r) \\
&= (\id_{(2)} \cotimes \mu)(\id \cotimes \tau \cotimes \id)\eta_V^{12}(\id_{(2)} \cotimes \rho_r) \\
&= (\id_{(2)} \cotimes \mu)(\id \cotimes \tau \cotimes \id)(\id_{(2)} \cotimes \rho_r)\eta_V^{12} \\
&= (\id \cotimes \rho_r)(\id \cotimes \mu \cotimes \id)\eta_V^{12} \\
&= (\id \cotimes \rho_r) \eta_V^{12} (\id \cotimes \mu \cotimes \id)
\end{align*}
which yields
$$
\eta_V^{13}(\id \cotimes \rho_r) = (\id \cotimes \rho_r) \eta_V^{12}
$$
since $ H $ is essential. Combining this with equation (\ref{rhogammhelp}) implies 
\begin{equation}\label{etalinearr}
(\id \cotimes \rho_r)\eta_V^{12}(\id \cotimes \gamma_r^{-1}) 
= \eta_V^{13}(\id \cotimes \rho_l)(\id \cotimes S \cotimes \id). 
\end{equation}
According to equations (\ref{rhogammproj}),(\ref{gammainvalt}), (\ref{etalinearr}) and (\ref{phiinvarproj}) we get 
\begin{align*}
(\zeta &\cotimes \id)\gamma_r = (\id \cotimes \phi \cotimes \id) (\eta_V \cotimes \id)(\sigma \xi \cotimes \id_{(2)})(\lambda \cotimes \id) \gamma_r \\
&= (\id \cotimes \phi \cotimes \id)(\eta_V \cotimes \id)(\id \otimes S \cotimes \id)(\id \cotimes \gamma_r)(\id \cotimes S^{-1} \cotimes \id)
(\sigma \xi \cotimes \id_{(2)})(\lambda \cotimes \id) \\
&= (\id \cotimes \id \cotimes \phi)\eta_V^{13}(\id \otimes S^{-1} \cotimes \id)(\id \cotimes \rho_l)(\id \cotimes S \cotimes \id)(\id \cotimes \tau) 
(\sigma \xi \cotimes \id_{(2)})(\lambda \cotimes \id) \\
&= (\id \cotimes \id \cotimes \phi)(\id \cotimes S^{-1} \cotimes \id)(\id \cotimes \rho_r)\eta^{12}_V(\id \cotimes \gamma_r^{-1}) 
(\id \cotimes \tau)(\sigma \xi \cotimes \id_{(2)})(\lambda \cotimes \id) \\
&= (\id \cotimes \id \cotimes \phi)(\id \cotimes \gamma_r)\eta^{12}_V(\id \cotimes \gamma_r^{-1}) (\id \cotimes \tau)
(\sigma \xi \cotimes \id_{(2)})(\lambda \cotimes \id) \\
&= (\id \cotimes \id \cotimes \phi)\eta^{12}_V \eta^{13}_V (\id \cotimes \tau)(\sigma \xi \cotimes \id_{(2)})(\lambda \cotimes \id) \\
&= \eta_V(\id \cotimes \phi \cotimes \id) (\eta_V \cotimes \id)(\sigma \xi \cotimes \id_{(2)})(\lambda \cotimes \id) \\
&= \eta_V (\zeta \cotimes \id)
\end{align*}
which shows that $ \zeta $ is $ H $-colinear. \\
Let us now prove that $ \zeta $ is a lifting for $ \xi $. Due to the fact that $ \pi $ is colinear the diagram
\begin{equation*} 
 \xymatrix{
     H \cotimes H \;\ar@{->}[r]^{\sigma \xi \cotimes \id} \ar@{->}[d]^\id & V \cotimes H \ar@{->}[r]^{\;\;\; \eta_V} \ar@{->}[d]^{\pi \cotimes \id}
     & V \cotimes H \ar@{->}[r]^{\;\;\;\id \cotimes \phi} \ar@{->}[d]^{\pi \cotimes \id} & V \ar@{->}[d]^\pi\\
     H \cotimes H \;\ar@{->}[r]^{\xi\cotimes \id} & W \cotimes H \ar@{->}[r]^{\;\;\; \eta_W} 
     & W \cotimes H \ar@{->}[r]^{\;\;\;\id \cotimes \phi} & W
     }
\end{equation*}
is commutative and moreover we have 
$$
\eta_W (\xi \cotimes \id) = (\xi \cotimes \id) \gamma_r 
$$
because $ \xi $ is colinear. As a consequence we get $ \pi \zeta = \xi $ since
$$
(\id \cotimes \phi)\gamma_r \lambda = (\id \cotimes \phi)\nu = \id
$$
by the definition of $ \lambda $. \qed \\
Consider the transposed right regular coaction $ \rho $ of $ H $ on itself 
given by $ \rho = \gamma_l^{-1} \tau $. Since equation (\ref{coop1}) for $ H^{op} $ yields
\begin{equation} \label{gammartogamml}
\gamma_r(S^{-1} \cotimes \id) = (S^{-1} \cotimes \id)\gamma_l^{-1} \tau
\end{equation}
we see that the transposed right regular coaction $ \rho $ corresponds to the right regular coaction $ \gamma_r $ 
under the linear automorphism of $ H $ given by $ S^{-1} $. It follows that the map $ \rho $ is indeed a 
coaction and that the coactions $ \rho $ and $ \gamma_r $ yield isomorphic comodules. In particular the comodule 
defined by $ \rho $ is projective due to proposition \ref{Hproj}. \\
As above let $ m: H \cotimes H \rightarrow H $ denote the map corresponding to the multiplication of $ \hat{H} $ under 
the isomorphism $ \mathcal{F}_l $. 
By definition of the right regular coaction we get
\begin{equation}\label{muhatrho}
m = (\id \cotimes \phi)\rho \tau.
\end{equation}
The pentagon relation for the operator $ \gamma_r^{cop} = \tau \gamma_l = \rho^{-1} $ can be written as 
\begin{equation}\label{pentrho}
\rho^{23} \rho^{13} \rho^{12} = \rho^{12} \rho^{23}
\end{equation}
and together with the formula $ (\id \cotimes \phi) \tau \rho^{-1} = (\id \otimes \phi)\gamma_l = (\id \cotimes \phi)\tau $ 
this shows
\begin{align*}
(m \cotimes \id)&\rho^{13} = (\id \cotimes \id \cotimes \phi)(\id \cotimes \tau)(\rho \tau \cotimes \id) \rho^{13} \\
&= (\id \cotimes \id \cotimes \phi)(\id \cotimes \tau \rho^{-1}) \rho^{12}(\tau \cotimes \id) \rho^{13} \\
&= (\id \cotimes \id \cotimes \phi)(\id \cotimes \tau)(\id \cotimes \rho^{-1}) \rho^{12} \rho^{23} (\tau \cotimes \id) \\
&= (\id \cotimes \id \cotimes \phi)(\id \cotimes \tau)\rho^{13}\rho^{12}(\tau \cotimes \id) \\
&= (\id \cotimes \id \cotimes \phi)(\rho \cotimes \id)(\id \cotimes \tau) \rho^{12}(\tau \cotimes \id) \\
&= \rho(m \otimes \id)
\end{align*}
which means that the map $ m $ is right $ H $-colinear if we view 
$ H \cotimes H $ as a right $ H $-comodule using the coaction $ \rho^{13} $. Since $ m $ has a bounded linear splitting we obtain a 
colinear splitting $ \sigma: H \rightarrow H \cotimes H $ due to proposition \ref{Hproj}. That is, we have 
$$
(\sigma \cotimes \id)\rho = \rho^{13}(\sigma \cotimes \id)
$$ 
which yields 
$$
\sigma m = (\id_{(2)} \cotimes \phi)(\sigma \cotimes \id)\rho \tau = (\id \cotimes \phi \cotimes \id)(\rho \tau \cotimes \id)(\id \cotimes \sigma) = 
(m \cotimes \id)(\id \cotimes \sigma).
$$
Translating this to $ \hat{H} $ using the isomorphism $ \mathcal{F}_l $ we see that there is a $ \hat{H} $-linear splitting for the 
multiplication map $ \hat{\mu} $ if $ \hat{H} $ acts by multiplication on the left tensor factor of $ \hat{H}\cotimes \hat{H} $. Using such 
a splitting it is straightforward to check that $ \hat{H} $ is an essential bornological algebra. \\
We define a linear form $ \hat{\psi} $ on $ \hat{H} $ by 
$$
\hat{\psi}(\mathcal{F}_l(f)) = \epsilon(f)
$$
and compute 
$$
\hat{\psi}(\mathcal{F}_l(f)\mathcal{F}_l(g)) = (\epsilon \cotimes \phi)\gamma_l^{-1}(f \otimes g) = \phi(S^{-1}(g) f)
$$
for all $ f, g \in H $ which implies that $ \hat{\psi} $ is faithful since $ \phi $ is faithful. Hence the algebra $ \hat{H} $ is 
equipped with a faithful bounded linear functional. We will see below that $ \hat{\psi} $ is 
right invariant for the comultiplication of $ \hat{H} $, however, first we have to construct this comultiplication of course. \\
In order to do this we define a bounded linear map $ \hat{\gamma}_r: \hat{H}\cotimes \hat{H} \rightarrow \hat{H} \cotimes \hat{H} $ 
by
$$
\hat{\gamma}_r(\mathcal{F}_l(f) \otimes \mathcal{F}_l(g)) = (\mathcal{F}_l \cotimes \mathcal{F}_l)\tau \gamma_l^{-1}(f \otimes g).
$$
It is evident that $ \hat{\gamma}_r $ is a bornological isomorphism. 
Let us show that $ \hat{\gamma}_r $ commutes with right multiplication on the second tensor factor. 
Using the pentagon relation (\ref{pentrho}) for $ \rho $ and $ (\phi \cotimes \id)\rho = \phi \cotimes \id $ we have 
\begin{align*}
(\id \cotimes &\id \cotimes \phi) (\tau \cotimes \id)(\gamma_l^{-1} \cotimes \id)(\id \cotimes \gamma_l^{-1})(f \otimes g \otimes h) \\
&= (\id \cotimes \phi \cotimes \id)\rho^{13} \rho^{12}(h \otimes g \otimes f) \\
&= (\id \cotimes \phi \cotimes \id)\rho^{23}\rho^{13} \rho^{12}(h \otimes g \otimes f) \\
&= (\id \cotimes \phi \cotimes \id)\rho^{12} \rho^{23}(h \otimes g \otimes f) \\
&= (\id \cotimes \id \cotimes \phi) (\id \cotimes \gamma_l^{-1})(\tau \cotimes \id)(\gamma_l^{-1} \cotimes \id)(f \otimes g \otimes h)
\end{align*}
and translating this using the map $ \mathcal{F}_l $ we obtain 
\begin{equation}
\hat{\gamma}_r(\id \cotimes \hat{\mu}) = (\id \cotimes \hat{\mu}) (\hat{\gamma}_r \cotimes \id)
\end{equation}
as desired. Similarly, we define a bornological automorphism $ \hat{\rho}_l $ of $ \hat{H}\cotimes \hat{H} $ by 
$$
\hat{\rho}_l(\mathcal{G}_r(g) \otimes \mathcal{G}_r(f)) = (\mathcal{G}_r \cotimes \mathcal{G}_r)\tau \rho_r^{-1}(g \otimes f). 
$$
and using formula (\ref{hatmugr}) we obtain
\begin{equation}
\hat{\rho}_l(\hat{\mu} \cotimes \id) = (\hat{\mu} \cotimes \id) (\id \cotimes \hat{\rho}_l).
\end{equation}
Due to equation (\ref{phiinvarproj}) and equation (\ref{psiinvarproj}) we have 
\begin{equation}\label{gammahatrev}
\hat{\gamma}_r (\omega \otimes \mathcal{F}_l(g))(x \otimes y) = (\omega \cotimes \phi)(\mu \cotimes \id) 
(\id \cotimes S^{-1} \cotimes \id)(\id \cotimes \rho_r)(x \otimes y \otimes g)
\end{equation}
and 
\begin{equation}\label{rhohatlev}
\hat{\rho}_l (\mathcal{G}_r(f) \otimes \omega)(x \otimes y) = (\psi \cotimes \omega)(\id \cotimes \mu) 
(\id \cotimes S^{-1} \cotimes \id)(\gamma_l \cotimes \id)(f \otimes x \otimes y)
\end{equation}
for all $ \omega \in \hat{H} $. Using equation (\ref{hatmul}) and equation (\ref{hatmur}) we calculate 
\begin{align*}
(\id \cotimes &\hat{\mu})(\hat{\rho}_l \cotimes \id)(\mathcal{G}_r(f) \otimes \omega \otimes \mathcal{F}_l(g))(x \otimes y) \\
&= (\psi \cotimes \omega \cotimes \phi)(\id \cotimes \mu \cotimes \id)(\id \cotimes S^{-1} \cotimes \id_{(2)})(\gamma_l \cotimes \id_{(2)}) \\
&\qquad (\id_{(2)} \cotimes S^{-1} \cotimes \id)(\id_{(2)} \cotimes \rho_r)(f \otimes x \otimes y \otimes g) \\
&= (\psi \cotimes \omega \cotimes \phi)(\id \cotimes \mu \cotimes \id)(\id_{(2)} \cotimes S^{-1} \cotimes \id)(\id_{(2)} \cotimes \rho_r) \\
&\qquad (\id \cotimes S^{-1} \cotimes \id_{(2)})(\gamma_l \cotimes \id_{(2)})(f \otimes x \otimes y \otimes g) \\
&= (\hat{\mu} \cotimes \id)(\id \cotimes \hat{\gamma}_r)(\mathcal{G}_r(f) \otimes \omega \otimes \mathcal{F}_l(g))(x \otimes y)
\end{align*}
which yields 
\begin{equation}
(\id \cotimes \hat{\mu})(\hat{\rho}_l \cotimes \id) = (\hat{\mu} \cotimes \id)(\id \cotimes \hat{\gamma}_r).
\end{equation}
Using the pentagon relation (\ref{pentrho}) for $ \rho $ we compute 
\begin{align*}
\hat{\gamma}_r(\hat{\mu} &\cotimes \id)(\mathcal{F}_l(f) \otimes \mathcal{F}_l(g) \otimes \mathcal{F}_l(h)) \\
&= (\mathcal{F}_l \cotimes \mathcal{F}_l)(\phi \cotimes \id \cotimes \id)(\id \cotimes \tau \gamma_l^{-1})(\tau \gamma_l^{-1} \cotimes \id)(f \otimes g \otimes h) \\
&= (\mathcal{F}_l \cotimes \mathcal{F}_l)(\id \cotimes \id \cotimes \phi)\rho^{12} \rho^{23}(h \otimes g \otimes f) \allowdisplaybreaks \\
&= (\mathcal{F}_l \cotimes \mathcal{F}_l)(\id \cotimes \id \cotimes \phi)\rho^{23} \rho^{13} \rho^{12}(h \otimes g \otimes f) \allowdisplaybreaks \\
&= (\mathcal{F}_l \cotimes \mathcal{F}_l)(\phi \cotimes \id \cotimes \id)(\tau \gamma_l^{-1} \cotimes \id)(\tau \cotimes \id)(\id \cotimes \tau\gamma_l^{-1}) \\
&\qquad (\tau \cotimes \id)(\id \cotimes \tau \gamma_l^{-1})(f \otimes g \otimes h) \\
&= (\hat{\mu} \cotimes \id)(\tau \cotimes \id)(\id \cotimes \hat{\gamma}_r)(\tau \cotimes \id)(\id \cotimes \hat{\gamma}_r)
(\mathcal{F}_l(f) \otimes \mathcal{F}_l(g) \otimes \mathcal{F}_l(h))
\end{align*}
which yields
\begin{equation}
\hat{\gamma}_r(\hat{\mu} \cotimes \id) = (\hat{\mu} \cotimes \id)(\tau \cotimes \id)(\id \cotimes \hat{\gamma}_r)(\tau \cotimes \id)(\id \cotimes \hat{\gamma}_r). 
\end{equation}
A similar computation shows 
\begin{equation}
\hat{\rho}_l(\id \cotimes \hat{\mu}) = (\id \cotimes \hat{\mu})(\id \cotimes \tau)(\hat{\rho}_l \cotimes \id)(\id \cotimes \tau)(\hat{\rho}_l \cotimes \id). 
\end{equation}
Finally, we have 
\begin{align*}
(\hat{\rho}_l &\cotimes \id)(\id \cotimes \hat{\gamma}_r)(\mathcal{G}_r(f) \otimes \omega \otimes \mathcal{F}_l(g))(x \otimes y \otimes z) \\
&= (\psi \cotimes \omega \cotimes \phi)(\id \cotimes \mu \cotimes \id)(\id_{(2)} \cotimes S^{-1} \cotimes \id)(\id_{(2)} \cotimes \rho_r) \\
&\qquad (\id \cotimes \mu \cotimes \id_{(2)})(\id \cotimes S^{-1} \cotimes \id_{(3)})(\gamma_l \cotimes \id_{(3)})(f \otimes x \otimes y \otimes z \otimes g) \\
&= (\psi \cotimes \omega \cotimes \phi)(\id \cotimes \mu \cotimes \id)(\id \cotimes S^{-1} \cotimes \id_{(2)})(\gamma_l \cotimes \id_{(2)}) \\
&\qquad (\id_{(2)} \cotimes \mu \cotimes \id)(\id_{(3)} \cotimes S^{-1} \cotimes \id)(\id_{(3)} \cotimes \rho_r)(f \otimes x \otimes y \otimes z \otimes g) \\
&= (\id \cotimes \hat{\gamma}_r)(\hat{\rho}_l \cotimes \id)(\mathcal{G}_r(f) \otimes \omega \otimes \mathcal{F}_l(g))(x \otimes y \otimes z) 
\end{align*}
according to equation (\ref{gammahatrev}) and equation (\ref{rhohatlev}) and hence
\begin{equation}
(\hat{\rho}_l \cotimes \id)(\id \cotimes \hat{\gamma}_r) = (\id \cotimes \hat{\gamma}_r)(\hat{\rho}_l \cotimes \id). 
\end{equation}
Using the properties of the maps $ \hat{\gamma}_r $ and $ \hat{\rho}_l $ obtained so far we shall construct the comultiplication for 
$ \hat{H} $ according to the following general result. 
\begin{prop} \label{galoistocomult}
Let $ K $ be an essential bornological algebra satisfying the approximation property equipped with a faithful bounded linear functional. 
If $ \gamma_r $ and $ \rho_l $ are bornological automorphisms of $ K \cotimes K $ such that 
\begin{bnum}
\item[a)] $ \gamma_r(\id \cotimes \mu) = (\id \cotimes \mu)(\gamma_r \cotimes \id) $ 
\item[b)] $ \gamma_r(\mu \cotimes \id) = (\mu \cotimes \id)(\tau \cotimes \id)(\id \cotimes \gamma_r)(\tau \cotimes \id)(\id \cotimes \gamma_r) $
\item[c)] $ \rho_l(\mu \cotimes \id) = (\mu \cotimes \id) (\id \cotimes \rho_l) $ 
\item[d)] $ \rho_l(\id \cotimes \mu) = (\id \cotimes \mu)(\id \cotimes \tau)(\rho_l \cotimes \id)(\id \cotimes \tau)(\rho_l \cotimes \id) $
\item[e)] $ (\id \cotimes \mu)(\rho_l \cotimes \id) = (\mu \cotimes \id)(\id \cotimes \gamma_r) $
\item[f)] $ (\rho_l \cotimes \id)(\id \cotimes \gamma_r) = (\id \cotimes \gamma_r)(\rho_l \cotimes \id) $
\end{bnum}
then there exists a unique comultiplication $ \Delta: K \rightarrow M(K \cotimes K) $ such that $ \gamma_r $ and 
$ \rho_l $ are the associated Galois maps. \\
In addition, if there exist bornological automorphisms $ \gamma_l $ and $ \rho_r $ of $ K \cotimes K $ such that 
\begin{bnum}
\item[g)] $ (\mu \cotimes \id)(\id \cotimes \tau) (\gamma_l \cotimes \id) = \gamma_l (\id \cotimes \mu) $ 
\item[h)] $ (\id \cotimes \mu) (\gamma_l \cotimes \id) = (\mu \cotimes \id)(\id \cotimes \tau)(\gamma_r \cotimes \id)(\id \cotimes \tau) $ 
\item[i)] $ (\id \cotimes \mu)(\tau \cotimes \id) (\id \cotimes \rho_r) = \rho_r (\mu \cotimes \id) $ 
\item[j)] $ (\mu \cotimes \id) (\id \cotimes \rho_r) = (\id \cotimes \mu)(\tau \cotimes \id)(\id \cotimes \rho_l)(\tau \cotimes \id) $
\end{bnum}
then these maps are the remaining Galois maps. In particular all Galois maps are isomorphisms in this case. 
\end{prop}
\proof Using condition a) it is straightforward to check that  
$$
\mu_{(2)}(\Delta_l \cotimes \id_{(2)}) = (\mu \cotimes \id)(\id \cotimes \tau)(\gamma_r \cotimes \id)(\id \cotimes \tau)
$$
defines a bounded linear map $ \Delta_l: K \rightarrow M_l(K \cotimes K) $. According to condition b) the map 
$ \Delta_l $ is actually a homomorphism. Similarly, 
$$
\mu_{(2)}(\id_{(2)} \cotimes \Delta_r) = (\id \cotimes \mu)(\tau \cotimes \id)(\id \cotimes \rho_l)(\tau \cotimes \id)
$$
defines a homomorphism $ \Delta_r: K \rightarrow M_r(K \cotimes K) $ due to conditions c) and d). Condition e) ensures that 
these maps combine to an algebra homomorphism $ \Delta: K \rightarrow M(K \cotimes K) $. It is straightforward to show that $ \gamma_r $ 
and $ \rho_l $ are the corresponding Galois maps. Moreover $ \Delta $ is uniquely determined by these maps. \\
We have to prove that the homomorphism $ \Delta $ is essential. Let us show that the natural map 
$ K \cotimes_K (K \cotimes K) \rightarrow K \cotimes K $ is 
an isomorphism where the module structure 
on $ K \cotimes K $ is given by $ \Delta $. Since $ K $ is essential we have $ K \cotimes_K K \cong K $ and we can 
identify the source of the previous map with $ K \cotimes_K(K \cotimes (K \cotimes_K K)) $ 
in a natural way. It is easy to check that 
$ \gamma_r^{13} $ descends to a bounded 
linear map 
$$ 
\xi: K \cotimes_K(K \cotimes (K \cotimes_K K)) 
\rightarrow (K \cotimes_K K) \cotimes (K \cotimes_K K) 
$$ 
and the composition of $ \xi $ with $ \mu \cotimes \mu $ can be identified with the map we are interested in. 
Hence it suffices to show that $ \xi $ is an isomorphism. Consider the maps 
$ p $ and $ q $ defined on the six-fold tensor product of $ K $ with itself by 
$$
p = \mu \cotimes \id \cotimes \mu \cotimes \id - (\id \cotimes \mu \cotimes \id_{(2)}) \gamma_r^{24} (\id_{(4)} \cotimes \mu)
$$
and
$$
q = \mu \cotimes \id \cotimes \mu \cotimes \id - \id \cotimes \mu \cotimes \id \cotimes \mu. 
$$
The source and target of $ \xi $ are the quotients of $ K^{\cotimes 4} $ by the closure of the image of $ p $ and $ q $, respectively. 
Using conditions a) and b) it is straightforward to verify the relation 
$$
q \kappa = \gamma_r^{13} p
$$
where $ \kappa $ is the bornological automorphism of $ K^{\cotimes 6} $ defined by 
$$
\kappa = (\id_{(2)} \cotimes \tau \cotimes \id_{(2)}) \gamma_r^{13} \gamma_r^{23} (\id_{(2)} \cotimes \tau \cotimes \id_{(2)}). 
$$
This relation shows that $ \xi $ is actually a bornological isomorphism. Using the 
map $ \rho_l $ one proves in a similar way that 
$ (K \cotimes K) \cotimes_K K \rightarrow K \cotimes K $ is an isomorphism. We conclude that $ \Delta $ is essential. 
Having established this, condition f) immediately yields that $ \Delta $ is coassociative. Hence $ \Delta $ is 
a comultiplication. \\
If there exists maps $ \gamma_l $ and $ \rho_r $ with the properties stated in conditions g) and h) then 
these maps describe the remaining Galois maps associated to $ \Delta $. It follows in particular 
that all Galois maps yield isomorphisms from $ K \cotimes K $ into itself in this case. \qed \\
We have already shown above that the maps $ \hat{\gamma}_r $ and $ \hat{\rho}_l $ satisfy the assumptions of 
proposition \ref{galoistocomult}. Let us write $ \hat{\Delta} $ for the comultiplication on $ \hat{H} $ defined in this way. 
By construction, the Galois maps $ \hat{\gamma}_r $ and $ \hat{\rho}_l $ associated to $ \hat{\Delta} $ 
are isomorphisms. \\
To treat the remaining Galois maps for $ \hat{\Delta} $ let us abstractly define 
$$
\hat{\gamma}_l(\mathcal{F}_r(f) \otimes \mathcal{F}_r(g)) = (\mathcal{F}_r \cotimes \mathcal{F}_r)\tau \rho_l^{-1}(f \otimes g)
$$
and 
$$
\hat{\rho}_r(\mathcal{G}_l(g) \otimes \mathcal{G}_l(f)) = (\mathcal{G}_l \cotimes \mathcal{G}_l) \tau \gamma_r^{-1}(g \otimes f).
$$
Applying the above discussion to $ H^{op} $ we see that $ \hat{\gamma}_l $ and $ \hat{\rho}_r $ satisfy conditions 
g) and i) in proposition \ref{galoistocomult}. 
Using equations (\ref{phiinvarproj}) and (\ref{psiinvarproj}) it is straightforward to obtain the formulas 
\begin{equation}\label{hatgammalex}
\hat{\gamma}_l(\omega \otimes \mathcal{F}_r(g))(x \otimes y) = (\omega \cotimes \phi)(\mu \cotimes \id)
(\id \cotimes \tau) (\rho_r \cotimes \id)(g \otimes x \otimes y)
\end{equation}
and 
\begin{equation}\label{hatrhorex}
\hat{\rho}_r(\mathcal{G}_l(f) \otimes \omega)(x \otimes y) = (\psi \cotimes \omega)(\id \cotimes \mu)(\tau \cotimes \id)(\id \cotimes \gamma_l)(x \otimes y \otimes f)
\end{equation}
for the maps $ \hat{\gamma}_l $ and $ \hat{\rho}_r $. 
According to the definition of $ \hat{\mu} $ and equations (\ref{hatgammalex}), (\ref{phiinvarproj}) and (\ref{gammahatrev}) we compute 
\begin{align*}
(\id \cotimes &\hat{\mu})(\hat{\gamma}_l \cotimes \id)(\omega \otimes \mathcal{F}_r(f) \otimes \mathcal{F}_l(g))(x \otimes y) = 
(\hat{\gamma}_l(\omega \otimes \mathcal{F}_r(f)) \cotimes \phi)(x \otimes \gamma_r(y \otimes g)) \\
&= (\omega \cotimes \phi \cotimes \phi)(\mu \cotimes \id_{(2)}) (\id \otimes \tau \cotimes \id)(\rho_r \cotimes \gamma_r)(f \otimes x \otimes y \otimes g) \\
&= (\omega \cotimes \phi \cotimes \phi)(\mu \cotimes \id_{(2)}) (\id \otimes \tau \cotimes \id)(S \cotimes \id \cotimes S^{-1} \cotimes \id)
(\gamma_r \cotimes \rho_r)(f \otimes x \otimes y \otimes g) \\
&= (\hat{\gamma}_r(\omega \otimes \mathcal{F}_l(g)) \cotimes \phi)(S \cotimes \id_{(2)})(\id \cotimes \tau)(\gamma_r \cotimes \id)(f \otimes x \otimes y)\\
&= (\hat{\mu} \cotimes \id)(\id \cotimes \tau)(\hat{\gamma}_r \cotimes \id)(\id \cotimes \tau)
(\omega \otimes \mathcal{F}_r(f) \otimes \mathcal{F}_l(g))(x \otimes y) 
\end{align*}
which yields condition h). Using equations (\ref{hatrhorex}), (\ref{psiinvarproj}) and (\ref{rhohatlev}) one obtains condition j) in a similar way. 
Hence it follows from proposition \ref{galoistocomult} that all Galois maps associated to $ \hat{\Delta} $ are isomorphisms. \\
It remains to exhibit the Haar functionals for the comultiplication $ \hat{\Delta} $. 
\begin{prop}\label{dualHaar}
Let $ H $ be a bornological quantum group. Then the linear form $ \hat{\psi} $ on $ \hat{H} $ defined by 
$$
\hat{\psi}(\mathcal{F}_l(f)) = \epsilon(f)
$$
is a faithful right invariant functional on $ \hat{H} $. Similarly, the linear form $ \hat{\phi} $ given by 
$$
\hat{\phi}(\mathcal{G}_r(f)) = \epsilon(f)
$$
is a faithful left invariant functional on $ \hat{H} $.
\end{prop}
\proof We have already checked above that $ \hat{\psi} $ is faithful. Moreover we compute
\begin{align*}
(\hat{\psi} \cotimes &\id)\hat{\gamma}_r(\mathcal{F}_l(f) \otimes \mathcal{F}_l(g)) = 
(\id \cotimes \hat{\psi}) (\mathcal{F}_l \cotimes \mathcal{F}_l)\gamma_l^{-1}(f \otimes g) \\
&= \mathcal{F}_l(\id \cotimes \epsilon)\gamma_l^{-1}(f \otimes g) = \epsilon(f) \mathcal{F}_l(g) = \hat{\psi}(\mathcal{F}_l(f)) \mathcal{F}_l(g)
\end{align*}
and we deduce that $ \hat{\psi} $ is right invariant. 
The assertions concerning $ \hat{\phi} $ are obtained in a similar way. \qed \\
We have now completed to proof of the following theorem.  
\begin{theorem} 
Let $ H $ be a bornological quantum group. Then $ \hat{H} $ with the structure maps described above is again a 
bornological quantum group. 
\end{theorem}
The bornological quantum group $ \hat{H} $ will be called the dual quantum group of $ H $. It is instructive to describe 
explicitly the counit and the antipode of $ \hat{H} $. Consider the map
$ \hat{\epsilon}: \hat{H} \rightarrow \mathbb{C} $ given by 
$$
\hat{\epsilon}(\omega) = \omega(1)
$$
where $ \hat{H} $ is viewed as a subspace of $ M(H)' $ according to the nondegenerate pairing $ \hat{H} \times M(H) \rightarrow \mathbb{C} $. 
The explicit formulas  
$$ 
\hat{\epsilon}(\mathcal{F}_l(f)) = \phi(f),\quad \hat{\epsilon}(\mathcal{F}_r(f)) = \phi(f), \quad \hat{\epsilon}(\mathcal{G}_l(f)) = \psi(f), \quad \hat{\epsilon}(\mathcal{G}_r(f)) = \psi(f)
$$ 
show that the map $ \hat{\epsilon} $ is bounded and nonzero. 
It is straightforward to check that $ \hat{\epsilon} $ is an algebra homomorphism and 
we calculate
\begin{align*}
(\hat{\epsilon} \cotimes &\id)\hat{\gamma}_r(\mathcal{F}_l(f) \otimes \mathcal{F}_l(h)) = (\mathcal{F}_l \cotimes \phi)\gamma_l^{-1}(f \otimes h) 
= \hat{\mu}(\mathcal{F}_l(f) \otimes \mathcal{F}_l(h))
\end{align*}
as well as 
\begin{align*}
(\id \cotimes &\hat{\epsilon})\hat{\rho}_l(\mathcal{G}_r(h) \otimes \mathcal{G}_r(f)) = (\psi \cotimes \mathcal{G}_r)\rho_r^{-1}(h \otimes f) = 
\hat{\mu}(\mathcal{G}_r(h) \otimes \mathcal{G}_r(f))
\end{align*}
which shows $ (\hat{\epsilon} \cotimes \id) \gamma_r = \hat{\mu} $ and $ (\id \cotimes \hat{\epsilon}) \rho_l = \hat{\mu} $. 
One can then proceed as in the proof of theorem \ref{bqchar} to show that 
$ \hat{\epsilon} $ is nondegenerate. By the uniqueness assertion of theorem \ref{bqchar} we see that the map $ \hat{\epsilon} $ 
is indeed the counit for $ \hat{H} $. \\
Similarly, we define $ \hat{S}: \hat{H} \rightarrow \hat{H} $ by 
$$
\hat{S}(\omega)(f) = \omega(S(f))
$$
and using $ \psi = S(\phi) $ we obtain the formulas 
$$ 
\hat{S}(\mathcal{F}_l(f)) = \mathcal{G}_r(S^{-1}(f)), \quad \hat{S}(\mathcal{F}_r(f)) = \mathcal{G}_l(S^{-1}(f)). 
$$ 
It follows that $ \hat{S} $ is a bounded linear automorphism of $ \hat{H} $. Using equation (\ref{phiinvarproj}) we compute 
\begin{align*}
\hat{\mu}(\hat{S} \cotimes \id)&(\mathcal{F}_l(f) \otimes \mathcal{F}_l(g))(x) = 
\hat{\mu}(\mathcal{G}_r(S^{-1}(f)) \otimes \mathcal{F}_l(g))(x) \\
&= (\psi \cotimes \phi)(\mu \cotimes \id)(S^{-1} \cotimes \id_{(2)})(\id \cotimes \gamma_r)(f \otimes x \otimes g) \\
&= (\phi \cotimes \phi)(\mu \cotimes \id)(\tau \cotimes \id)(\id \cotimes S \cotimes \id)(\id \cotimes \gamma_r)(f \otimes x \otimes g) \\
&= (\phi \cotimes \phi)(\mu \cotimes \id)(\tau \cotimes \id)(\id \cotimes \rho_r)(f \otimes x \otimes g) \\
&= (\hat{\epsilon} \cotimes \id) (\mathcal{F}_l \cotimes \mathcal{F}_l)\gamma_l \tau (f \otimes g)(x) \\
&= (\hat{\epsilon} \cotimes \id) \hat{\gamma}_r^{-1}(\mathcal{F}_l(f) \otimes \mathcal{F}_l(g))(x)
\end{align*}
which shows $ \hat{\mu}(\hat{S}\cotimes \id) = (\hat{\epsilon} \cotimes \id) \hat{\gamma}_r^{-1} $. In a 
similar way one obtains the relation $ \hat{\mu}(\id \cotimes \hat{S}) = (\id \cotimes \hat{\epsilon})\hat{\rho}_l^{-1} $. Inspecting 
the constructions in the proof of theorem \ref{bqchar} we see that $ \hat{S} $ is the antipode of $ \hat{H} $. \\
Let us now prove the Pontrjagin duality theorem. 
\begin{theorem} \label{pontrjagin}
Let $ H $ be a bornological quantum group. Then the double dual quantum group of $ H $ is canonically 
isomorphic to $ H $.
\end{theorem}
\proof We define a linear map $ P: H \rightarrow (\hat{H})' $ by 
$$
P(f)(\omega) = \omega(f) 
$$
for all $ f \in H $ and $ \omega \in \hat{H} $. According to proposition \ref{dualHaar} we compute
\begin{align*}
(\hat{\mathcal{G}}_l \mathcal{F}_l(f))(\mathcal{F}_l(h)) = \hat{\psi}(\mathcal{F}_l(h) &\mathcal{F}_l(f)) 
= \hat{\psi}(\mathcal{F}_l \cotimes \phi)\gamma_l^{-1}(f \otimes h) \\
&= (\epsilon \cotimes \phi)\gamma_l^{-1}(f \otimes h) = \phi(S^{-1}(f)h) = \mathcal{F}_l(h)(S^{-1}(f))
\end{align*}
for all $ f, g \in H $ where $ \hat{\mathcal{G}}_l $ is the map $ \mathcal{G}_l $ for $ \hat{H} $. This implies 
$$
P(f) = \hat{\mathcal{G}}_l \mathcal{F}_l(S(f))
$$
and shows that $ P $ defines a bornological isomorphism from $ H $ to $ \hat{\hat{H}} $. In a similar way one has 
$$
P(f) = \hat{\mathcal{F}}_r \mathcal{G}_r(S(f)).  
$$
Let us also remark that using $ \psi = S^{-1}(\phi) $ one calculates
\begin{align*}
(\hat{\mathcal{F}}_l S \mathcal{G}_l(f))&(\mathcal{G}_r(h)) = \hat{\phi}(\mathcal{G}_r(h) \mathcal{F}_r(S^{-1}(f))) 
= \hat{\phi}(\phi \cotimes \mathcal{G}_r)\rho_r^{-1} (h \otimes S^{-1}(f)) \\
&= (\phi \cotimes \epsilon)\rho_r^{-1}(h \otimes S^{-1}(f)) = \phi(S^{-1}(f) S^{-1}(h)) = \mathcal{G}_r(h)(f)
\end{align*}
which shows $ P = \hat{\mathcal{F}}_l S \mathcal{G}_l $. \\
Next consider the transpose $ \mu^*: H' \rightarrow (H \cotimes H)' $ of the multiplication map given by 
$$
\mu^*(\omega)(f \otimes g) = \omega \mu(f \otimes g) = \omega(fg) 
$$
for all $ f, g \in H $. In particular, we obtain a bounded linear map $ \mu^*: \hat{H} \rightarrow (H \cotimes H)' $ 
by restriction. Using the isomorphism $ P $ we can view $ \hat{\Delta} $ as a map 
from $ \hat{H} \rightarrow (H \cotimes H)' $ as well. Equivalently, we have bounded linear maps 
from $ \hat{H} $ into $ \Hom(H, H') $ given by 
$$ 
\mu^*(\omega)(g)(f) = \omega \mu(f \otimes g) 
$$ 
and likewise for $ \hat{\Delta} $. Using equation (\ref{gammartogamml}) we calculate 
\begin{align*}
\hat{\Delta}&(\mathcal{F}_l(h))(g) = (\id \cotimes \hat{\psi})\hat{\gamma}_r(\mathcal{F}_l(h) \otimes \mathcal{F}_l S(g)) \\
&= (\id \cotimes \hat{\psi})(\mathcal{F}_l \cotimes \mathcal{F}_l) \tau \gamma_l^{-1} (\id \cotimes S)(h \otimes g) \\
&= (\epsilon \cotimes \mathcal{F}_l) (S \cotimes \id) \gamma_r \tau (h \otimes g) \\
&= \mathcal{F}_l(gh)
\end{align*}
and thus obtain using the definition of $ S^{-1} $ 
\begin{align*}
\hat{\Delta}&(\mathcal{F}_l(h))(g)(f) = \hat{\psi}(\mathcal{F}_l(gh)\mathcal{F}_l S(f)) = \hat{\psi}(\mathcal{F}_l \cotimes \phi) \gamma_l^{-1}(gh \otimes S(f)) \\
&= (\epsilon \cotimes \phi) \gamma_l^{-1}(gh \otimes S(f)) = \phi(fgh) = \mu^*(\mathcal{F}_l(h))(g)(f) 
\end{align*}
which shows that $ \hat{\Delta} $ can be identified with the transpose $ \mu^* $ of the multiplication. Similarly, we 
have seen in the constructions above that $ \hat{\mu} $ can be identified with the transpose $ \Delta^* $ of the 
comultiplication. \\
With this in mind it is straightforward to check that $ P $ is an algebra homomorphism and a coalgebra homomorphism.
Hence $ P $ is an isomorphism of bornological quantum groups. \qed 

\section{Duality for modules and comodules}\label{secdualmod}

In this section we study the duality between essential modules and comodules over a bornological quantum group and 
its dual. \\
Let $ H $ be a bornological quantum group and let $ \eta: V \rightarrow \Hom_H(H, V \cotimes H) $ be an essential $ H $-comodule. We define a 
bounded linear map $ D(\eta): \hat{H} \cotimes V \rightarrow V $ by
$$
D(\eta)(\mathcal{F}_l(f) \otimes v) = (\id \cotimes \phi)\eta(v \otimes f). 
$$
For later use we need another description of this map. Since $ H $ is an essential algebra we may view $ \eta $ as a bounded linear map from 
$ V $ into $ \Hom_H(H, V \cotimes H) \cong \Hom_H(H \cotimes_H H, V \cotimes H) $. Under the latter isomorphism $ \eta(v) $ 
corresponds to the map $ (\id \cotimes \mu)(\eta(v) \cotimes \id) $. Moreover, using notation and results from section 
\ref{secmodhaar} we have 
$$
\phi(hg \nu(f)) = \psi(hgf ) = \phi(h \nu(gf)) 
$$
for all $ f,g, h \in H $ which implies that $ \nu $ is left $ H $-linear. Together with the relation 
$$
(\id \cotimes \phi)(\id \cotimes \mu)(\eta(v)(g) \otimes \nu(f)) = (\id \cotimes \psi)(\id \cotimes \mu)(\eta(v)(g) \otimes f)
$$
we thus obtain 
$$
D(\eta)(\mathcal{G}_l(f) \otimes v) = (\id \cotimes \psi)\eta(v \otimes f). 
$$
Next we compute 
\begin{align*}
D(\eta) (\id \cotimes &D(\eta))(\mathcal{F}_l(f) \otimes \mathcal{F}_l(g) \otimes v) = (\id \cotimes \phi)\eta(D(\eta)(\mathcal{F}_l(g) \otimes v) \otimes f) \\
&= (\id \cotimes \phi)\eta (\id \cotimes \id \cotimes \phi) \eta_{13} (v \otimes f \otimes g) \\
&= (\id \cotimes \phi \cotimes \phi) \eta_{12} \eta_{13} (v \otimes f \otimes g) \\
&= (\id \cotimes \phi \cotimes \phi) (\id \otimes \gamma_r) \eta_{12} (\id \otimes \gamma_r^{-1})(v \otimes f \otimes g). 
\end{align*}
According to equation (\ref{phiinvarproj}) we have 
$$
(\phi \cotimes \phi)\gamma_r = (\phi \cotimes \phi)(S^{-1} \cotimes \id)\rho_r 
$$
and using equation (\ref{coop12}) and equation (\ref{etalinearr})
we obtain
\begin{align*}
D(\eta) (\id \cotimes &D(\eta))(\mathcal{F}_l(f) \otimes \mathcal{F}_l(g) \otimes v) \\
&= (\id \cotimes \phi \cotimes \phi) (\id \cotimes S^{-1} \cotimes \id)(\id \cotimes \rho_r) \eta_{12}(\id \cotimes \gamma_r^{-1})(v \otimes f \otimes g) \\
&= (\id \cotimes \phi \cotimes \phi) \eta_{13}(\id \cotimes S^{-1} \cotimes \id)(\id \cotimes \rho_l)(\id \cotimes S \cotimes \id) (v \otimes f \otimes g) \\
&= (\id \cotimes \phi \cotimes \phi) \eta_{12} (\id \cotimes \gamma_l^{-1}) (v \otimes f \otimes g) \\
&= D(\eta)(\hat{\mu} \cotimes \id) (\mathcal{F}_l(f) \otimes \mathcal{F}_l(g) \otimes v) 
\end{align*}
which shows that $ V $ becomes a left $ \hat{H} $-module in this way. \\
We want to show that $ V $ is actually an essential $ \hat{H} $-module. 
In order to do this it is convenient to work with the map $ (\id \cotimes \phi) \eta $ 
instead of $ D(\eta) $. There is an evident bounded linear splitting $ \sigma: V \rightarrow V \cotimes H $ of this map given by 
$ \sigma(v) = \eta^{-1}(v \otimes h) $ where $ h $ is chosen such that $ \phi(h) = 1 $.  
If we identify $ \hat{H} \cotimes_{\hat{H}} V $ accordingly with a quotient $ Q $ of $ V \cotimes H $ 
we have the relation
$$
(\id \cotimes \id \cotimes \phi)\eta_{13} = (\id \cotimes \id \cotimes \phi)(\id \cotimes \gamma_l^{-1})
$$
in this quotient. Now we see as in the proof of proposition \ref{Hproj} and according to formula (\ref{gammartogamml})
\begin{align*}
(\id \cotimes \id \cotimes &\phi)\eta_{12}\eta_{13} = (\id \cotimes \phi \cotimes \id) \eta_{12}(\id \cotimes S \cotimes \id)(\id \cotimes \gamma_r)
(\id \cotimes S^{-1}\cotimes \id)(\id \cotimes \tau) \\
&= (\id \cotimes \phi \cotimes \id)\eta_{12}(\id \otimes \gamma_l^{-1})
\end{align*}
which implies
\begin{align*}
\eta^{-1}(\id \cotimes \phi \cotimes \id) \eta_{12} &= \eta^{-1}(\id \cotimes \id \cotimes \phi) \eta_{12} \eta_{13}(\id \cotimes \gamma_l) \\
&= (\id \cotimes \id \cotimes \phi) \eta_{13} (\id \cotimes \gamma_l) = \id \cotimes \id \cotimes \phi
\end{align*}
in $ Q $. It follows from this relation that $ \sigma(\id \otimes \phi)\eta $ is the identity map on $ Q $. 
Translating this back to $ \hat{H} \cotimes_{\hat{H}} V $ we deduce that $ V $ is an essential module. \\
An $ H $-colinear map $ f: V \rightarrow W $ is easily seen to be $ \hat{H} $-linear for the module structures 
defined in this way. Hence we have proved the following statement. 
\begin{prop}\label{comodmod}
Let $ H $ be a bornological quantum group and let $ \hat{H} $ be the dual quantum group. The previous construction 
defines a functor $ D $ from $ \SComodR H $ to $ \hat{H} \LSMod $. 
\end{prop}
Conversely, let $ \lambda: H \cotimes V \rightarrow V $ be an essential left $ H $-module. 
By slight abuse of notation we write $ \lambda^{-1} $ for the inverse of the isomorphism $ H \cotimes_H V \cong V $ induced 
by $ \lambda $. We define a bounded linear map $ D(\lambda): V \cotimes \hat{H} \rightarrow V \cotimes \hat{H} $ by 
$$
D(\lambda)(v \otimes \mathcal{F}_l(f)) = (\id \cotimes \mathcal{F}_l)\tau(\id \cotimes \lambda)(\gamma_l^{-1} \tau \cotimes \id) 
(\id \cotimes \lambda^{-1})(f \otimes v) 
$$
which is seen to be well-defined since $ \gamma_l^{-1} \tau $ is right $ H $-linear for the action by multiplication on the 
second tensor factor. It is evident that $ D(\lambda) $ is an isomorphism. Since $ \lambda $ is left 
$ H $-linear we calculate with $ \rho = \gamma_l^{-1} \tau $ 
\begin{align*}
(\id \cotimes &\lambda)(\gamma_l^{-1}\tau \cotimes \id)(\id \cotimes \lambda^{-1})(\id \cotimes \phi \cotimes \id)(\gamma_l^{-1} \cotimes \id)
(f \otimes g \otimes v) \\
&= (\phi \cotimes \id \cotimes \id)(\id_{(2)} \cotimes \lambda)(\id \cotimes \gamma_l^{-1}\tau \cotimes \id)(\id_{(2)} \cotimes \lambda^{-1})
(\tau\gamma_l^{-1} \cotimes \id)(f \otimes g \otimes v) \\
&= (\id \cotimes \phi \cotimes \id)(\id_{(2)} \cotimes \lambda)(\tau \cotimes \id_{(2)})\rho^{23} \\
&\qquad (\tau \cotimes \id_{(2)})(\id_{(2)} \cotimes \lambda^{-1}) \rho^{12}(\tau \cotimes \id)(f \otimes g \otimes v) \allowdisplaybreaks \\
&= (\id \cotimes \phi \cotimes \id)(\id_{(2)} \cotimes \lambda)\rho^{13} \rho^{12} (\tau \cotimes \id_{(2)})(\id_{(2)} \cotimes \lambda^{-1})(f \otimes g \otimes v)\\
&= (\id \cotimes \phi \cotimes \id)(\id_{(2)} \cotimes \lambda)\rho^{12} \rho^{23} (\tau \cotimes \id_{(2)})(\id_{(2)} \cotimes \lambda^{-1})(f \otimes g \otimes v)\\
&= (\id \cotimes \phi \cotimes \id)(\gamma_l^{-1}\cotimes \id)(\tau \cotimes \id)(\id_{(2)} \cotimes \lambda)(\id \cotimes \gamma_l^{-1} \tau \cotimes \id) \\
&\qquad (\id_{(2)} \cotimes \lambda^{-1})(\tau \cotimes \id)(f \otimes g \otimes v)
\end{align*}
using $(\phi \cotimes \id) = (\phi \cotimes \id)\rho $ as well as the pentagon relation (\ref{pentrho}) for the map $ \rho $. 
This shows 
$$
D(\lambda)(\id \cotimes \hat{\mu})(v \otimes \mathcal{F}_l(f) \otimes \mathcal{F}_l(g)) = 
(\id \cotimes \hat{\mu})(D(\lambda) \cotimes \id)(v \otimes \mathcal{F}_l(f) \otimes \mathcal{F}_l(g))
$$
which means that $ D(\lambda) $ is right $ \hat{H} $-linear. Again by the pentagon relation for $ \rho $ and the fact that 
$ \lambda $ is $ H $-linear we have 
\begin{align*}
(&\id_{(2)} \cotimes \lambda)(\tau \cotimes \id_{(2)})(\id \cotimes \gamma_l^{-1} \tau \cotimes \id)(\tau \cotimes \id_{(2)}) \\
&\qquad 
(\id \cotimes \gamma_l^{-1} \tau \cotimes \id)(\id_{(2)} \cotimes \lambda^{-1})(\tau \gamma_l^{-1} \cotimes \id)(f \otimes g \otimes v) \\ 
&= (\id_{(2)} \cotimes \lambda) (\tau \cotimes \id_{(2)}) \rho^{23} \rho^{13} \rho^{12}(\tau \cotimes \id_{(2)})(\id_{(2)} \cotimes \lambda^{-1})
(f \otimes g \otimes v) \\
&= (\id_{(2)} \cotimes \lambda) (\tau \cotimes \id_{(2)}) \rho^{12} \rho^{23}(\tau \cotimes \id_{(2)})(\id_{(2)} \cotimes \lambda^{-1})
(f \otimes g \otimes v) \\
&= (\tau \gamma_l^{-1} \cotimes \id)(\tau \cotimes \id)(\id_{(2)} \cotimes \lambda)(\id \cotimes \gamma_l^{-1} \tau \cotimes \id)
(\tau \cotimes \id_{(2)})(\id_{(2)} \cotimes \lambda^{-1})(f \otimes g \otimes v)
\end{align*}
which shows 
$$
D(\lambda)_{12} D(\lambda)_{13}(\id \cotimes \hat{\gamma}_r)(v \otimes \mathcal{F}(f) \otimes \mathcal{F}(g)) = 
(\id \cotimes \hat{\gamma}_r) D(\lambda)_{12}(v \otimes \mathcal{F}(f) \otimes \mathcal{F}(g)). 
$$
Hence $ D(\lambda) $ is a right coaction of $ \hat{H} $ on $ V $. It is easy to check that an $ H $-equivariant 
map $ f: V \rightarrow W $ between $ H $-modules defines an $ \hat{H} $-colinear map between the associated comodules. 
\begin{prop}\label{modcomod}
Let $ H $ be a bornological quantum group and let $ \hat{H} $ be the dual quantum group. There is a natural functor 
from $ H\LSMod $ to $ \SComodR \hat{H} $ which will again be denoted by $ D $. 
\end{prop}
We obtain the following duality theorem for modules and comodules. 
\begin{theorem} \label{modcomod}
Let $ H $ be a bornological quantum group. 
Every essential left $ H $-module is an essential right $ \hat{H} $-comodule in a natural way and vice versa.  
This yields inverse isomorphisms between the category of essential $ H $-modules and the category 
of essential $ \hat{H} $-comodules. These isomorphisms are compatible with tensor products. 
\end{theorem}
\proof Let us check that the functors defined above are inverse to each other if we take into account the 
Pontrjagin duality theorem \ref{pontrjagin}. According to equation (\ref{coop2}) we have 
$$
\gamma_l^{-1} = \tau(S^{-1} \cotimes \id)\rho_l(S \cotimes \id)
$$
and hence
$$
(\epsilon \cotimes \id)\gamma_l^{-1}\tau(S \cotimes \id) = (\id \cotimes \epsilon)(S^{-1} \cotimes \id)\rho_l(S \cotimes S) \tau = 
S^{-1} \mu(S \cotimes S) \tau = \mu.
$$
Using the definition of the right Haar functional $ \hat{\psi} $ on $ \hat{H} $ we thus compute for an essential 
$ H $-module $ \lambda: H \cotimes V \rightarrow V $ 
\begin{align*}
(\id &\cotimes \hat{\psi})D(\lambda)(v \otimes \mathcal{F}_l(S(f))) \\
&= (\hat{\psi} \cotimes \id)(\mathcal{F}_l \cotimes \id)(\id \cotimes \lambda)
(\gamma_l^{-1}\tau \cotimes \id)(\id \cotimes \lambda^{-1})(S \cotimes \id)(f \otimes v) \\
&= (\epsilon \cotimes \id)(\id \cotimes \lambda)
(\gamma_l^{-1}\tau \cotimes \id)(S \cotimes \id_{(2)})(\id \cotimes \lambda^{-1})(f \otimes v) \\
&= \lambda(\epsilon \cotimes \id_{(2)})(\gamma_l^{-1}\tau \cotimes \id)(S \cotimes \id_{(2)})(\id \cotimes \lambda^{-1})(f \otimes v) \\
&= \lambda(\mu \cotimes \id)(\id \cotimes \lambda^{-1})(f \otimes v) \\
&= \lambda(\id \cotimes \lambda)(\id \cotimes \lambda^{-1})(f \otimes v) = \lambda(f \otimes v).
\end{align*}
Consequently we have 
$$
DD(\lambda)(\hat{\mathcal{G}}_l \mathcal{F}_l(S(f)) \otimes v) = \lambda(f \otimes v)
$$
and according to Pontrjagin duality this shows that the module structure $ DD(\lambda) $ can be identified with $ \lambda $. \\
Conversely, let $ \eta: V \rightarrow \Hom_H(V, V \cotimes H) $ be an essential $ H $-comodule. 
Using equation (\ref{gammartogamml}) for $ \hat{H} $ we compute 
\begin{align*}
DD(\eta)(v \otimes &\hat{\mathcal{F}}_l S(\omega)) = (\id \cotimes \hat{\mathcal{F}}_l)\tau(\id \cotimes D(\eta)) (\hat{\gamma}_l^{-1} \tau \cotimes \id)
(\id\cotimes D(\eta)^{-1})(S(\omega) \otimes v) \\
&= (\id \cotimes \hat{\mathcal{F}}_l) \tau (\id \cotimes D(\eta))(S \cotimes \id_{(2)}) (\hat{\gamma}_r \cotimes \id)(\id \cotimes D(\eta)^{-1})(\omega \otimes v) \\
&= \tau (\hat{\mathcal{F}}_l S \cotimes \id)(\id \cotimes D(\eta))(\hat{\gamma}_r \cotimes \id)(\id \cotimes D(\eta)^{-1})(\omega \otimes v) 
\end{align*}
and thus 
obtain 
\begin{align*}
DD(\eta)&(D(\eta)(\mathcal{F}_l(g) \otimes v) \otimes \hat{\mathcal{F}}_l S \mathcal{G}_l(f)) \\ 
&= \tau (\hat{\mathcal{F}}_l S \cotimes \id)(\id \cotimes D(\eta)) (\hat{\gamma}_r \cotimes \id)(\mathcal{G}_l(f) \otimes \mathcal{F}_l(g) \otimes v) \\
&= (\id \cotimes \hat{\mathcal{F}}_l S)\tau (\id \cotimes D(\eta)) (\mathcal{G}_l \cotimes \mathcal{F}_l \cotimes \id)(\tau \gamma_l^{-1} \cotimes \id)
(f \otimes g \otimes v) \\
&= (\id \cotimes \hat{\mathcal{F}}_l S)(D(\eta) \cotimes \id)(\tau \cotimes \id)(\id \cotimes \mathcal{F}_l \cotimes \mathcal{G}_l)(\id \cotimes \gamma_l^{-1})
(v \otimes f \otimes g) \allowdisplaybreaks \\
&= (\id \cotimes \hat{\mathcal{F}}_l S \mathcal{G}_l) (D(\eta) \cotimes \id)
(\tau \cotimes \id)(\id \cotimes \mathcal{F}_l \cotimes \id) (\id \cotimes \gamma_l^{-1})(v \otimes f \otimes g) \allowdisplaybreaks \\
&= (\id \cotimes \hat{\mathcal{F}}_l S \mathcal{G}_l) (\id \cotimes \phi \cotimes \id) \eta^{12} (\id \cotimes \gamma_l^{-1})(v \otimes f \otimes g) \allowdisplaybreaks \\
&= (\id \cotimes \hat{\mathcal{F}}_l S \mathcal{G}_l) (\id \cotimes \id \cotimes \phi) \eta^{12} \eta^{13} (v \otimes f \otimes g) \\
&= (\id \cotimes \hat{\mathcal{F}}_l S \mathcal{G}_l) \eta(D(\eta)(\mathcal{F}_l(g) \otimes v) \otimes f)
\end{align*}
which implies
$$
DD(\eta)(v \otimes \hat{\mathcal{F}}_l S\mathcal{G}_l(f)) = (\id \cotimes \hat{\mathcal{F}}_l S \mathcal{G}_l) \eta(v \otimes f)
$$
since $ D(\eta) $ is an essential $ \hat{H} $-module. 
Again by Pontrjagin duality this shows that $ DD(\eta) $ is isomorphic to $ \eta $. \\
Consider $ H \cotimes_H(H \cotimes (H \cotimes_H H)) $ and $ (H \cotimes_H H) \cotimes (H \cotimes_H H) $ as $ H $-modules by multiplication 
on the first tensor factor and by the diagonal action on the first and third tensor factors, respectively. 
Then the isomorphism $ \xi $ used in the proof of proposition \ref{galoistocomult} is 
$ H $-linear. Using this observation it is straightforward to check that the functor $ D $ from $ H\LSMod $ to $ \SComodR H $ is 
compatible with tensor products. \qed \\
Of course there is an analogue of theorem \ref{modcomod} for right modules and left comodules. 
Let us use the above duality results to construct the dual of a morphism between bornological quantum groups. 
\begin{prop}\label{morphismdual}
Let $ \alpha: H \rightarrow M(K) $ be a morphism of bornological quantum groups. Then there exists 
a unique morphism $ \hat{\alpha}: \hat{K} \rightarrow M(\hat{H}) $ such that 
$$
\bra \alpha(f), \omega \ket = \bra f, \hat{\alpha}(\omega) \ket 
$$
for all $ f \in H $ and $ \omega \in \hat{K} $. 
\end{prop}
\proof Uniqueness of $ \hat{\alpha} $ follows immediately from the nondegeneracy of the pairing between $ H $ and $ M(\hat{H}) $. 
Consider the transposed right regular coaction $ \rho = \gamma_l^{-1} \tau $ on $ H $. The dual action of the pushforward coaction $ \alpha_*(\rho) $ 
yields a left $ \hat{K} $-module structure on $ H $. Using the linear isomorphism $ \mathcal{F}_l $ we may view this 
as a $ \hat{K} $-module structure on $ \hat{H} $. Associativity of the multiplication in $ \hat{H} $ and equation (\ref{muhatrho}) shows that  
we obtain in fact a bounded linear map 
$ \hat{\alpha}_l: \hat{K} \rightarrow M_l(\hat{H}) $. 
Similarly, the map $ \gamma = \rho_r^{-1} \tau $ defines a left coaction 
of $ H $ on itself, and the dual action of the corresponding pushforward coaction determines a right $ \hat{K} $-module structure 
on $ H $. This action yields a homomorphism $ \hat{\alpha}_r: \hat{K} \rightarrow M(\hat{H}) $. Using lemma \ref{hopflemma1} for $ H^{cop} $ we obtain 
$$
(\id \cotimes \rho)(\gamma \cotimes \id) = (\gamma \cotimes \id)(\id \cotimes \rho)
$$
and hence the resulting left and right $ \hat{K} $-module structures on $ \hat{H} $ commute. Consequently, the maps $ \hat{\alpha}_l $ and 
$ \hat{\alpha}_r $ yield a nondegenerate homomorphism $ \hat{\alpha}: \hat{K} \rightarrow M(\hat{H}) $. \\
Consider also the transpose $ \alpha^*: \hat{K} \rightarrow H' $ of 
$ \alpha $ given by 
$$
\alpha^*(\omega)(f) = \omega(\alpha(f)).  
$$
Then we have 
$$
\bra f, \alpha^*(\omega) \ket = \bra \alpha(f), \omega \ket.
$$
Moreover the calculation after equation (\ref{coop111}) for $ H^{cop} $ gives 
$$
(\mu \cotimes \id)(\id \cotimes \rho) = (\id \cotimes \mu)(\id \cotimes S^{-1} \cotimes \id)(\tau \rho_r \cotimes \id)
$$ 
and using equation (\ref{phiinvarproj}) and the definition of $ \hat{\alpha}_l $ as well as the definition of $ \hat{\mu}: H' \cotimes \hat{H} 
\rightarrow H' $ we calculate 
\begin{align*}
\hat{\mu}(\hat{\alpha}_l &\cotimes\id)(\mathcal{F}_l(k) \otimes \mathcal{F}_l(f))(h) \\
&= (\phi_H \cotimes \phi_K)(\mu \cotimes \id)(\id_{(2)} \cotimes \lambda_l)(\id \cotimes \rho \cotimes \id)(\id_{(2)} \cotimes \lambda_l^{-1})
(h \otimes f \otimes k) \\
&= (\phi_H \cotimes \phi_K)(\id \cotimes \lambda_l)(\id \cotimes \mu \cotimes \id)(\id \cotimes S^{-1} \cotimes \id_{(2)})(\tau\rho_r \cotimes \id_{(2)}) \\
&\qquad (\id_{(2)} \cotimes \lambda_l^{-1})(h \otimes f \otimes k) \allowdisplaybreaks \\
&= (\phi_H \cotimes \phi_K)(\id \cotimes \lambda_l)(\id \cotimes \mu \cotimes \id)(\tau \cotimes \id_{(2)})(\gamma_r \cotimes \id_{(2)}) \\
&\qquad (\id_{(2)} \cotimes \lambda_l^{-1})(h \otimes f \otimes k) \\
&= (\phi_H \cotimes \phi_K)(\id \cotimes \mu)(\id \cotimes \alpha \cotimes \id)(\tau \cotimes \id)(\gamma_r \cotimes \id)(h \otimes f \otimes k) \\
&= \hat{\mu}(\alpha^* \cotimes\id)(\mathcal{F}_l(k) \cotimes \mathcal{F}_l(f))(h) 
\end{align*}
where $ \lambda_l $ denotes the isomorphism $ H \cotimes_H K \cong K $ induced by $ \alpha $. This shows
$$
\hat{\mu}(\hat{\alpha}_l(\mathcal{F}_l(k)) \otimes \mathcal{F}_l(f)) = \hat{\mu}(\alpha^*(\mathcal{F}_l(k)) \otimes \mathcal{F}_l(f))
$$ 
for all $ k \in K $ and $ f \in H $. Similarly we have
$$
\hat{\mu}(\mathcal{G}_r(f)) \otimes \hat{\alpha}_r(\mathcal{G}_r(k))) = \hat{\mu}(\mathcal{G}_r(f)) \otimes \alpha^*(\mathcal{G}_r(k)))
$$ 
and we obtain 
$$ 
\hat{\alpha}(\omega) = \alpha^*(\omega) 
$$ 
for all $ \omega \in \hat{K} $. \\
We shall only sketch how to show that $ \hat{\alpha} $ is a coalgebra homomorphism. Using equation (\ref{phiinvarproj}) one obtains 
\begin{align*}
(\phi \cotimes \phi)& \mu_{(2)}(\rho_l \cotimes \id_{(2)}) = (\phi \cotimes \phi)(\mu \cotimes \id)(\id \cotimes \mu \cotimes \id)\gamma_r^{24} \\
&= (\phi \cotimes \phi)(\mu \cotimes \id)(\id \cotimes \mu \cotimes \id)(\id \cotimes S^{-1} \cotimes \id_{(2)})\rho_r^{24} \\
&= (\phi \cotimes \phi) \mu_{(2)}(\id_{(2)} \cotimes \tau \gamma_l^{-1})
\end{align*}
which shows 
$$
\bra \rho_l(f \otimes g), \mathcal{F}_l(h) \otimes \mathcal{F}_l(k) \ket = \bra f \otimes g, \hat{\gamma}_r(\mathcal{F}_l(h) \otimes \mathcal{F}_l(k)) \ket 
$$
for all $ f,g,h,k \in H $. This relation extends to the case where $ f $ and $ g $ are multipliers of $ H $ and we have similar 
statements involving other Galois maps. Based on this we calculate
\begin{align*}
\bra f \otimes g,& (\hat{\alpha} \cotimes \hat{\alpha})\hat{\gamma}_r (\mathcal{F}_l(k) \otimes \mathcal{F}_l(l)) \ket = 
\bra (\alpha \cotimes \alpha)(f \otimes g), \hat{\gamma}_r(\mathcal{F}_l(k) \otimes \mathcal{F}_l(l)) \ket \\ 
&= \bra \rho_l(\alpha \cotimes \alpha)(f \otimes g), \mathcal{F}_l(k) \otimes \mathcal{F}_l(l) \ket \\ 
&= \bra (\alpha \cotimes \alpha) \rho_l(f \otimes g), \mathcal{F}_l(k) \otimes \mathcal{F}_l(l) \ket \\ 
&= \bra \rho_l(f \otimes g), (\hat{\alpha} \cotimes \hat{\alpha})(\mathcal{F}_l(k) \otimes \mathcal{F}_l(l)) \ket  \\ 
&= \bra f \otimes g, \hat{\gamma}_r (\hat{\alpha} \cotimes \hat{\alpha})(\mathcal{F}_l(k) \otimes \mathcal{F}_l(l)) \ket
\end{align*}
and deduce 
$$
(\hat{\alpha} \cotimes \hat{\alpha})\hat{\gamma}_r = \hat{\gamma}_r (\hat{\alpha} \cotimes \hat{\alpha}) 
$$ 
which easily implies that $ \hat{\alpha} $ is compatible with the comultiplication. \qed 

\section{Bornological quantum groups associated to Lie groups}\label{seclie}

In this section we describe a dual pair of bornological quantum groups associated naturally 
to every Lie group. These bornological quantum groups are generalizations of 
the Hopf algebra of functions $ C(G) $ and the group algebra $ \mathbb{C}G $ of 
a finite group $ G $. As a matter of fact, one can extend the constructions described below to 
arbitrary locally compact groups. We will comment on this at the end of this section. \\
If $ M $ is a smooth manifold we let $ \D(M) $ be the space of smooth functions on $ M $ with compact 
support. The space $ \D(M) $ is equipped with the bornology associated to its natural LF-topology. 
We need the following assertion which is straightforward to prove. 
\begin{lemma} \label{multog}
Let $ M $ be a smooth manifold. The multiplier algebra of the 
algebra $ \D(M) $ of smooth functions with compact support with pointwise multiplication 
is the algebra $ \E(M) $ of all smooth functions. 
\end{lemma}
Now let $ G $ be a Lie group. We choose a left Haar measure 
$ dt $ and denote the modular function of $ G $ by $ \delta $. Then we have 
$$
\int_G f(ts) dt = \delta(s) \int_G f(t) dt, \qquad \int_G f(t^{-1}) dt = \int_G \delta(t) f(t) dt 
$$
for all $ f \in \D(G) $. \\
Let us write $ C^\infty_c(G) $ for the bornological algebra of smooth functions on 
$ G $ with pointwise multiplication. Using lemma \ref{multog} one defines the comultiplication $ \Delta: C^\infty_c(G) 
\rightarrow M(C^\infty_c(G \times G)) $ by 
$$ 
\Delta(f)(r,s) = f(rs). 
$$ 
This homomorphism is easily seen to be nondegenerate and coassociative. 
\begin{prop} Let $ G $ be a Lie group. Then the algebra $ C^\infty_c(G) $ of smooth functions 
with compact support on $ G $ is a bornological Hopf algebra.
\end{prop} 
\proof It is straightforward to check that all Galois maps associated 
to $ \Delta $ are isomorphisms. A left invariant integral $ \phi $ for $ C^\infty_c(G) $ is given by 
integration, that is, 
$$
\phi(f) = \int_G f(t) dt 
$$
for all $ f \in C^\infty_c(G) $. \qed \\
Let us also consider the counit and the antipode for $ C^\infty_c(G) $. The counit $ \epsilon: C^\infty_c(G) \rightarrow \mathbb{C} $ 
is given by 
$$ 
\epsilon(f) = f(e) 
$$ 
where $ e $ is the unit element of $ G $. The antipode $ S: C^\infty_c(G) \rightarrow C^\infty_c(G) $ is defined by 
$$ 
S(f)(t) = f(t^{-1}) 
$$
for all $ f \in C^\infty_c(G) $. Evidently the relation $ S^2 = \id $ holds. The modular element in $ M(C^\infty_c(G)) $ 
is given by the modular function $ \delta $.  \\
Let us explicitly describe the dual of $ C^\infty_c(G) $. We write $ \D(G) $ for this bornological quantum group 
and refer to it as the smooth group algebra of $ G $. The underlying bornological vector space is of course
again the space of smooth functions with compact support on $ G $. Multiplication is given by the convolution product
\begin{equation*}
(f * g)(t) = \int_G f(s) g(s^{-1}t) ds
\end{equation*}
which turns $ \D(G) $ into a bornological algebra. Note that $ \D(G) $ does not have a unit unless $ G $ is discrete. 
The corresponding multiplier algebra is determined in \cite{Meyersmoothrep}. 
\begin{prop}\label{multHG} Let $ G $ be a Lie group. The multiplier algebra of 
the smooth group algebra $ \D(G) $ is the algebra $ \mathcal{E}'(G) $ of distributions on $ G $ with compact support. 
\end{prop}
Explicitly, a left multiplier $ L $ of $ \D(G) $ defines a distribution $ D_L $ on $ G $ by the formula
$$
D_L(f) = L(f)(e).
$$
Remark that the complex group ring $ \mathbb{C}G $ is contained in 
$ M(\D(G)) = \mathcal{E}'(G) $ as the subalgebra spanned by the 
Dirac distributions $ \delta_s $ for $ s \in G $. \\
Using proposition \ref{multHG} one may describe the comultiplication $ \Delta: \D(G) \rightarrow \mathcal{E}'(G \times G) $ by 
\begin{equation*}
\Delta(f)(h) = \int_G f(s) h(s,s) ds.
\end{equation*}
The counit $ \epsilon: \D(G) \rightarrow \mathbb{C} $ is defined by 
\begin{equation*}
\epsilon(f) = \int_G f(s) ds. 
\end{equation*}
Finally, the antipode $ S: \D(G) \rightarrow \D(G) $ is given by 
\begin{equation*}
S(f)(t) = \delta(t) f(t^{-1})
\end{equation*}
and we have again $ S^2 = \id $. 
The general theory developped in the previous sections yields immediately the following result. 
\begin{prop} Let $ G $ be a Lie group. Then the smooth group algebra 
$ \D(G) $ of $ G $ is a bornological quantum group.
\end{prop}
A left and right invariant integral $ \phi $ for $ \D(G) $ is given by evaluation at the 
identity,
\begin{equation*}
\phi(f) = \delta_e(f) = f(e).
\end{equation*}
We remark that $ \phi $ is not a trace. More precisely, we have 
$$
\phi(f * g) = \int_G f(t) g(t^{-1}) dt = \int \delta(t) g(t) f(t^{-1}) dt = \phi((\delta \cdot g) * f)
$$
for all $ f, g \in \D(G) $. \\
As mentioned above, one may as well consider smooth functions on arbitrary locally compact groups $ G $ 
and obtain corresponding bornological quantum groups $ C^\infty_c(G) $ and $ \D(G) $. The definition of 
the space of smooth functions in this setting involves the structure theory of locally 
compact groups \cite{MZ} in order to reduce to the case of Lie groups. More information can be 
found in \cite{Meyersmoothrep} where smooth representations of locally compact groups on bornological 
vector spaces are studied. \\
Actually, it is immediate from the definitions that a smooth representation of the group $ G $ is the same 
thing as an essential comodule over $ C^\infty_c(G) $. In \cite{Meyersmoothrep} it is shown that the category 
of smooth representations of $ G $ is naturally isomorphic to the category of essential modules over $ \D(G) $. 
This statement may be viewed as a special case of theorem \ref{modcomod} and explains the motivation for the 
general definitions of essential modules and comodules given in section \ref{secmodcomod}. \\ 
In the context of locally compact groups it is more natural to work with continuous functions than to 
consider smooth functions. Eventually, this leads to the study of 
quantum groups in the setting of $ C^* $-algebras. The most satisfactory definition of such quantum groups 
is due to Kustermans and Vaes \cite{KV}. Although their definition resembles the definition of a bornological quantum 
group to some extent, it has to be emphasized that the theory of locally compact quantum groups is technically much 
more involved. For instance, basic examples show that the counit and the antipode of a quantum group do not exist 
on the level of $ C^* $-algebras in general. 

\section{Schwartz algebras and discrete groups} \label{secschwartz}

In this section we describe bornological quantum groups arising from Schwartz algebras of certain Lie groups 
as well as from algebras of functions satisfying various decay conditions on finitely generated discrete groups. \\
We begin with the abelian Lie group $ G = \mathbb{R}^n $. Let $ \S(\mathbb{R}^n) $ be the Schwartz space of rapidly 
decreasing smooth functions on $ \mathbb{R}^n $. The topology of this nuclear Fr\'echet space is defined by 
the seminorms 
$$
p^k_{\alpha}(f) = \sup_{x \in \mathbb{R}^n} \bigg{|} \frac{\partial^\alpha f(x)}{\partial x_\alpha}(1 + |x|)^k \biggr{|} 
$$
for any multiindex $ \alpha $ and any nonnegative integer $ k $ where 
$$
|x| = \sqrt{x_1^2 + \cdots + x_n^2}.
$$ 
We write $ \hat{\S}(\mathbb{R}^n) $ 
for the essential bornological algebra obtained by equipping $ \S(\mathbb{R}^n) $ with the pointwise multiplication of functions. 
In order to identify the corresponding multiplier algebra recall that a function $ f \in C^\infty(\mathbb{R}^n) $ is called slowly increasing if for every 
multiindex $ \alpha $ there exists an integer $ k $ such that 
$$
\sup_{x \in \mathbb{R}^n} \biggl{|}\frac{1}{(1 + |x|)^k} \frac{\partial^\alpha f(x)}{\partial x_\alpha} \biggr{|} < \infty. 
$$
Slowly increasing functions on $ \mathbb{R}^n $ form an algebra under pointwise multiplication. 
\begin{lemma}\label{multsrn}
The multiplier algebra $ M(\hat{\S}(\mathbb{R}^n)) $ is the algebra of slowly increasing functions on $ \mathbb{R}^n $. 
\end{lemma}
\proof It is evident that every multiplier of $ \hat{\S}(\mathbb{R}^n) $ is given by a smooth function $ f $ on $ \mathbb{R}^n $. 
Multiplication by such a function induces a continuous linear map from $ \S(\mathbb{R}^n) $ 
to itself iff $ f $ is slowly increasing \cite{Treves}. \qed \\
To define the quantum group structure of $ \hat{\S}(\mathbb{R}^n) $ the formulas for $ C^\infty_c(\mathbb{R}^n) $ 
carry over. The comultiplication is an essential homomorphism and the associated 
Galois maps are isomorphisms. Moreover the Lebesgue integral defines a faithful Haar integral for 
$ \hat{\S}(\mathbb{R}^n) $. Hence we obtain the following result. 
\begin{prop}\label{schwartzrn}
The algebra $ \hat{\S}(\mathbb{R}^n) $ of rapidly decreasing functions on $ \mathbb{R}^n $ is a bornological quantum group. 
\end{prop}
Let us also describe the dual of $ \hat{\S}(\mathbb{R}^n) $. We will denote 
this quantum group by $ \S(\mathbb{R}^n) $ and call it the tempered group algebra of $ \mathbb{R}^n $. 
The underlying algebra structure is given by $ \S(\mathbb{R}^n) $ with convolution multiplication. 
In order to determine the multiplier algebra of $ \S(\mathbb{R}^n) $ let us denote by $ B(\mathbb{R}^n) $ 
the space of all smooth functions $ f $ on $ \mathbb{R}^n $ such that 
all derivatives of $ f $ are bounded. The topology on $ B(\mathbb{R}^n) $ is given by 
uniform convergence of all derivatives. By definition, a bounded distribution is a continuous linear form on the space 
$ B(\mathbb{R}^n) $. A distribution $ T \in \D'(\mathbb{R}^n) $ has rapid decay if the distribution $ T (1 + |x|)^k $ is bounded 
for all $ k \geq 0 $. 
\begin{lemma}\label{multtempgroup}
The multiplier algebra $ M(\S(\mathbb{R}^n)) $ of $ \S(\mathbb{R}^n) $ is the algebra of distributions with rapid decay. 
\end{lemma}
\proof The Fourier transform defines an algebra isomorphism $ \S(\mathbb{R}^n) \cong \hat{\S}(\mathbb{R}^n) $. 
In particular, the multiplier algebras of $ \S(\mathbb{R}^n) $ and $ \hat{\S}(\mathbb{R}^n) $ are isomorphic. 
According to lemma \ref{multsrn} the multiplier algebra of $ \hat{\S}(\mathbb{R}^n) $ is the algebra of slowly increasing 
functions. Under Fourier transform, slowly increasing functions correspond to rapidly decreasing distributions \cite{Treves}. \qed \\
The comultiplication, counit, antipode and the Haar integral for $ \S(\mathbb{R}^n) $ can be determined in the same way as for the 
smooth group algebra $ \D(\mathbb{R}^n) $. Remark that the classical 
Fourier transform can be viewed as an isomorphism of bornological quantum groups $ \hat{\S}(\mathbb{R}^n) \cong \S(\hat{\mathbb{R}}^n) $ 
where $ \hat{\mathbb{R}}^n $ is the dual group of $ \mathbb{R}^n $. \\
The tempered group algebra $ \S(\mathbb{R}) $ and its dual as well as corresponding 
crossed products have been considered by Elliot, Natsume and Nest in their work on the cyclic cohomology of one-parameter crossed
products \cite{ENN}. \\
Let us also explain how the abelian case treated above can be extended to nilpotent Lie groups. 
The algebra $ \S(G) $ of Schwartz functions on a nilpotent Lie group $ G $ has been considered by Natsume and Nest \cite{NN} 
in connection with their study of the cyclic cohomology of the Heisenberg group. \\
Let $ G $ be an $ n $-dimensional connected and 
simply connected Lie group with Lie algebra $ \mathfrak{g} $. Fix a Jordan-H\"older sequence 
$$ 
0 = \mathfrak{g}_0 \subset \mathfrak{g}_1 \subset \cdots \subset \mathfrak{g}_n = \mathfrak{g} 
$$
for the Lie algebra of $ G $ and a basis $ X_1, \dots, X_n $ such that $ X_j \in \mathfrak{g}_j\setminus \mathfrak{g}_{j - 1} $. 
Then one has a diffeomorphism $ \phi: \mathbb{R}^n \rightarrow G $ given by 
$$
\phi(t_1, \dots, t_n) = \exp(t_1 X_1) \cdots \exp(t_n X_n)
$$
and this diffeomorphism may be used to define the Schwartz space $ \S(G) $ as the space of smooth functions on $ G $ corresponding to the 
Schwartz space of $ \mathbb{R}^n $. The space $ \S(G) $ is independent of the choice of Jordan-H\"older basis. 
We denote by $ \hat{\S}(G) $ the bornological algebra of Schwartz functions on $ G $ with pointwise multiplication and 
write $ \S(G) $ for the bornological algebra obtained by considering the convolution product. 
As a generalization of proposition \ref{schwartzrn} and the previous discussion we obtain the following statement. 
\begin{prop}
Let $ G $ be a connected simply connected nilpotent Lie group. Then $ \hat{\S}(G) $ and $ \S(G) $ define a dual 
pair of bornological quantum groups in a natural way. 
\end{prop}
Now let $ \Gamma $ be a finitely generated discrete group equipped with a word metric. We denote by $ L $ the associated length function 
on $ \Gamma $. The function $ L $ satisfies 
$$
L(e) = 0, \qquad  L(t) = L(t^{-1}) \qquad \text{and}\qquad  L(st) \leq L(s) + L(t) 
$$
for all $ s,t \in \Gamma $. Following the notation in \cite{Meyercomb} we define several function spaces associated to $ \Gamma $. 
For every $ k \in \mathbb{R} $ consider the norm 
$$
|| f ||^k = \sum_{t \in \Gamma} |f(t)|(1 + L(t))^k 
$$
on the complex group ring $ \mathbb{C}\Gamma $ and denote by $ \S^k(\Gamma) $ the corresponding Banach space completion. 
We write also $ l^1(\Gamma) $ instead of $ \S^0(\Gamma) $. 
Moreover let $ \S(\Gamma) $ be the completion of $ \mathbb{C}\Gamma $ with respect to the family of 
norms $ ||-||^k $ for all $ k \in \mathbb{N} $. The natural map $ \S^{k + 1}(\Gamma) \rightarrow S^k(\Gamma) $ is compact for 
all $ k \in \mathbb{N} $ and hence $ \S(\Gamma) $ is a Fr\'echet Schwartz space. In particular, the bounded and the precompact 
bornology on $ \S(\Gamma) $ agree. We call $ \S(\Gamma) $ the space of Schwartz functions on $ \Gamma $. 
Remark that for $ \Gamma = \mathbb{Z}^n $ with its natural length function we reobtain the usual definition of the 
space $ \S(\mathbb{Z}^n) $ of Schwartz functions. \\
Consider moreover the norm
$$
|| f ||_\alpha = \sum_{t \in \Gamma} |f(t)| \alpha^{L(t)}
$$
for $ \alpha > 1 $. We write $ l^1(\Gamma, \alpha) $ for the completion of $ \mathbb{C}\Gamma $ with respect 
to this norm and $ \mathcal{O}(\Gamma) $ for the completion with respect to the 
family $ ||-||_n $ for $ n \in \mathbb{N} $. The space $ \mathcal{O}(\Gamma) $ is again a Fr\'echet Schwartz space.
Moreover let $ \S^\omega(\Gamma) $ be the direct limit of the Banach spaces
$ l^1(\Gamma, \alpha) $ for $ \alpha > 1 $. This space is a Silva space, that is, a bornological vector 
space which is the direct limit of a sequence of Banach spaces with injective and compact structure maps. 
All these function spaces do not depend on the choice of the word metric. \\
It is easy to check that the norms $ ||-||^k $ and $ ||-||_\alpha $ are submultiplicative with respect to 
the convolution product. As a consequence, all function spaces considered above become bornological algebras 
in a natural way. Moreover the antipode, the counit and the Haar functional of $ \mathbb{C}\Gamma $ extend 
continuously to the completions. For every $ f \in \mathbb{C}\Gamma $ we have the estimates 
\begin{align*}
||-||^k \cotimes_\pi ||-||^l &(\Delta(f)) \leq \sum_{t \in \Gamma} |f(t)|\, || \delta_t ||^k || \delta_t ||^l \\
&= \sum_{t \in \Gamma} |f(t)| (1 + L(t))^k (1 + L(t))^l = || f ||^{k + l} 
\end{align*}
of the projective tensor product norms where $ \delta_t $ denotes the characteristic function located in $ t $. Similary,
$$
||-||_\alpha \cotimes_\pi ||-||_\beta (\Delta(f)) \leq \sum_{t \in \Gamma} |f(t)|\, || \delta_t ||_\alpha|| \delta_t ||_\beta = 
|| f ||_{\alpha \beta} 
$$
for all $ \alpha, \beta > 1 $. We obtain the following statement. 
\begin{prop}
Let $ \Gamma $ be a finitely generated discrete group. Then the algebras $ l^1(\Gamma), \S(\Gamma), \mathcal{O}(\Gamma) $ and $ \S^\omega(\Gamma) $ 
are bornological quantum groups in a natural way. 
\end{prop}
The algebra structure of the corresponding dual quantum groups is obtained by equipping the above spaces of functions with
pointwise multiplication. 

\section{Rieffel deformation}\label{secrieffel}

In the monograph \cite{Rieffelmemoirs} Rieffel studies deformation quantization for Poisson brackets 
arising from actions of $ \mathbb{R}^d $. A basic example of such a deformation is the 
Moyal product for functions on $ \mathbb{R}^{2n} $. Although the main focus in \cite{Rieffelmemoirs} 
is on the study of the $ C^* $-algebras arising in this way, a large part of the theory 
is carried out in the setting of Fr\'echet spaces. \\
If the underlying manifold is a Lie group one may restrict attention to 
those deformations which are compatible with the group structure. As it turns out, one obtains quantum groups 
in the setting of $ C^* $-algebras in this way \cite{Rieffelcqg}, \cite{Rieffelncqg}. 
We shall only consider the case of compact Lie groups. A remarkable feature of the corresponding compact 
quantum groups is that they arise from deformations of the algebra of all smooth functions and not only of the algebra of 
representative functions. The deformed algebras of smooth functions fit naturally into the framework 
of bornological quantum groups. In this section we shall discuss this point, however, our exposition will be brief since all the 
necessary work is already done in the papers by Rieffel. \\
Let $ G $ be a compact Lie group and let $ T $ be an $ n $-dimensional torus in $ G $ with Lie algebra $ \mathfrak{t} $. 
We identify $ \mathfrak{t} $ with $ \mathbb{R}^n  $ and set $ V = \mathbb{R}^n \times \mathbb{R}^n $. 
Let $ \exp: \mathfrak{t} \rightarrow T $ denote the exponential map. 
Moreover let $ J $ be a skew-symmetric operator on $ V $ with respect to the standard inner product. 
In order to obtain a Poisson bracket which is compatible with the group structure of $ G $ 
we shall assume that the operator $ J $ is of the form 
$ J = K \oplus (-K) $ where $ K $ is a skew-symmetric operator on $ \mathbb{R}^n $. \\
Using this data, the deformed product of $ f, g \in C^\infty(G) $ is defined by
$$
(f \star_K g)(x) = \int f(\exp(-Ks)x \exp(-Ku))\, g(\exp(-t)x \exp(v))\, e^{2 \pi i (\bra s, t \ket + \bra u,v\ket)}
$$ 
where the variabes of integration range over $ \mathbb{R}^n $. This formula yields a continuous and associative 
multiplication on $ C^\infty(G) $ equipped with its natural Fr\'echet topology. We write $ C^\infty(G)_K $ for the 
corresponding bornological algebra. \\
The algebra $ C^\infty(G)_K $ together with the ordinary comultiplication, antipode, counit and Haar integral of 
$ C^\infty(G) $ given by 
$$ 
\Delta(f)(s,t) = f(st), \qquad S(f)(t) = f(t^{-1}), \qquad  \epsilon(f) = f(e), \qquad \phi(f) = \int f(t) dt 
$$
becomes a bornological quantum group. In particular, the classical Haar functional 
is also faithful with respect to the deformed multiplication. \\
We have thus the following statement. 
\begin{prop}
Let $ G $ be a compact Lie group and let $ T $ be a torus in $ G $ with Lie algebra $ \mathfrak{t} $. 
For every skew-symmetric matrix $ K $ on $ \mathfrak{t} $ there exists a bornological quantum group 
$ C^\infty(G)_K $ with structure as described above. 
\end{prop}
We refer to \cite{Rieffelmemoirs}, \cite{Rieffelcqg} for more information and concrete examples. 

\bibliographystyle{plain}

\end{document}